\PassOptionsToPackage{dvipsnames}{xcolor}
\documentclass[12pt]{amsart}
\usepackage{amscd}     
\usepackage{amssymb}
\usepackage{amsmath, amsthm, graphics}
\usepackage{graphicx}
\usepackage{xypic}     
\LaTeXdiagrams        
\usepackage[all]{xy}
\xyoption{2cell} \UseAllTwocells \xyoption{frame} \CompileMatrices
\allowdisplaybreaks[3]
\usepackage{amsfonts}
\usepackage{tikz}
\usepackage{tikz-cd}

\usepackage{mathtools, mathrsfs}
\usepackage{verbatim} 

\usepackage{booktabs}

\usepackage[margin=4cm]{geometry}

\usepackage{latexsym}
\usepackage{epsfig}
\usepackage{amsfonts}
\usepackage{enumerate}
\usepackage{times,float}
\usepackage{graphics}

\usepackage{xcolor}

\usepackage[nodisplayskipstretch]{setspace}
\usepackage[colorlinks=true, linkcolor=blue, citecolor=blue, pagebackref=true]{hyperref}
\newtheorem{theorem}{Theorem}[section]
\newtheorem{corollary}[theorem]{Corollary}
\newtheorem{lemma}[theorem]{Lemma}
\newtheorem{example}[theorem]{Example}

\newtheorem{proposition}[theorem]{Proposition}

\theoremstyle{remark}

\theoremstyle{remark}

\newtheorem{remark}[subsection]{Remark}

\numberwithin{equation}{section}

\newcommand{\Mbar}{\overline{\M}}

\newcommand{\M}{\mathcal{M}}

\newcommand{\F}{\mathcal{F}}

\newcommand{\T}{\mathsf{T}}
\newcommand{\sfK}{\mathsf{K}}

\newcommand{\blang}{\big\langle}
\newcommand{\brang}{\big\rangle}

\newcommand{\CC}{\mathbb{C}}

\def\End{{\rm End}}
\def\b1{{\mathbf 1}}
\def\Hom{{\rm Hom}}

\def\id{{\text{id}}}
\def\GL{{\text{GL}}}
\def\Tr{{\text{Tr}}}
\def\loc{{\text{loc}}}

\newcommand{\Res}{\operatorname*{Res}}

\title{Elliptic Gromov-Witten invariants of instanton moduli}

\author[Hsian-Hua Tseng]{Hsian-Hua Tseng}
\address{Department of Mathematics\\ Ohio State University\\ 100 Math Tower, 231 West 18th Ave.\\ Columbus, OH 43210\\ USA}
\email{hhtseng@math.ohio-state.edu}

\author[Xiaohan Yan]{Xiaohan Yan}
\address{Institut de math\'ematiques de Toulouse\\ 118 route de Narbonne\\ 31400 Toulouse\\ France}
\email{xiaohan.yan@math.univ-toulouse.fr}

\date{\today}



\begin{document}
\raggedbottom


\begin{abstract}
We calculate the genus $1$ Gromov-Witten invariants with one tautological divisor insertion of the moduli space $\M(r,n)$ of rank-$r$, $c_2=n$ torsion-free sheaves on $\mathbb{P}^2$ framed along the infinity divisor. As rational functions of the framing parameters, we determine the polar parts of these invariants by exploring the equivariant geometry of the Gieseker-Uhlenbeck map $\M(r,n)\to\M_0(r,n)$. We determine the polynomial part of these invariants by relating their large framing limits to the results \cite{ILPT} for the rank $1$ case, namely $\text{Hilb}^n(\mathbb{C}^2)$. The degree-$0$ invariants with one divisor insertion are evaluated by studying the combinatorics of multi-partitions.

We also prove reconstruction results for invariants with multiple insertions. Under the assumption that the eigenvalues of the divisor operator are linearly independent as functions of the Novikov variable, these results determine the full genus $1$ theory from the one-divisor series and the genus $0$ divisor operator. 

\end{abstract}

\maketitle

\setcounter{tocdepth}{1}
\tableofcontents

\setcounter{section}{-1}

{\large\part{Introduction}}

There is a deep connection between enumerative geometry and geometric representation theory, see \cite{Ok} for an expert's account. Gromov-Witten theory of Nakajima varieties is a prominent example of this connection that has seen important developments: see \cite{op}, \cite{MOb} for results in genus $0$ of Hilbert schemes of points and \cite{MOk} for results in genus $0$ for more general Nakajima varieties. 

The instanton moduli spaces $\M(r,n)$ form a natural class of targets in which to study this deep connection. They parametrize framed torsion-free sheaves of rank $r$ and second Chern class $n$ on $\mathbb{P}^2$, and include the Hilbert schemes $\M(1,n)=\text{Hilb}^n(\mathbb C^2)$ as their rank-one special case. Their genus-zero equivariant Gromov–Witten theory has a particularly rich algebraic structure: the work of Maulik–Okounkov \cite[Part II]{MOk} relates quantum multiplication to Yangian actions and gives explicit formulas for quantum divisor operators. 

It is natural to try to further deepen the connection between enumerative geometry and geometric representation theory by studying {\em higher genus} Gromov-Witten theory of Nakajima varieties. In genus $1$, the case of Hilbert schemes of points have been worked out in \cite{ILPT} and \cite{ILPTY}. The goal of the present paper is to study genus $1$ Gromov-Witten theory of the next interesting class of Nakajima varieties, the moduli spaces of instantons $\M(r,n)$.

\section{Instanton moduli}
We begin with describing the target geometry.

\subsection{Framed sheaves}
Let 
$$\M(r,n)$$
be the moduli space of torsion-free sheaves on $\mathbb{P}^2$ of rank $r$ and $c_2=n$ framed along the divisor at infinity. More precisely, let
\[
 \ell_\infty=\{[z_0=0:z_1 :z_2]\}\subset\mathbb P^2, \qquad \mathbb C^2=\mathbb P^2\setminus\ell_\infty.
\]
The variety $\M(r,n)$ parametrizes pairs $(\mathcal E,\Phi)$, where $\mathcal E$ is a rank-$r$ torsion-free sheaf on $\mathbb P^2$, locally free near $\ell_\infty$, with\footnote{The framing implies $c_1(\mathcal E)=0$.}
\[
 c_1(\mathcal E)=0,\qquad c_2(\mathcal E)=n,
\]
and
\[
 \Phi:\mathcal E|_{\ell_\infty}
 \xrightarrow{\ \sim\ }\mathcal O_{\ell_\infty}^{\oplus r}
\]
is a framing.

In Physics, $\M(r,n)$ is a moduli space of framed Yang-Mills instantons on $\mathbb{R}^4$. Following this point of view, the number $n$ in $\M(r,n)$ will be called {\em charge}.

In the case rank $r=1$, since $\mathcal{E}$ is torsion free and $c_1(\mathcal{E})=0$, it follows that $\mathcal{E}\hookrightarrow\mathcal{E}^{\vee\vee}\simeq \mathcal{O}_{\mathbb{P}^2}$. It can be deduced from this that $\M(1,n)\simeq \text{Hilb}^n(\mathbb{C}^2)$. 

The tangent space of $\M(r,n)$ at $(\mathcal{E},\Phi)$ is
\[
T_{(\mathcal E,\Phi)}\M(r,n)
=\operatorname{Ext}^1_{\mathbb P^2}\bigl(\mathcal E,\mathcal E(-\ell_\infty)\bigr).
\]
The corresponding Ext groups in degrees $0$ and $2$ vanish:
\[
 \operatorname{Hom}_{\mathbb{P}^2}(\mathcal E,\mathcal E(-\ell_\infty))=0,
 \qquad
 \operatorname{Ext}_{\mathbb{P}^2}^2(\mathcal E,\mathcal E(-\ell_\infty))=0.
\]
The second equality follows from Serre duality and the vanishing of
$\operatorname{Hom}(\mathcal E,\mathcal E(-2\ell_\infty))$. Therefore, $\M(r,n)$ is a smooth variety whose dimension, calculated by Riemann-Roch, is 
\[
 \dim_{\mathbb C}\M(r,n)=2rn.
\]
\subsection{ADHM construction}\label{sec:ADHM_construction}
We recall the ADHM description of $\M(r,n)$, see \cite{ADHM} and \cite[Chapter 2]{NakB}. Let $V=\mathbb C^n$, $W=\mathbb C^r$, and consider 
\[
 \mathbf M_{r,n} =\End(V)^{\oplus2}\oplus\Hom(W,V)\oplus\Hom(V,W).
\]
Consider the group $$G=\GL(V)$$ action on $\mathbf{M}_{r,n}$ given by
\[
 g\cdot(B_1,B_2,I,J)=(gB_1g^{-1},gB_2g^{-1},gI,Jg^{-1}).
\]
The complex moment map for this action is given by 
\[
 \mu: \mathbf{M}_{r,n}\to \End(V), \quad \mu(B_1,B_2,I,J)=[B_1,B_2]+IJ.
\]
An ADHM monad $(B_1,B_2, I,J)$ is stable if no proper subspace $S\subsetneq V$ contains $\text{Im}(I)$ and is invariant under both $B_1$ and $B_2$. Denote by $\mu^{-1}(0)^{\mathrm{st}}\subset \mu^{-1}(0)$ the locus of stable objects. The ADHM construction states that
\[
 \M(r,n)\simeq\mu^{-1}(0)^{\mathrm{st}}/G=:\mu^{-1}(0)\mathbin{/\!\!/}_{\mathrm{st}}G,
\]
where the second expression denotes the GIT quotient in the determinant chamber. It is known that stability is equivalent to the property that $V$ is generated from $\text{Im}(I)$ by repeated applications of $B_1$, $B_2$. The $G$-action on the stable locus is free.

We recall the identification between ADHM monad and framed sheaves. Consider homogeneous coordinates $[z_0:z_1:z_2]\in \mathbb{P}^2$. For $(B_1, B_2, I,J)\in\mathbf{M}_{r,n}$, define
\[
 \alpha=
 \begin{pmatrix}
 z_0B_1-z_1\id_V\\
 z_0B_2-z_2\id_V\\
 z_0J
 \end{pmatrix},
 \qquad
 \beta=
 \begin{pmatrix}
 -z_0B_2+z_2\id_V&
 z_0B_1-z_1\id_V&
 z_0I
 \end{pmatrix}.
\]
We have $\beta\alpha=z_0^2\bigl([B_1,B_2]+IJ\bigr)$. If $(B_1,B_2, I,J)\in \mu^{-1}(0)$, then $\beta\alpha=0$, i.e. we get a complex
\[
 V\otimes\mathcal O_{\mathbb P^2}(-1)
 \xrightarrow{\ \alpha\ }
 (V\oplus V\oplus W)\otimes\mathcal O_{\mathbb P^2}
 \xrightarrow{\ \beta\ }
 V\otimes\mathcal O_{\mathbb P^2}(1).
\]
Define $\mathcal{E}$ to be the middle cohomology of this complex.  On $\ell_\infty$, the two $V$-summands form the exact Koszul complex and the remaining cohomology is
$W\otimes\mathcal O_{\ell_\infty}$. This gives the framing. Conversely, every framed sheaf above arises from a unique stable ADHM orbit.

The ADHM construction shows that $\M(r,n)$ is irreducible and quasi-projective. Furthermore, the form
\begin{equation}
 \omega=\operatorname{tr}  \bigl(dB_1\wedge dB_2+dI\wedge dJ\bigr)
\end{equation}
descends to a holomorphic symplectic form on $\M(r,n)$.

\subsection{Tautological bundle}
The principal $G$-bundle on the quotient $\mu^{-1}(0)^{\text{st}}/G$ gives the rank-$n$ tautological bundle $\mathcal V$ over $\M(r,n)$. In the sheaf description, if $\mathscr{E}$ is the universal framed sheaf and $p_M:\mathbb P^2\times \M(r,n)\to \M(r,n)$ is the projection, then
\[
 \mathcal V=R^1(p_M)_*\bigl(\mathscr E(-\ell_\infty)\bigr).
\]
The tautological divisor
\[
 H=c_1(\det\mathcal V)
\]
plays a crucial role in our study. In particular, when $r,n\geq1$ and $(r,n)\neq (1,1)$, we have\footnote{We exclude the following cases. $\M(1,1)=\mathbb A^2$. Also, we set $\M(r,0)=\mathrm{pt}$ In both cases, there are no nonzero curve classes and every Gromov-Witten invariant of positive degree is defined to be $0$.} 
$$\text{Pic}(\M(r,n))=\mathbb{Z}H, \quad H_2(\M(r,n),\mathbb{Z})/\text{tor}=\mathbb{Z}\beta,\quad \text{and }\int_\beta H=1.$$  

\subsection{Gieseker-Uhlenbeck map}
Forgetting the stability condition, we have the quotient
\[
 \M_0(r,n)=\mu^{-1}(0)/G.
\]
This is the Uhlenbeck compactification of spaces of instantons. Variation of GIT quotients gives the following projective birational morphism, called the Gieseker-Uhlenbeck morphism, 
\begin{equation}\label{eqn:GU_map}
 \pi_n:\M(r,n)\longrightarrow \M_0(r,n).
\end{equation}
It can be seen that $\pi_n$ is induced from the inclusion $\mu^{-1}(0)^{\text{st}}\subset\mu^{-1}(0)$, and $\pi_n$ is a symplectic resolution. In the rank $r=1$ case, the Gieseker-Uhlenbeck map $\pi_n$ is the Hilbert-Chow morphism $\M(1,n)=\text{Hilb}^n(\mathbb{C}^2)\to\text{Sym}^n(\mathbb{C}^2)=\M_0(1,n)$.

Set-theoretically, the Uhlenbeck space has the following stratification
\[
 \M_0(r,n)
=\coprod_{k=0}^n\M_0^{\mathrm{reg}}(r,n-k)\times\operatorname{Sym}^k(\mathbb C^2).
\]
Here $\M_0^{\mathrm{reg}}(r,n-k)$ parametrizes framed locally free sheaves. On closed points, $\pi_n$ sends $(\mathcal E,\Phi)$ to
\[
 \left((\mathcal E^{\vee\vee},\Phi),
 \operatorname{Supp}(\mathcal E^{\vee\vee}/\mathcal E)\right),
\]
where the support is counted with multiplicity.   

\subsection{Torus action}
Consider the following torus,
\[
 \T=(\mathbb C^*)^2_{\mathrm{coord}}\times(\mathbb C^*)^r_{\mathrm{framing}}.
\]
We denote the multiplicative and additive characters of $\T$-action as follows
\begin{equation}\label{eqn:parameters}
 t_i=e^{\epsilon_i},\,\, i=1,2;\qquad
 u_\alpha=e^{a_\alpha},\,\, \alpha=1,...r.
\end{equation} 
We call $a_1,...,a_r$ the {\em framing parameters} and call charaters of $(\mathbb C^*)^r_{\mathrm{framing}}$ the {\em framing characters}. We also call indices of the framing parameters {\em (framing) colors}, see Remark \ref{rmk:color terminology} below for further explanations. 

Throughout the paper, we use the following notations: 
\begin{equation}\label{eqn:parameter_expressions}
\hbar=t_1t_2=e^h,
 \qquad h=\epsilon_1+\epsilon_2,
 \qquad p=\epsilon_1\epsilon_2.
\end{equation}

We consider the $\T$-action on $\M(r,n)$ given by the $(\mathbb C^*)^2_{\mathrm{coord}}$-action on $\mathbb{P}^2$: 
\[
 [z_0:z_1:z_2]\longmapsto[z_0:t_1z_1:t_2z_2]
\]
together with the diagonal $(\mathbb C^*)^r_{\mathrm{framing}}$-action that changes the framing. In the ADHM description, the $\T$-action is given by 
\[
 (t_1,t_2;\mathbf u)\cdot(B_1,B_2,I,J)
 =\bigl(t_1B_1,t_2B_2,I\mathbf u^{-1},\hbar\,\mathbf uJ\bigr), 
\quad (t_1,t_2; \mathbf{u})\in (\mathbb C^*)^2_{\mathrm{coord}}
   \times(\mathbb C^*)^r_{\mathrm{framing}}.
\]
Here we view $\mathbf u=\operatorname{diag}(u_1,\ldots,u_r)\in\operatorname{GL}(W)$.
The $\T$-action preserves stability and takes $\mu$ to $\hbar\mu$.  

The diagonal framing torus 
$$\mathbb{C}_{\text{diag}}^*=\{\text{diag}(u,...,u)\ |\ u\in \mathbb{C}^*\} \subset (\mathbb{C}^*)_{\text{framing}}^r,$$
acts trivially on $\M(r,n)$. Tangent characters therefore depend on the framing parameters only through
\[
 a_{\alpha\beta}=a_\alpha-a_\beta.
\]
The tautological divisor $H$ admits a natural $\T$-linearization. The $\T$-character of the symplectic form $\omega$ is $\hbar$. Therefore the tangent and cotangent bundles satisfy $T_{\M(r,n)}=\hbar\,T_{\M(r,n)}^\vee$ in equivariant $K$-theory.

\subsection{Fixed points}\label{sec:fixed-pts}
For a partition $\lambda=(\lambda_1\ge\lambda_2\ge\cdots)$, the diagram of $\lambda$ may be described as
\[
 \mathsf{D}_\lambda=\{(i,j) \ |\ i\ge1,\ 1\le j\le\lambda_i\}.
\]
The monomials $x^{i-1}y^{j-1}\in \mathbb{C}[x,y], (i,j)\in \mathsf{D}_\lambda$ form a basis of the quotient $\mathbb{C}[x,y]/I_\lambda$ by the corresponding monomial ideal $I_\lambda=(x^{i-1}y^{\lambda_i}, i\geq 1)$.

The $\T$-fixed points of $\M(r,n)$ are isolated and are indexed by $r$-tuples
of partitions
\begin{equation}\label{eqn:multi-partition}
 \boldsymbol\lambda=(\lambda^1,\ldots,\lambda^r),
 \qquad |\boldsymbol\lambda|=\sum_{\alpha=1}^r|\lambda^\alpha|=n.
\end{equation}
The number of $\T$-fixed points of $\M(r,n)$ is $p_r(n)$, i.e. the number of $r$-colored partitions of $n$. The framed sheaf indexed by $\boldsymbol{\lambda}$ is a direct sum of ideal sheaves:
\[
 \mathcal E_{\boldsymbol\lambda}=\bigoplus_{\alpha=1}^r I_{\lambda^\alpha},
\]
where the summand is assigned the $\alpha$-th framing character $u_\alpha$. In ADHM description, we have $J=0$, the image of $I$ generates the corner box of each nonempty diagram, and $B_1,B_2$ move through the coordinate directions.

At the fixed point $\boldsymbol{\lambda}$, the fiber of the tautological bundle $\mathcal{V}$ has $\T$-character
\[
 \mathsf V_{\boldsymbol\lambda}
 :=\left.\mathcal V\right|_{\boldsymbol\lambda}
 =\sum_{\alpha=1}^r u_\alpha\sum_{(i,j)\in\lambda^\alpha}t_1^{1-i}t_2^{1-j},
\]
and $\T$-weights 
\[
 a_\alpha-(i-1)\epsilon_1-(j-1)\epsilon_2, \qquad (i,j)\in\lambda^\alpha.
\]
The fiber of the tautological divisor $H$ has $\T$-character  
\[
 \left.H\right|_{\boldsymbol\lambda}
 =\sum_{\alpha=1}^r\left(|\lambda^\alpha|a_\alpha
 -\epsilon_1n(\lambda^\alpha)-\epsilon_2n((\lambda^\alpha)')\right),
 \quad n(\lambda)=\sum_{i\ge1}(i-1)\lambda_i.
\]

The fixed loci of the effective framing torus
\begin{equation}\label{eqn:eff_frame_torus}
 \mathsf{A}=(\mathbb C^*)_{\text{framing}}^r/\mathbb C^*_{\mathrm{diag}}
\end{equation}
are given by
\begin{equation}\label{eqn:framed_torus_fixed_loci}
 \M(r,n)^{\mathsf{A}}
=\coprod_{N_1+\cdots+N_r=n}\prod_{\alpha=1}^r\text{Hilb}^{N_\alpha}(\mathbb C^2).
\end{equation}
\begin{remark}\label{rmk:color terminology}
We view $\boldsymbol\lambda$ as an \(r\)-colored partition: every box of \(\lambda^\alpha\) has color \(\alpha\). We call an expression \emph{same-color} if it involves a single index \(\alpha\), and \emph{cross-color} if it involves distinct indices \(\alpha\neq\beta\).

The integer $N_\alpha=|\lambda^\alpha|$ is the {\em charge}\footnote{This choice of terminology is consistent with calling $n$ in $\M(r,n)$ the charge: indeed, in (\ref{eqn:framed_torus_fixed_loci}) the color $\alpha$ indexes the factor $\text{Hilb}^{N_\alpha}(\mathbb{C}^2)=\M(1, N_\alpha)$.} carried by color \(\alpha\), and $\mathbf N=(N_1,\ldots,N_r), \sum_\alpha N_\alpha=n,$ is the \emph{color-charge vector}. A color-charge block consists of all partition tuples with these fixed component sizes.
\end{remark}

\subsection{Tangent weights}\label{sec:tangent-weights}
We describe the $\T$-weights of the tangent spaces $T_{\M(r,n)}|_{\boldsymbol{\lambda}}$. Consider the deformation complex of a stable ADHM monad:
\[
 \End(V)\longrightarrow
 \begin{matrix}
 t_1\End(V)\oplus t_2\End(V)\oplus\Hom(W,V)
 \oplus\hbar\Hom(V,W)
 \end{matrix}
 \longrightarrow\hbar\End(V).
\]
Its maps are
\[
 \xi\longmapsto \bigl([\xi,B_1],[\xi,B_2],\xi I,-J\xi\bigr),
 \quad 
 (b_1,b_2,i,j)\longmapsto [b_1,B_2]+[B_1,b_2]+iJ+Ij.
\]
Stability makes the first map injective and the second surjective. The middle cohomology is the desired tangent space. Hence its $\T$-character is
\[
\left.T_{\boldsymbol\lambda}\M(r,n):=T_{\M(r,n)}\right|_{\boldsymbol\lambda}
 =\mathsf W^\vee\mathsf V_{\boldsymbol\lambda}
 +\hbar\mathsf V_{\boldsymbol\lambda}^\vee\mathsf W
 -(1-t_1)(1-t_2)\mathsf V_{\boldsymbol\lambda}^\vee\mathsf V_{\boldsymbol\lambda}, 
  \quad \mathsf W:=\sum_{\alpha=1}^r u_\alpha.
\]

For $s=(i,j)$, define the arm and relative leg by
\[
 A_\lambda(s)=\lambda_i-j, \qquad L_\mu(s)=\mu'_j-i.
\]
The relative leg may be negative. Expanding the preceding character gives the relative arm-leg formula
\[
\begin{split}
 T_{\M(r,n)}\big|_{\boldsymbol\lambda}
 =\sum_{\alpha,\beta=1}^r\frac{u_\beta}{u_\alpha}
 \Bigg[
 \sum_{s\in\lambda^\alpha}t_1^{-L_{\lambda^\beta}(s)}t_2^{A_{\lambda^\alpha}(s)+1}
 +\sum_{s\in\lambda^\beta}t_1^{L_{\lambda^\alpha}(s)+1}t_2^{-A_{\lambda^\beta}(s)}
 \Bigg].
\end{split}
\]
Equivalently, the additive tangent weights of $T_{\M(r,n)}|_{\boldsymbol{\lambda}}$ are
\begin{equation}\label{eqn:full_tangent_weights}
\begin{aligned}
 &a_\beta-a_\alpha
 -L_{\lambda^\beta}(s)\epsilon_1
 +(A_{\lambda^\alpha}(s)+1)\epsilon_2,
 &&s\in\lambda^\alpha,\\
 &a_\beta-a_\alpha
 +(L_{\lambda^\alpha}(s)+1)\epsilon_1
 -A_{\lambda^\beta}(s)\epsilon_2,
 &&s\in\lambda^\beta,
\end{aligned}
\end{equation}
where $(\alpha,\beta)$ runs through all ordered pairs. There are $2rn$ weights. The first weight for $(\alpha,\beta,s)$ and the second weight for $(\beta,\alpha,s)$ add to $h$, as required by the symplectic form. For $\alpha=\beta$, the two weights reduce to
\begin{equation}\label{eqn:tangent_weights}
 -L(s)\epsilon_1+(A(s)+1)\epsilon_2,
 \qquad
 (L(s)+1)\epsilon_1-A(s)\epsilon_2.
\end{equation}

\section{Results}
The main result is the evaluation of the genus $1$ Gromov-Witten invariants of $\M(r,n)$ with one tautological divisor insertion, in all degrees.

\subsection{Gromov-Witten theory}
We are interested in Gromov-Witten theory of $\M(r,n)$. Let 
\begin{equation*}
\Mbar_{g,k}(\M(r,n),d)    
\end{equation*}
be the moduli stack of stable maps $f:(C, p_1,...,p_k)\to \M(r,n)$ where $C$ is a nodal curve of arithmetic genus $g$, $p_1,...,p_k\in C$ are smooth points, and the morphism $f$ is stable and $f_*[C]=d\beta$. The moduli stack $\Mbar_{g,k}(\M(r,n),d)$ admits a perfect obstruction theory. It is not proper, but the fixed loci of the $\T$-action on $\Mbar_{g,k}(\M(r,n),d)$ induced from the $\T$-action of $\M(r,n)$ is proper. Consequently, (primary) $\T$-equivariant genus $g$ Gromov-Witten invariants of $\M(r,n)$ may be defined by virtual localization:
\begin{equation*}
\langle \phi_1,...,\phi_k\rangle_{g,k,d}^{\M(r,n)}:=\int_{[\Mbar_{g,k}(\M(r,n),d)^\T]^{\mathrm{vir}}}\frac{\prod_{i=1}^kev_i^*(\phi_i)}{e_\T(N^{\mathrm{vir}})}\in\mathbb{Q}(\epsilon_1,\epsilon_2, a_1,...,a_r)=:K,
\end{equation*}
where $\phi_1,...,\phi_k\in H^*_\T(\M(r,n))_\loc$ are equivariant cohomology classes and 
\begin{equation*}
ev_i:\Mbar_{g,k}(\M(r,n),d)\to\M(r,n), \quad (f:(C, p_1,...,p_k)\to \M(r,n))\mapsto f|_{p_i},    
\end{equation*}
$i=1,...,k$ are the evaluation maps.

The genus $0$ Gromov-Witten invariants of $\M(r,n)$ are fully determined by the work of \cite[Part II]{MOk}. We summarize how this works. Consider the operator 
$$H\star: H_\T^*(\M(r,n))_{\text{loc}}[\![q]\!]\to H_\T^*(\M(r,n))_{\text{loc}}[\![q]\!]$$ 
of small quantum multiplication by the divisor $H$, defined by 
\begin{equation*}
\int_{\M(r,n)}^\T (H\star \phi)\cup \varphi=\sum_{d\geq 0}q^d\langle H, \phi,\varphi\rangle_{0,3,d}^{\M(r,n)}, \quad \phi, \varphi\in H_\T^*(\M(r,n))_\loc.    
\end{equation*}
The operator $H\star$ is explicitly calculated in \cite[Section 15.1]{MOk}. Set $D=\dim_K H_\T^*(\M(r,n))_{\loc}$. By \cite[Theorem 15.2.1]{MOk}, after passing to the field $K(\!(q)\!)$, the quantum cohomology ring is generated by $H$ and has an additive basis $1,H,\ldots,H^{\star(D-1)}$. Since the divisor matrix is rational in $q$, this basis is also defined over $K(q)$. The WDVV reconstruction argument in \cite[Section 4]{op} can then be used to explicitly and effectively determine all genus $0$ generating functions of $\M(r,n)$.  

Introduce the canonical sign modification\footnote{A polarization compatible with the tangent character is
\[
 T^{1/2}=(t_1-1)\End(\mathcal V)+\Hom(W,\mathcal V),
 \qquad T_{\M(r,n)}=T^{1/2}+\hbar(T^{1/2})^\vee.
\]
Its nonequivariant first Chern class satisfies
$\kappa_{\M(r,n)}=c_1(T^{1/2})=rH\pmod2$. This explains the sign choice.}
 following \cite{MOk}:
\begin{equation}\label{eqn:signed_novikov}
 Q=(-1)^rq.
\end{equation}
We write $H\star_Q$ for the operator $H\star$ viewed as a function of $Q$. Thus every series expressed in terms of ordinary geometric invariants and $Q$ uses
\[
 \sum_{d\geq0}q^d\langle\cdots\rangle^{\M(r,n)}_{g,k,d}
 =\sum_{d\geq0}(-1)^{rd}Q^d\langle\cdots\rangle^{\M(r,n)}_{g,k,d}.
\]
This convention applies also to big quantum products and their  derivatives in flat vector fields. 

\subsection{Divisor insertion}
Consider the following generating functions of genus $1$ Gromov-Witten invariants of $\M(r,n)$:
\begin{equation}\label{eqn:g1_gen_function}
 \F^+_{r,n}(Q)
=\sum_{d\ge1}(-1)^{rd}Q^d\left\langle H\right\rangle^{\M(r,n)}_{1,1,d\beta},
\quad \F_{r,n}(Q)=\langle H\rangle_{1,1,0}^{\M(r,n)}+\F^+_{r,n}(Q). 
\end{equation}

For $m\ge1$, set
\[
 \sigma_{-1}(m)=\sum_{\delta\mid m}\frac1\delta,
 \quad
\text{and } x_{\delta,m}=\delta\epsilon _1+\frac m\delta\epsilon _2 \,\,\text{ if } \delta | m. 
\]
Define the coefficients $b_{r,m}$ by
\begin{equation}\label{eqn:coefficients}
 b_{r,0}=1,
\qquad  
b_{r,m}
=r\sigma_{-1}(m)
+\frac{2p}{h}\sum_{1\le\alpha<\beta\le r}\ 
\sum_{\delta\mid m}\frac{x_{\delta,m}}{x_{\delta,m}^{\,2}-a_{\alpha\beta}^{\,2}},
 \quad m>0.
\end{equation}
Finally, let
\begin{equation}\label{eqn:traces}
\tau_{r,N}(Q)=\text{Tr}_{H_\T^\bullet(\M(r,N))_{\text{loc}}}\bigl(H\star_Q-H\cup\bigr),
 \qquad \tau_{r,0}=0.
\end{equation}
The scalar normalization in this trace is essential in rank one.

\begin{theorem}[Genus-one, divisor insertion]\label{thm:main}
For every $r\ge1$, $n\ge1$, we have
\[
\F^+_{r,n}(Q)=-\frac h{24p}\sum_{m=0}^{n-1}b_{r,m}\tau_{r,n-m}(Q).
\]
\end{theorem}
Together with the evaluation of degree $0$ invariants in Theorem \ref{thm:closed-main}, this determines all invariants $\langle H\rangle_{1,1,d}^{\M(r,n)}$.

\subsection{Strategy of Proof}
The proof of Theorem \ref{thm:main} is based on the following: for $d>0$, as a rational function of $a_1,...,a_r$, we have 
\begin{equation}
\langle H\rangle_{1,1,d}^{\M(r,n)}=(\text{polar part})+(\text{polynomial part}).     
\end{equation}
We identify the polar part and polynomial part explicitly.

\subsubsection{Polar part}
We consider $\langle H\rangle_{1,1,d}^{\M(r,n)}$ as a rational function of the framing parameters over $\mathbb Q(\epsilon_1,\epsilon_2)$. Consider $\T$-weights of the form:
\begin{equation}\label{eqn:linear_factors}
w:=a_\alpha-a_\beta-i\epsilon_1-j\epsilon_2 , 
\quad 1\leq \alpha\neq \beta\leq r,\,\, i,j\geq 1.   
\end{equation}
The weight $w$ corresponds to the $\T$-character 
\begin{equation}\label{eqn:char_w}
 \chi_w:\T\longrightarrow\mathbb C^*,\quad
 (\mathbb C^*)^2_{\rm coord}\times(\mathbb C^*)^r_{\rm framing}\ni (t_1,t_2;u_1,\ldots,u_r)\longmapsto\frac{u_\alpha}{u_\beta}\,t_1^{-i}t_2^{-j}.
\end{equation}
Define the following connected codimension-$1$ subtorus associated to $\chi_w$:
\begin{equation}\label{eqn:torus_w}
 \sfK_w:=\ker(\chi_w)
 =\left\{(t_1,t_2;u_1,\ldots,u_r): \frac{u_\alpha}{u_\beta}=t_1^it_2^j\right\}\subset \T.
\end{equation}
Clearly $\mathbb{C}^*_{\text{diag}}\subset \sfK_w$, so we obtain the same codimension-$1$ subtorus after passing to the effective quotient. There is an exact sequence of tori
\[
 1\longrightarrow \sfK_w\longrightarrow \T
 \xrightarrow{\,\chi_w\,}\mathbb C^*\longrightarrow 1,
 \qquad \T/\sfK_w\cong\mathbb C^*.
\]
On character lattices, this becomes
\[
 0\longrightarrow\mathbb Z w\longrightarrow \text{Char}(\T)
 \xrightarrow{\,\mathrm{res}\,}\text{Char}(\sfK_w)
 \longrightarrow 0.
\]
Hence a $\T$-weight $\lambda$ restricts to the trivial character of $\sfK_w$ if
and only if $\lambda$ is an integral multiple of $w$. In additive notation, we write this as $\lambda|_{\sfK_w}=0$.

We process $\langle H\rangle_{1,1,d}^{\M(r,n)}$ by factoring the $\T$-equivariant pushforward to a point through the Gieseker-Uhlenbeck map $\pi_n$ in (\ref{eqn:GU_map}) and, for each $w$, perform localization with respect to $\T$ in two stages: first localize with respect to $\sfK_w$, then with respect to $\T/\sfK_w$. The analysis of $\sfK_w$-fixed loci of $\M_0(r,n)$ shows that the expression $w$ appears in the denominator exactly when localizing with respect to the quotient torus $\T/\sfK_w$. Further analysis shows the following
\begin{theorem}[Polar part]\label{prop:polar_part} 
The polar part of $\langle H\rangle_{1,1,d}^{\M(r,n)}$ is of the form
\begin{equation}\label{eqn:polar_part}
\sum_w \frac{\text{Res}_w}{w}    
\end{equation}
where the sum runs through 
\begin{equation}\label{eqn:poles}
w=a_\alpha-a_\beta-i\epsilon_1-j\epsilon_2, 
\quad 1\leq \alpha\neq \beta\leq r,\,\, i,j\geq 1, \,\, ij\leq n.   
\end{equation}
Furthermore, we have 
\begin{equation}\label{eqn:residue}
\text{Res}_w=\frac{1}{24}(-1)^{rd}[Q^d]\tau_{r,n-ij}(Q).    
\end{equation}
\end{theorem}
There are three nontrivial statements here: the first is that only $w$ of the form (\ref{eqn:poles}) occur; the second is that $w$ appears in the denominator with exponent $1$ (``simple poles''); the third is the condition $ij\leq n$. To identify $\text{Res}_w$ for each possible $w$, we work with a description of the map $\pi_n$ over an \'etale neighborhood of a $\sfK_w$-fixed point in $\M_0(r,n)$ due to \cite{BC}.

\subsubsection{Polynomial part}
Let $t\in \mathbb{C}$ and define $\{b_\alpha\}$ by
\begin{equation}\label{eqn:parameter_scaling}
a_\alpha=tb_\alpha, \quad \alpha=1,...,r.    
\end{equation}
Our approach to studying the polynomial part is to make the substitution (\ref{eqn:parameter_scaling}) and consider the limit $t\to \infty$ for generic values of $\{b_\alpha\}$. For this purpose, it is more convenient to set $$s=t^{-1}$$ and consider the $s=0$ specialization.\footnote{Intuitively, the $s=0$ specialization corresponds to the loci (\ref{eqn:framed_torus_fixed_loci}) fixed by the action of the effective framing torus $\mathsf{A}$ given by (\ref{eqn:eff_frame_torus}).} We call this process the {\em large framing limit}.

To study the $s=0$ specialization of $\F_{r,n}^+(Q)$, we use Givental's formula \cite{g1} for genus $1$ Gromov-Witten generating functions:
\begin{equation}\label{eqn:givental_formula}
\F_{r,n}(Q)= \frac1{48}Q\frac{d}{dQ}\left(\log\prod_i\Delta_i\right)
  +\frac12\sum_i(\mathsf{R}_1)_{ii}\rho_i(Q),
\end{equation}
see Section \ref{sec:limit_genus1} for details. To study the $s=0$ specialization of (\ref{eqn:givental_formula}), we first identify the $s=0$ specialization of the Frobenius structure on $QH_\T^*(\M(r,n))_{\text{loc}}$. This uses the construction of stable envelopes. See Section \ref{sec:Frob} for details. The Frobenius structure is semisimple. We want to show that the canonical idempotents become canonical idempotents in the $s=0$ specialization. To do this, we use the relations of the divisor operator $H\star_Q$ with Baxter algebra. See Section \ref{sec:projector} for details. We also need to identify the $s=0$ specialization of the $\mathsf{R}$-matrix coefficients occurring in (\ref{eqn:givental_formula}). This is done by directly studying the recursive construction of these $\mathsf{R}$-matrix coefficients. See Section \ref{sec:limit_genus1} for details. 

The outcome of this analysis shows that the $s=0$ specialization of $\F_{r,n}^+(Q)$ is a generating function of genus $1$ Gromov-Witten invariants of $\M(r,n)^{\mathsf{A}}$. We explicitly evaluate this function in Theorem \ref{thm:coulomb-limit} using the results of \cite{ILPT}, the answer is
\begin{equation*}
 \lim_{t\to\infty}\F^+_{r,n}(tb;Q)
 =-\frac h{24p}\left(\tau_{r,n}(Q)
 +r\sum_{m=1}^{n-1}\sigma_{-1}(m)\tau_{r,n-m}(Q)\right).
\end{equation*}

Since the $s=0$ specialization (i.e. the $t\to \infty$ limit) of $\F_{r,n}^+(q)$ exists, it follows that the part of $\F_{r,n}^+(Q)$ that is polynomial in $\{a_\alpha\}$ must be constant in $\{a_\alpha\}$, and the constant is given by the $s=0$ specialization already determined. The detailed argument is presented in Section \ref{sec:pf_main_thm}.

\begin{remark}
Theorem \ref{thm:main} includes the $\text{Hilb}^n(\mathbb{C}^2)$ result of \cite{ILPT} as the rank $r=1$ case. However, as outlined above, our proof {\em uses} the result in rank $r=1$ and does not give another proof of the result of \cite{ILPT}. The rank $r=1$ result of \cite{ILPT} is proved by appealing to the relation with Gromov-Witten theory of families of local curves. We do not know how to extend this approach to the rank $r>1$ case.    
\end{remark}

\subsection{Degree zero invariants}
We compute the genus $1$ degree $0$ $\T$-equivariant invariants $\langle H \rangle_{1,1,0}^{\mathcal{M}(r,n)}$ separately. Theorem \ref{thm:closed-main} gives a closed formula for their generating series. 

By Getzler's formula \cite[Proposition 2.7]{Getzler} and fixed point localization, we reduce the computation to a sum over $\T$-fixed points in $\mathcal{M}(r,n)$, indexed by $r$-colored partitions, of rational expressions in the equivariant parameters of $\T$. The single-color terms in the sum come essentially from the $r=1$ case computed in \cite{ILPT}. The ingredients genuinely new to the $r>1$ case are certain cross-color terms, namely terms involving two distinct directions in the framing space, which we compute in Lemma \ref{lemma: norm Phi a} and Lemma \ref{lem:pointed-two-partition}.

For Lemma \ref{lemma: norm Phi a}, the key observation is that in spite of its cross-color nature, a K-theoretic model of what we need to compute appears already in the $r=1$ case as a correspondence between a pair of fixed points. It is studied in \cite{co} via an arm-and-leg formula. As we adapt the computation from K-theoretic characters to cohomological weights, further cancellation occurs between dual characters, giving rise to a relatively simple form of the result.

For Lemma \ref{lem:pointed-two-partition}, we take a more direct approach. Similar to the $d>0$ case, what we need to compute is a rational function which splits into a part with simple poles and a regular part. For the polar part, we find the residues by identifying them with weighted counts of certain intersection configurations of boundaries of translated Young diagrams. Certain cancellation is then observed based on transformations of these boundaries. In the end, we determine the regular part using the limiting behavior of the target rational function at infinity. 

\subsection{Reconstruction}
For classes $\phi_1,...,\phi_k\in H^*_\T(\M(r,n))_\loc$, consider the generating functions of genus $1$ Gromov-Witten invariants,
\begin{equation}\label{eqn:multi-pt-series}
\langle \phi_1,...,\phi_k\rangle_{g}^{\M(r,n)}=\sum_{d\geq 0}q^d\langle \phi_1,...,\phi_k\rangle_{g,k,d}^{\M(r,n)}\in\mathbb{Q}(\epsilon_1,\epsilon_2, a_1,...,a_r)[\![q]\!].
\end{equation}
We have the following reconstruction result\footnote{Here we use the variable $q$.} of (\ref{eqn:multi-pt-series}).

\begin{theorem}\label{Red1}
To every genus $1$ series $\blang H^{\star k_1},...,H^{\star k_\ell}\brang^{\M(r,n)}_{1}$, there are canonically associated functions
\begin{equation}\label{ccck}
\big\{ C_{k,m} \big\}_{0\leq k \leq \text{dim}\, H_\T^*(\M(r,n))_\loc-1\, ,\,  0\leq m\leq \ell-1} \subset  \mathbb{Q}(\epsilon_1, \epsilon_2, a_1,...,a_r,q)    
\end{equation}
such that the following equation holds:
\begin{equation*}
\blang H^{\star k_1},...,H^{\star k_\ell}\rangle_{1}^{\M(r,n)}=\sum_{k=0}^{\text{dim}\, H_\T^*(\M(r,n))_\loc-1}\, \sum_{m=0}^{\ell-1} \, C_{k,m}\cdot \left(q\frac{d}{dq}\right)^m\blang H^{\star k}\brang_{1}^{\M(r,n)}\, .    
\end{equation*}
\end{theorem}

The functions \eqref{ccck} are constructed as rational functions of the matrix coefficients of $H\star$. The proof of Theorem \ref{Red1}, which is the same as the proof of \cite[Theorem 4]{ILPT}, uses the WDVV relations in genus $0$ and Getzler's relation in genus $1$ \cite{Getzler}. 

The arguments of \cite[Sections 4.4-4.5]{ILPT} apply to prove the following

\begin{theorem}[Conditional reconstruction]\label{Red3}
Let $\theta=q\,d/dq$, and let $\rho_1,\ldots,\rho_D$ be the eigenvalues of $H\star$ in a differential splitting field $L$ over $K(q)$. Assume that their Wronskian is nonzero:
\begin{equation}\label{eqn:reconstruction-wronskian}
 \det\bigl(\theta^{j-1}\rho_i\bigr)_{1\leq j,i\leq D}\ne0.
\end{equation}
Equivalently, these functions are linearly independent over the constant field $L^\theta=\ker\theta$. Then all primary genus $1$ generating functions of $\M(r,n)$ can be effectively reconstructed from $\blang H\brang_1^{\M(r,n)}$ and the genus $0$ theory determined by $H\star$.
\end{theorem}

\begin{proof}
We spell out the role of the hypothesis in the argument of \cite[Sections 4.4-4.5]{ILPT}. In a normalized canonical frame, set $r_i=(\mathsf R_1)_{ii}$. Givental's formula expresses the known one-divisor series, after subtracting its known genus $0$ metric term $\frac12\sum_i\rho_i r_i$. The first $\mathsf R$-matrix recursion determines the off-diagonal entries from genus $0$ data, and the diagonal recursion then, determines $\theta r_i$ and all its further $\theta$-derivatives from that data. Applying $1,\theta,\ldots,\theta^{D-1}$ to the one-divisor identity and moving these known derivatives to the right gives a linear system for the $r_i$, with coefficient matrix
$(\theta^{j-1}\rho_i)_{j,i}$. Hypothesis \eqref{eqn:reconstruction-wronskian} makes this system invertible. The genus-$1$ one-form consequently determines every one-point primary series. Theorem~\ref{Red1} reconstructs the marked primary series.
\end{proof}

(\ref{eqn:reconstruction-wronskian}) will be discussed elsewhere \cite{ILPST}.

\subsection{Acknowledgment}
In the process of developing the results presented here and writing this paper, the authors benefited from the assistance of GPT-5.6 Sol and GPT-6 Astra. The work of XY was supported by the ERC starting grant 101164820.

{\large\part{The polar part}}
The key tool in our study of the polar part is the equivariant geometry of the Gieseker-Uhlenbeck map (\ref{eqn:GU_map}). 

\section{Fixed loci of Uhlenbeck space}
The Uhlenbeck space $\M_0(r,n)$ is affine. We use the ADHM description to describe its coordinate ring and study $\sfK_w$-fixed loci. 

\subsection{Description of the coordinate ring}\label{sec:coord_ring}
Decompose the framing space into its $1$-dimensional weight spaces:
\[
 W=\bigoplus_{\gamma=1}^r W_\gamma,
 \qquad
 W_\gamma=\mathbb C e_\gamma.
\]
Accordingly, write
\[
 I=\sum_{\gamma=1}^r I_\gamma, \,\,  I_\gamma\in\Hom(W_\gamma,V); 
 \quad 
 J=(J_1,\ldots,J_r),\,\, J_\gamma\in\Hom(V,W_\gamma), \quad \text{so } IJ=\sum_{\gamma=1}^r I_\gamma J_\gamma.
\]

The quiver relevant to the affine ADHM quotient is the frozen framed doubled Jordan quiver. It has one unfrozen vertex \(0\) carrying the \(n\)-dimensional vector space \(V\) and \(r\) frozen vertices \(\infty_1,\ldots,\infty_r\) with \(\infty_\gamma\) carrying the $1$-dimensional vector space \(W_\gamma\). There are two loops
\begin{equation*}
\xymatrix{
0\ar@(ul,dl)_{b_1}\ar@(ur,dr)^{b_2}
}    
\end{equation*}
and, for every \(\gamma=1,\ldots,r\), arrows
\begin{equation*}
\xymatrix{
\infty_\gamma\ar@/^0.7pc/[r]^{i_\gamma} & 0.\ar@/^0.7pc/[l]^{j_\gamma}
}    
\end{equation*}
In the ADHM description, these arrows are represented by
\[
 b_1\longmapsto B_1,\qquad
 b_2\longmapsto B_2,\qquad
 i_\gamma\longmapsto I_\gamma,\qquad
 j_\gamma\longmapsto J_\gamma.
\]
Only the vertex \(0\) is gauged with gauge group $G=\GL(V)$, and the vertices \(\infty_\gamma\) are frozen: we do not quotient by \(\GL(W_\gamma)\)'s.

The relation at the gauge vertex is
\begin{equation}\label{eqn:quiver_relation}
 [B_1,B_2]+\sum_{\gamma=1}^r I_\gamma J_\gamma=0.
\end{equation}
Consider the coordinate ring of the space of ADHM monads:
\[
R=\mathcal O(\mathbf M_{r,n})
=\operatorname{Sym}\left(
      \left(\End(V)^{\oplus2}
       \oplus\bigoplus_{\gamma=1}^r\Hom(W_\gamma,V)
       \oplus\bigoplus_{\gamma=1}^r\Hom(V,W_\gamma)\right)^\vee
     \right),
\]
and let \(\mathfrak a\subset R\) be the ideal generated by the matrix entries of the left-hand side of (\ref{eqn:quiver_relation}). Then
\[
 \mathcal O\bigl(\mu^{-1}(0)\bigr)=R/\mathfrak a,
 \qquad
 \mathcal O\bigl(\M_0(r,n)\bigr)=(R/\mathfrak a)^G.
\]

We now describe generators of the coordinate ring of $\M_0(r,n)$. For a word
\begin{equation}\label{eqn:word}
 \omega=\ell_1\cdots\ell_s, \quad \ell_k\in\{1,2\},
\end{equation}
put
\begin{equation}\label{eqn:B_of_word}
 B_\omega=B_{\ell_1}\cdots B_{\ell_s}, \quad B_\emptyset=\id_V.
\end{equation}
There are two types of invariant functions.
\begin{enumerate}
\item For every nonempty word \(\omega\), there is a trace generator
\[
 \tau_\omega=\Tr(B_\omega).
\]
Cyclic permutations of \(\omega\) do not change the trace, and \(\Tr(B_\emptyset)=n\) is constant.

\item For every (possibly empty) word \(\omega\), and every ordered pair \((\gamma,\delta)\), there is a framed-path generator
\[
 p_{\gamma\delta;\omega}:=J_\gamma B_\omega I_\delta.
\]
This is the scalar corresponding to the path
\begin{equation*}
\begin{tikzcd}
\infty_\delta \arrow[dr, "i_\delta"] & \\
 &0 \arrow[dl,"j_\gamma"]\arrow[loop, out=20,  in=340, distance=1em, "b_1"]
   \arrow[loop, out=25,  in=335, distance=2.5em, "b_2"]
   \arrow[loop, out=30, in=330, distance=4em, draw=none, no head, "\dots" {shift={(-0.05,0.05)}}]
   \arrow[loop, out=35,  in=325, distance=5.5em, "b_{\ell_s}"]\\
\infty_\gamma. &
\end{tikzcd} 
\end{equation*}
Thus the path starts at the frozen vertex \(\infty_\delta\) and ends at the frozen vertex \(\infty_\gamma\). For the empty word, we have 
\[
 p_{\gamma\delta;\emptyset}=J_\gamma I_\delta, \,\, \text{corresponding to the path } \infty_\delta\xrightarrow{i_\delta}0\xrightarrow{j_\gamma}\infty_\gamma.
\]    
\end{enumerate}

The frozen-vertex Le Bruyn-Procesi theorem gives
\begin{equation}\label{eqn:invariant_ring_description}
 R^G=
 \mathbb C\left[\{\Tr(B_\omega),\,
  J_\gamma B_\omega I_\delta\ |\ \omega\text{ is a word in }B_1,B_2, 1\leq\gamma,\delta\leq r\} 
 \right]
\end{equation}
More precisely, this is the specialization of \cite[Theorem~2.2]{CY} in which the vertex \(0\) is unfrozen and the vertices \(\infty_\gamma\) are frozen. No other types of generator are needed. Indeed, because every weight space \(W_\eta\) is $1$-dimensional, a path that passes through a frozen vertex in its interior factors into a product of shorter framed paths. For example,
\[
 J_\gamma B_U I_\eta J_\eta B_V I_\delta
 =\bigl(J_\gamma B_U I_\eta\bigr)\bigl(J_\eta B_V I_\delta\bigr).
\]
Likewise, a cycle that passes through a frozen vertex reduces to a framed path. For example,
\[
\Tr\bigl(B_U I_\gamma J_\gamma B_V\bigr)=J_\gamma B_VB_UI_\gamma.
\]
The collection of generators in (\ref{eqn:invariant_ring_description}) is generally redundant, but it is sufficient for the fixed-locus calculation.

Since \(G=\GL(V)\) is linearly reductive over \(\mathbb C\), taking \(G\)-invariants is an exact functor. Taking $G$-invariants of
\[
 0\longrightarrow\mathfrak a
 \longrightarrow R
 \longrightarrow R/\mathfrak a
 \longrightarrow0
\]
yields a surjective map
\[
 R^G\longrightarrow(R/\mathfrak a)^G.
\]
Consequently, imposing the ADHM relation introduces relations among the generators in (\ref{eqn:invariant_ring_description}), but introduces no new types of generators. This is also the content of \cite[Proposition~3.2]{CY}.

One relation will be used repeatedly. For arbitrary words \(U,V\), the ADHM quiver relation (\ref{eqn:quiver_relation}) gives
\begin{equation}\label{eqn:relation1}
 \begin{aligned}
 J_\gamma B_U(B_1B_2-B_2B_1)B_VI_\delta
 =-J_\gamma B_UIJB_VI_\delta
 =-\sum_{\eta=1}^r \bigl(J_\gamma B_UI_\eta\bigr)\bigl(J_\eta B_VI_\delta\bigr).
 \end{aligned}
\end{equation}

\subsection{$\sfK_w$-fixed loci}\label{sec:K_w-fixed}
Fix a weight $w=a_\alpha-a_\beta-i\epsilon_1-j\epsilon_2$ as in (\ref{eqn:linear_factors}). Recall the definition of $\chi_w$ in (\ref{eqn:char_w}) and \(\sfK_w=\ker(\chi_w)\). Our goal here is to prove the following scheme-theoretic description of the $\sfK_w$-fixed locus:
\begin{equation}\label{eqn:K_w-fixed loci final}
 \M_0(r,n)^{\sfK_w}\cong
 \begin{cases}
 \mathbb A^1_z,&ij\leq n,\\[3pt]
 \text{a finite scheme supported at }z=0,&ij>n.
 \end{cases}
\end{equation}

Consider the coordinate ring of the Uhlenbeck space $\M_0(r,n)$ (which is affine),
\[
 A_{r,n}=\mathcal O\bigl(\M_0(r,n)\bigr).
\]
The \(\sfK_w\)-action decomposes \(A_{r,n}\) into character spaces:
\[
 A_{r,n}=\bigoplus_{\eta\in \text{Char}(\sfK_w)}A_\eta.
\]
The scheme-theoretic fixed locus is
\begin{equation}\label{eqn:K_w-fixed_loci1}
 \M_0(r,n)^{\sfK_w}=
\operatorname{Spec}\left(A_{r,n}\big/\left\langle A_\eta:\eta\neq0\right\rangle\right).
\end{equation}
Indeed, if \(f\in A_\eta\) with \(\eta\neq0\), choose \(k\in \sfK_w\) such that \(\eta(k)\neq1\). On the equalizer of the action of \(k\) and the identity, one has $f=\eta(k)f$, and hence \(f=0\). Conversely, after all nontrivial character spaces are eliminated, the \(\sfK_w\)-action is trivial. Thus (\ref{eqn:K_w-fixed_loci1}) represents the fixed-point functor scheme-theoretically.

It is therefore enough to determine the \(\sfK_w\)-weights of the homogeneous generators in (\ref{eqn:invariant_ring_description}). Although a product of generators of opposite nonzero weights may have total weight zero, it vanishes in the quotient (\ref{eqn:K_w-fixed_loci1}), because each of its nonzero weight factors already belongs to the fixed-locus ideal.

We use the convention that a regular function \(f\) has multiplicative \(\T\)-weight \(\chi_f\) if
\[
 f(t\cdot x)=\chi_f(t)f(x).
\]
$\T$ acts on ADHM monads as 
\[
 B_1\longmapsto t_1B_1,\qquad
 B_2\longmapsto t_2B_2,\qquad
 I_\delta\longmapsto u_\delta^{-1}I_\delta,\qquad
 J_\gamma\longmapsto t_1t_2u_\gamma J_\gamma.
\]
Suppose that the word \(\omega\) in (\ref{eqn:word}) contains \(p\) copies of \(1\) and \(q\) copies of \(2\). We say that the {\em bidegree} of $\omega$ is
\[
 \deg(\omega)=(p,q).
\]
The corresponding trace generator $\Tr(B_\omega)$, whose bidegree is defined to be that of $\omega$, has weight
\[
 \chi_{\tau_\omega}=t_1^pt_2^q,
\]
whereas the framed-path generator $J_\gamma B_\omega I_\delta$, whose bidegree is likewise defined to be that of $\omega$, has weight
\begin{equation}\label{eqn:wt_framed-path-generator}
 \chi_{p_{\gamma\delta;\omega}}
 =\frac{u_\gamma}{u_\delta}\,t_1^{p+1}t_2^{q+1}.
\end{equation}
The shift from \((p,q)\) to \((p+1,q+1)\) comes from the factor \(\hbar=t_1t_2\) in the action on \(J\).

Since \(\chi_w\) is primitive, a \(\T\)-character \(\chi\) restricts trivially to \(\sfK_w\) if and only if
\[
 \chi=\chi_w^m\quad \text{ for some } m\in\mathbb Z. 
\]
We now identify $\sfK_w$-invariant generators. 
\begin{enumerate}
\item
If a trace generator of bidegree \((p,q)\) is \(\sfK_w\)-invariant, then
\[
 t_1^pt_2^q=\left(\frac{u_\alpha}{u_\beta}t_1^{-i}t_2^{-j}\right)^m.
\]
Comparing $u_\alpha's$ implies \(m=0\), because \(\alpha\neq\beta\). It then follows that \(p=q=0\). Thus the only \(\sfK_w\)-invariant trace word is the empty trace $\Tr(\id_V)=n$, which is constant. Every nonconstant trace generator vanishes on \(\M_0(r,n)^{\sfK_w}\).

\item
A framed-path generator of bidegree \((p,q)\) is \(\sfK_w\)-invariant precisely when
\begin{equation}\label{eqn:framed_gen_invariance}
 \frac{u_\gamma}{u_\delta}t_1^{p+1}t_2^{q+1}
=\left(\frac{u_\alpha}{u_\beta}t_1^{-i}t_2^{-j}\right)^m, \quad \text{for some }m\in \mathbb{Z}.
\end{equation}

If \(\gamma=\delta\), we again get \(m=0\). The coordinate characters would then require
\[
 p+1=q+1=0,
\]
which is impossible.

If \(\gamma\neq\delta\). In the character lattice of the effective framing torus $\mathsf{A}$ in (\ref{eqn:eff_frame_torus}), we have 
\[
 e_\gamma-e_\delta=m(e_\alpha-e_\beta).
\]
Both nonzero sides are roots of type \(A_{r-1}\), so necessarily \(m=\pm1\). If \(m=1\), then
\[
 (\gamma,\delta)=(\alpha,\beta), \qquad p+1=-i,\quad q+1=-j,
\]
which is impossible because \(p,q\geq0\) and \(i,j\geq1\). Therefore \(m=-1\), and we get 
\[
 \gamma=\beta,\qquad
 \delta=\alpha,\qquad
 p=i-1,\qquad
 q=j-1.
\]
\end{enumerate}
Consequently, the only nonconstant invariant generators that survive on the \(\sfK_w\)-fixed subscheme are
\[
 J_\beta B_{\ell_1}\cdots B_{\ell_{i+j-2}}I_\alpha,
\]
given by a word that contains exactly \(i-1\) copies of \(1\) and \(j-1\) copies of \(2\). Their common \(\T\)-character is
\[
 \frac{u_\beta}{u_\alpha}t_1^it_2^j=\chi_w^{-1},
\]
whose additive weight is \(-w\). The corresponding path starts at \(\infty_\alpha\), enters \(V\) over $0$ through \(I_\alpha\), and ends at \(\infty_\beta\) through \(J_\beta\).

It remains to show that the value on the fixed loci of such a path is independent of the ordering of the 
\(B_1\)'s and \(B_2\)'s. Any two words having the same bidegree are related by adjacent transpositions 
\(B_1B_2\leftrightarrow B_2B_1\). For an adjacent transposition, the relation (\ref{eqn:relation1}) gives
\begin{equation}\label{eqn:relation2}
 \begin{aligned}
 J_\beta B_U(B_1B_2-B_2B_1)B_VI_\alpha
 =-J_\beta B_UIJB_VI_\alpha
 =-\sum_{\eta=1}^r \bigl(J_\beta B_UI_\eta\bigr)\bigl(J_\eta B_VI_\alpha\bigr).
 \end{aligned}
\end{equation}
If the original word has bidegree \((i-1,j-1)\), then after removing the adjacent \(B_1B_2\), the total bidegree of \(U\) and \(V\) is
\[
 \deg(U)+\deg(V)=(i-2,j-2).
\]
In particular, neither factor on the right-hand side of (\ref{eqn:relation2}) can be the unique surviving framed generator of bidegree \((i-1,j-1)\). Each factor therefore has nontrivial \(\sfK_w\)-character and vanishes on the fixed loci. Hence the two words related by the adjacent transposition have the same value on the fixed loci. If \(i=1\) or \(j=1\), there is only one possible ordering and the assertion is trivial.

Thus all words of bidegree \((i-1,j-1)\) determine a single element in the invariant ring
\begin{equation}\label{eqn:invariant_generator}
 z=\left[ J_\beta B_1^{i-1}B_2^{j-1}I_\alpha\right]
 \in
 \mathcal O\bigl(\M_0(r,n)^{\sfK_w}\bigr)=:B^w_{r,n}.
\end{equation}
Here the brackets mean the class in the coordinate ring of the fixed loci, and any ordering of the \(B_1\)'s and \(B_2\)'s can be used in the definition.

Since \(A_{r,n}\) is generated by the homogeneous trace and framed-path generators described above, and all generators except \(z\) and the constants vanish in \(B^w_{r,n}\), the ring \(B^w_{r,n}\) is generated by \(z\). Hence there is a surjection
\begin{equation}\label{eqn:quotient_inv_ring}
 q:\mathbb C[Z]\twoheadrightarrow B_{r,n}^w, \qquad q(Z)=z.
\end{equation}

We now prove the numerical restriction \(ij\leq n\) whenever \(z\) is nonzero, and prove scheme-theoretic finiteness when \(ij>n\).

For $0\leq u<i, 0\leq v<j$, define, in the ADHM description,
\[
 v_{uv}=B_1^uB_2^vI_\alpha(1)\in V,\quad \text{and } \phi_{uv}
 = J_\beta B_1^{i-1-u}B_2^{j-1-v}\in V^\vee.
\]
Form the \(ij\times ij\) matrix $\mathsf{M}$ whose matrix coefficients are the following invariant regular functions:
\[
 \mathsf M_{(u,v),(u',v')}
 =
 \phi_{uv}(v_{u'v'})
 =
 J_\beta B_1^{i-1-u}B_2^{j-1-v} B_1^{u'}B_2^{v'}I_\alpha.
\]
The bidegree of the corresponding framed path is
\[
 \bigl(i-1-u+u',\,j-1-v+v'\bigr),
\]
which is equal to \((i-1,j-1)\) if and only if
\[
 u'=u,\qquad v'=v.
\]
Consequently, 
\begin{equation}\label{eqn:inv_functions}
 \mathsf M_{(u,v),(u',v')}
 =
 \begin{cases}
 z,&(u,v)=(u',v'),\\
 0,&(u,v)\neq(u',v').
 \end{cases}
\end{equation}
Indeed, every off-diagonal entry has bidegree different from \((i-1,j-1)\), and therefore has nontrivial \(\sfK_w\)-character; every diagonal entry has bidegree \((i-1,j-1)\), and the ordering-independence proved above identifies it with \(z\). Thus $\mathsf{M}$ is a multiple of the identity:
\begin{equation}\label{eqn:inv_matrices}
 \mathsf M=zI_{ij} \quad\text{in }\operatorname{Mat}_{ij}(B_{r,n}^w).
\end{equation}

Before passing to the quotient, the matrix \(\mathsf M\) factors through the \(n\)-dimensional vector space \(V\):
\[
 \mathbb C^{ij} \xrightarrow{\ P\ } V \xrightarrow{\ Q\ } \mathbb C^{ij}, 
 \quad \text{where } P(e_{uv})=v_{uv},\,\,(Qv)_{uv}=\phi_{uv}(v).
\]
Therefore
\[
 \text{rank}(\mathsf M)\leq n,
\]
and all minors of size \(n+1\) vanish identically. Since the entries of \(\mathsf M\) are \(G\)-invariant, these determinantal identities descend first to \(\M_0(r,n)\) and then to \(\M_0(r,n)^{\sfK_w}\).

Suppose that \(ij>n\). Choose any \(n+1\) of the pairs \((u,v)\). By (\ref{eqn:inv_matrices}), the corresponding principal minor is $z^{n+1}$. The rank bound therefore gives
\begin{equation}\label{eqn:nilpotent_z}
 z^{n+1}=0 \quad\text{in }B_{r,n}^w.
\end{equation}
Since \(B_{r,n}^w\) is generated by \(z\), it follows that \(B_{r,n}^w\) is a quotient of
\[
 \mathbb C[Z]/(Z^{n+1}).
\]
Thus \(B_{r,n}^w\) is finite-dimensional over \(\mathbb C\), and \(\M_0(r,n)^{\sfK_w}=\text{Spec}\, B_{r,n}^w\) is a finite, possibly nonreduced, scheme supported at \(z=0\).

The same matrix argument shows that if \(z\neq0\) at a geometric point of $\text{Spec}\,B_{r,n}^w$, then 
\(\mathsf M=zI_{ij}\) is invertible. Hence the \(ij\) vectors (the column vectors of $P: \mathbb{C}^{ij}\to V$ above)
\[
 v_{uv}=B_1^uB_2^vI_\alpha(1), \quad 0\leq u<i, 0\leq v<j
\]
are linearly independent, and therefore
\[
 n\geq ij.
\]

In the case $ij\leq n$, we construct the entire $\sfK_w$-fixed axis. Let
\[
 \mathcal R_{i,j} = \mathbb C[x,y]/(x^i,y^j).
\]
Then
\[
\text{dim}\,\mathcal{R}_{i,j}=ij,\quad \text{with a basis } x^uy^v,\quad 0\leq u<i,\,\, 0\leq v<j.
\]
 Put
\[
 V=\mathcal R_{i,j}\oplus V_{\mathrm{sp}},
 \qquad
 \dim V_{\mathrm{sp}}=n-ij.
\]
We refer to the second summand as the spectator charge.

Let \(s\) be the coordinate on \(\mathbb A^1_s\). Define an $\mathbb{A}^1_s$-family of ADHM monads on the above $V$ by
\[
 B_1=m_x\oplus0,
 \qquad
 B_2=m_y\oplus0,
\]
where \(m_x,m_y\) denote multiplications by \(x,y\) on \(\mathcal R_{i,j}\). Define the framing maps $I,J$ by
\[
 I_\alpha(1)=1\in\mathcal R_{i,j},
 \quad
 I_\gamma=0,\quad \gamma\neq\alpha, \text{ and }
 J_\beta=s\,\lambda_{\mathrm{soc}},
 \quad
 J_\gamma=0,\quad \gamma\neq\beta,
\]
where \(J_\gamma\) is $0$ on \(V_{\mathrm{sp}}\), and where
\[
 \lambda_{\mathrm{soc}}:\mathcal R_{i,j}\longrightarrow\mathbb C
\]
is the socle functional determined by
\[
 \lambda_{\mathrm{soc}}(x^{i-1}y^{j-1})=1,
 \qquad
 \lambda_{\mathrm{soc}}(x^uy^v)=0
 \quad\text{for }(u,v)\neq(i-1,j-1).
\]
This family satisfies the ADHM quiver relation (\ref{eqn:quiver_relation}). Indeed, $[B_1,B_2]=0$ because \(m_x\) and \(m_y\) commute. For $IJ$, note that \(I_\gamma\) is nonzero only for \(\gamma=\alpha\), while \(J_\gamma\) is nonzero only for \(\gamma=\beta\), and \(\alpha\neq\beta\). Thus each composition \(I_\gamma J_\gamma\) is zero, and
\[
 \sum_{\gamma=1}^r I_\gamma J_\gamma=0.
\]

The values of the invariant generators on this family are explicit. Every nonempty trace word vanishes because \(B_1\) and \(B_2\) are nilpotent. A framed path gives nonzero value only if its initial color is \(\alpha\) and its terminal color is \(\beta\). In that case,
\[
 J_\beta B_\omega I_\alpha
 =
 \begin{cases}
 s,&\deg(\omega)=(i-1,j-1),\\
 0,&\deg(\omega)\neq(i-1,j-1).
 \end{cases}
\]
Since \(m_x\) and \(m_y\) commute, this value depends only on the bidegree and not on the ordering of the word.

It follows that the family defines a morphism
\[
 \varphi:\mathbb A^1_s\longrightarrow \M_0(r,n)^{\sfK_w}.
\]
On coordinate rings, $\varphi^*:B_{r,n}^w\longrightarrow\mathbb C[s]$ satisfies $\varphi^*(z)=s$. Composing with the surjection (\ref{eqn:quotient_inv_ring}), we obtain
\[
 \mathbb C[Z] \xrightarrow{\, q\, } B_{r,n}^w \xrightarrow{\, \varphi^* \,} \mathbb C[s],
 \qquad
 Z\longmapsto s.
\]
The composition is an isomorphism after the harmless renaming \(Z=s\), and in particular it is injective. Hence \(q\) is injective. Since \(q\) is also surjective, it follows that
\[
 B_{r,n}^w\cong\mathbb C[z].
\]
Therefore
\[
 \M_0(r,n)^{\sfK_w}\cong\mathbb A^1_z,  \quad\text{when }ij\leq n.
\]

We can describe the quiver representations\footnote{$\M(r,n)$ are the Nakajima variety defined using the Jordan quiver and certain dimension and framing vectors. The ADHM construction in Section \ref{sec:ADHM_construction} is the specialization to this case of the construction of Nakajima varieties as spaces of quiver representations.} parametrized by this line. The charge-\(ij\) rectangular block $\mathcal{R}_{i,j}$ is stable for every value of \(s\): the vector $I_\alpha(1)=1$ generates \(\mathcal R_{i,j}\) under multiplication by \(x\) and \(y\). For \(s\neq0\), it is also costable. Indeed, any nonzero
\(\mathbb C[x,y]\)-submodule
\[
 0\neq N\subset\mathcal R_{i,j}
\]
meets the socle nontrivially. Since the socle of $\mathcal{R}_{i,j}$ is $\mathbb C x^{i-1}y^{j-1}$, we have $x^{i-1}y^{j-1}\in N$ up to a nonzero scalar, and
\[
 J_\beta(x^{i-1}y^{j-1})=s\neq0.
\]
Thus no nonzero \(B_1,B_2\)-invariant subspace is contained in \(\ker J\).

When \(n>ij\), the direct sum with \(V_{\mathrm{sp}}\) is not stable in the determinant GIT chamber; this is irrelevant here because the construction takes place on the affine quotient \(\M_0(r,n)\). The spectator summand $V_\text{sp}$ is the direct sum of \(n-ij\) copies of the zero unframed simple. For \(s\neq0\), the rectangular summand \(\mathcal R_{i,j}\) is stable and costable and hence is a
closed framed simple, while the spectator summand \(V_{\mathrm{sp}}\) is semisimple. Their direct sum is therefore a closed representative of the corresponding point in the affine quotient.

We have proved (\ref{eqn:K_w-fixed loci final}). Finally, the weight formula (\ref{eqn:wt_framed-path-generator}) explains why no additional $w$ occur. A surviving framed generator necessarily
has coordinate exponents
\[
 p+1>0,\qquad q+1>0,
\]
and hence can occur only when $a_\alpha-a_\beta=i\epsilon_1+j\epsilon_2, \ i,j\geq1.$ The same calculation excludes the cases $\alpha=\beta$ and $ij<0$.

\section{Proof of Theorem \ref{prop:polar_part}}

We first justify the proper pushforwards used in the two-stage localization.

\begin{lemma}[Properness over the affine quotient]
\label{lem:stable-maps-proper-affine-base}
For fixed $g,k,d$, the morphism
\[
 b:\Mbar_{g,k}(\M(r,n),d)\longrightarrow\M_0(r,n)
\]
obtained by composing with $\pi_n$ is proper. The same assertion holds for the lower-charge instanton moduli with any framing representation, specifically $\M(W^{\text{sl}},m)$ below. Formation of the morphism $b$ commutes with arbitrary base change, namely the following diagram is cartesian for any $B'\to \M_0(r,n)$:
\begin{equation*}
\xymatrix{
\overline{\mathcal{M}}_{g,k}(\M(r,n)\times_{\M_0(r,n)}B',d)\ar[r]\ar[d] & B'\ar[d]\\
\Mbar_{g,k}(\M(r,n),d)\ar[r]^{b} &\M_0(r,n),
}    
\end{equation*}
where the two maps from $\overline{\mathcal{M}}_{g,k}(\M(r,n)\times_{\M_0(r,n)}B',d)$ are formed by composing with the two projections from $\M(r,n)\times_{\M_0(r,n)}B'$ to its factors.
\end{lemma}

\begin{proof}
For a family $\pi_C:C\to S$ of connected nodal curves of arithmetic genus $g$, cohomology and base change gives $(\pi_C)_*\mathcal O_C=\mathcal O_S$, universally under base change. Consequently a morphism from $C$ to the affine scheme $\M_0(r,n)$ factors uniquely through $S$. Thus $\Mbar_{g,k}(\M(r,n),d)$ is precisely the stack of stable maps to the fibers of the projective morphism $\pi_n$. This assertion holds for families, including nonreduced bases.

Choose a relative projective embedding for $\pi_n$. Its degree on $d\beta$ is fixed, since the relative polarization is a positive multiple of $H$. The relative stable-map stack is the closed substack of the fixed-degree stable-map stack of the resulting projective-space bundle whose universal map factors through $\M(r,n)$. The latter stack is proper over $\M_0(r,n)$ because moduli of stable maps to projective targets are proper, so the former is proper as well. Equivalently, the usual stable reduction proof applies to the projective family $\pi_n$. Finally, a map to a fiber product $\M(r,n)\times_{\M_0(r,n)}B'$ is exactly a stable map together with a lift of its base point to $B'$. This proves the base-change assertion; stability and the degree are unchanged. The proof uses only the projectivity of the Gieseker-Uhlenbeck morphism, so it also applies to $\M(W^{\text{sl}},m)$ below.
\end{proof}

Write the invariant as a $\T$-equivariant pushforward to a point:
\begin{equation}
\langle H\rangle_{1,1,d}^{\M(r,n)}=(pt)^\T_*(ev^*H\cap [\Mbar_{1,1}(\M(r,n),d)]^{\mathrm{vir}}).
\end{equation}
Consider the factorization of the map $pt: \Mbar_{1,1}(\M(r,n),d)\to \text{pt}$,
\begin{equation}
\Mbar_{1,1}(\M(r,n),d)\xrightarrow{\, ev\,}\M(r,n)\xrightarrow{\,\pi_n\,}\M_0(r,n)\xrightarrow{\, pt\,} \text{pt}.  \end{equation}

Consider an arbitrary possible framing-dependent denominator. Virtual localization expresses the invariant with denominators that are products of nonzero torus weights. Thus every irreducible finite framing-dependent denominator gives a hyperplane 
$$\{\lambda=0\}\subset \text{Spec} \left(\mathbb{C}(\epsilon_1, \epsilon_2)[a_1,...,a_r]\right)$$ 
for a primitive integral weight $\lambda$ with a nonzero framing part. The exponential $\exp(\lambda):\T\to\mathbb{C}^*$ is a character. Put $\sfK_\lambda=\ker(\exp(\lambda))\subset \T$. Work in the equivariant coefficient ring localized at the prime $(\lambda)$: all characters nontrivial on $\sfK_\lambda$ are units there. The first, $\sfK_\lambda$-localization therefore has no denominator divisible by $\lambda$.

The study of generators of coordinate rings in Sections~\ref{sec:coord_ring}-\ref{sec:K_w-fixed} exclude every hyperplane $\lambda=0$ except those in (\ref{eqn:linear_factors}). Indeed, the fixed subscheme of an affine torus scheme is obtained by setting to $0$ all homogeneous functions of nontrivial weight, rather than by taking the ring of invariant functions. A nonconstant trace generator has no framing part, and so its weight cannot be a nonzero multiple of $\lambda$. A framed path of bidegree $(p,q)$ has weight
\[
 a_\gamma-a_\delta+(p+1)\epsilon_1+(q+1)\epsilon_2,
 \qquad p,q\geq0.
\]
For $\gamma=\delta$ the same observation applies. For distinct colors, its two framing coefficients are $1$ and $-1$, so it is primitive. It can be a multiple of $\lambda$ only if, up to sign, $\lambda$ is one of the characters in (\ref{eqn:linear_factors}), with $i=p+1$ and $j=q+1$. If no nonconstant generator survives, the fixed affine base is the origin scheme-theoretically. Lemma~\ref{lem:stable-maps-proper-affine-base} then makes the $\sfK_\lambda$-fixed stable-map locus proper. Its localized pushforward is regular at $(\lambda)$, so no such pole occurs.

For a remaining $w$, perform localization first with respect to $\sfK_w$ and then with respect to $\T/\sfK_w$. If $ij>n$, (\ref{eqn:K_w-fixed loci final}) says that the fixed affine base is finite. Properness over that base again excludes a pole. If $ij\leq n$, the only possible pole comes from the reduced affine axis in (\ref{eqn:K_w-fixed loci final}). The calculation below proves that this pole is at most simple and computes its coefficient. Characters depending only on $\epsilon_1,\epsilon_2$ belong to the coefficient field and do not define finite poles in framing parameters.

For each $w$ in (\ref{eqn:poles}), we need to identify the numerator, which is denoted in Theorem \ref{prop:polar_part} as $\text{Res}_w$ and is obtained by localization with respect to $\sfK_w$. To evaluate this, we need a description of the map $\pi_n$ over \'etale neighborhoods of suitable points of $\M_0(r,n)$ that is compatible with the $\sfK_w$-action. 

\subsection{Local description of Gieseker-Uhlenbeck map}\label{sec:local_description_Uhlenbeck}

Fix a $w$ from (\ref{eqn:poles}). To simplify notation, put
\[
 m=n-ij.
\]
We work over the generic point of the wall 
$$\{w=0\}\subset \text{Spec}\,\left(\CC[\epsilon_1,\epsilon_2,a_1,...,a_r]\right)$$
in the space of framing parameters. We choose a point \(0\neq z\in\mathbb A^1_z\simeq\M_0(r,n)^{\sfK_w}\). Let
\[
 x_z\in \M_0(r,n)^{\sfK_w}\subset \M_0(r,n)
\]
be the corresponding closed point in the affine quotient. By the discussions in Section \ref{sec:K_w-fixed}, \(x_z\) is represented by the polystable ADHM monad
\begin{equation}\label{eqn:polystable_pt}
 R_z\oplus S^{\oplus m},
\end{equation}
where \(R_z\) is the \(ij\)-dimensional rectangular framed simple
\[
 R_z=\mathbb C[x,y]/(x^i,y^j),\quad \text{ on which }
 B_1=m_x,\,\, B_2=m_y,\,\,
 I_\alpha(1)=1,\,\,
 J_\beta(x^{i-1}y^{j-1})=z,
\]
and where \(S\) is the $1$-dimensional unframed simple supported at the origin: on $S$, we have 
\[
 B_1=B_2=I=J=0.
\]
The summand $R_z$ is stable and costable for \(z\neq0\), and hence is a closed framed simple. Its stabilizer\footnote{When \(m=0\), we understand that \(\GL_0=\{1\}\) and 
\(\M(r,0)=\M_0(r,0)=\operatorname{Spec}\mathbb C\).} is trivial. The stabilizer of (\ref{eqn:polystable_pt}) is therefore $G_{x_z}\cong\GL_m(\mathbb{C})$.

We use the Bellamy-Craw relative slice theorem \cite{BC} to describe $\pi_n$. First, we need to identify the quiver occurring there following \cite[Section 3.2]{BC}: for convenience, we call it the slice quiver. Let \(\mathbf{v}\) be the dimension vector of the unframed simple \(S\). For the Jordan quiver, the symmetric Euler form satisfies
\[
 (\mathbf{v},\mathbf{v})=0.
\]
Equivalently, the diagonal contribution at the gauge vertex is canceled by the contribution of the Jordan loop. Thus
\[
 p(\mathbf{v}):=1-\frac12(\mathbf{v},\mathbf{v})=1.
\]
Consequently, the slice quiver has one vertex corresponding to \(S\) and one loop at that vertex. After doubling, this loop gives the two endomorphisms denoted by \(B'_1,B'_2\).

Let \(\mathbf{v}_R\) denote the dimension vector of the framed simple \(R_z\), including the framing vertex. The framing multiplicity at the vertex of the slice quiver is
\begin{equation}\label{eqn:framing_multiplicity}
 -(\mathbf{v}_R,\mathbf{v})=r.
\end{equation}
Indeed, the gauge-vertex diagonal term is again cancelled by the Jordan loop, while the \(r\) framing arrows contribute \(-r\) to the symmetric Euler pairing. Thus the (transverse) representation of the slice quiver consists of
\[
 B'_1,B'_2\in\End(V'),\,
 I'\in\Hom(W^{\mathrm{sl}},V'),\,
 J'\in\Hom(V',W^{\mathrm{sl}}),\,\, \text{where }
 \dim V'=m,\, \dim W^{\mathrm{sl}}=r,
\]
and its moment-map equation is
\[
 [B'_1,B'_2]+I'J'=0.
\]
The restriction of the original determinant stability character to \(G_{x_z}=\GL(V')\) is the usual positive character in the ADHM construction. Therefore the quiver variety appearing in the slice theorem is the charge-$m$ instanton moduli
\[
 \M(W^{\mathrm{sl}},m),\quad \text{with affine quotient }\,\, \M_0(W^{\mathrm{sl}},m).
\]
As an ordinary variety, \(\M(W^{\mathrm{sl}},m)\) is \(\M(r,m)\); the notation \(W^{\mathrm{sl}}\) here indicates the modified equivariant framing weights, which will be computed below.

The other factor in the \cite[Theorem 3.2]{BC} is a smooth affine space \(U\). 
The relative slice theorem \cite[Theorem~3.2]{BC} therefore gives, \'etale locally at \(x_z\) and at the origin of the affine transverse slice $\M_0(W^{\mathrm{sl}},m)$, a common Cartesian \'etale presentation of the projective morphisms (\ref{eqn:GU_map})
\[
 \pi_n:\M(r,n)\longrightarrow \M_0(r,n)
\]
and
\begin{equation}\label{eqn:pi_n_local_model}
 \pi_m^{\mathrm{sl}}\times\id_U:
 \M(W^{\mathrm{sl}},m)\times U \longrightarrow \M_0(W^{\mathrm{sl}},m)\times U.
\end{equation}
More precisely, after shrinking around the relevant points, there exist a scheme \(\hat{Z}\), \'etale maps
\begin{equation}\label{eqn:etale_charts}
 p:\hat{Z}\longrightarrow \M_0(r,n),
 \qquad
 q:\hat{Z}\longrightarrow \M_0(W^{\mathrm{sl}},m)\times U,
\end{equation}
and an isomorphism of projective \(\hat{Z}\)-schemes
\begin{equation}\label{eqn:local_isom}
 \Psi:
 \M(r,n)\times_{\M_0(r,n)}\hat{Z}
 \xrightarrow{\ \sim\ }
 \bigl(\M(W^{\mathrm{sl}},m)\times U\bigr)\times_{\M_0(W^{\mathrm{sl}},m)\times U}\hat{Z}.
\end{equation}
There is a point \(\mathbf{0}\in \hat{Z}\) such that $p(\mathbf{0})=x_z\in \M_0(r,n)$ and \(q(\mathbf{0})=(0,0)\in\M_0(W^{\mathrm{sl}},m)\times U\). Corollary~3.4 and diagram (3.5) of \cite{BC} give the corresponding identification of formal completions, while Proposition~3.8 of \cite{BC} identifies the tautological bundles after pulling back to the common fiber product in (\ref{eqn:local_isom}).

\subsubsection{Equivariant refinement}

For the residue calculation, (\ref{eqn:etale_charts})--(\ref{eqn:local_isom}) must be made equivariant with respect to \(\sfK_w\). We now prove the needed refinement.

\begin{lemma}[Commuting-torus refinement of the relative slice]
\label{lem:equivariant-relative-slice}
In the setting \cite[Theorem~3.2]{BC} of quiver representations, let \(G\) be the gauge group acting on the representation space \(\mathbf M\) equipped with the symplectic structure $\omega$, with moment map
\[
 \mu:\mathbf M\longrightarrow\mathfrak g^\vee.
\]

Let \(\mathsf S\) be a torus acting linearly on \(\mathbf M\). Suppose that the $\mathsf{S}$-action commutes with the \(G\)-action, preserves the chosen GIT character, and satisfies
\begin{equation}\label{eqn:S_action}
 \sigma^*\omega=\chi(\sigma)\omega,
 \qquad
 \mu(\sigma x)=\chi(\sigma)\mu(x)
 \quad
 \text{for }\sigma\in\mathsf S,
\end{equation}
for an algebraic character \(\chi:\mathsf S\to\mathbb C^*\).

Let \(x\in\mu^{-1}(0)\). Suppose that 
\begin{enumerate}
\item 
the $G$-orbit of $x$ is closed;
\item
the  point \([x]\in \mu^{-1}(0)/G\) is fixed by $\mathsf{S}$ and admits an algebraic compensating homomorphism
\begin{equation}\label{eqn:rho_property}
 \rho:\mathsf S\longrightarrow G, \quad \sigma\cdot x=\rho(\sigma)\cdot x;
\end{equation}
\item     
\(\rho\) preserves the chosen decomposition of \(x\) into a direct sum of polystables.
\end{enumerate}

Then the common Cartesian \'etale cover in \cite[Theorem~3.2]{BC} can be chosen \(\mathsf S\)-equivariantly. The formal-completion isomorphisms of \cite[Corollary~3.4 and diagram~(3.5)]{BC} and the tautological-bundle isomorphism of \cite[Proposition~3.8]{BC} can likewise be chosen \(\mathsf S\)-equivariantly.
\end{lemma}

\begin{proof}
Define the compensated action of \(\mathsf S\) on $\mathbf{M}$ by
\begin{equation}\label{eqn:compensating action}
\sigma\star y := \rho(\sigma)^{-1}\cdot(\sigma\cdot y).
\end{equation}
Equation (\ref{eqn:rho_property}) says that \(x\) is fixed by this compensated action. Since multiplication by a gauge transformation does not change a point of the affine or projective quotient, the compensated action induces the same \(\mathsf S\)-action on the quotient by $G$ as the original action.

Because the orbit containing \(x\) is closed, the stabilizer group \(G_x=\{g\in G \ |\ g\cdot x=x\}\) is reductive. The compensated \(\mathsf S\)-action normalizes \(G_x\), and hence we obtain a reductive semidirect product
\[
H=G_x\rtimes\mathsf S.
\]
All tangent spaces and linear maps used in the construction of the relative slice are \(H\)-equivariant. In particular, the tangent space of the orbit containing $x$, namely
\[
 \mathfrak g\cdot x\subset T_x\mathbf M,
\]
the kernel and image of \(d\mu_x\), the derivative of $\mu$ at $x$, and the relevant orthogonal complements with respect to $\omega$, are all \(H\)-stable. 

The symplectic form $\omega$ is an \(H\)-equivariant pairing with values in the $1$-dimensional representation \(\mathbb C_\chi\):
\[
 \omega:T_x\mathbf M\otimes T_x\mathbf M \longrightarrow\mathbb C_\chi.
\]
Here \(\chi\) is trivial on \(G_x\). We explain explicitly why all $\omega$-orthogonal complements can be chosen \(H\)-equivariantly. Let \(E\) be an \(H\)-module equipped with a nondegenerate alternating pairing
\[
 E\otimes E\longrightarrow\mathbb C_\chi,
\]
and let \(A\subset E\) be an \(H\)-stable isotropic subspace. Consider $A^\perp\subset E$ defined by the pairing. Since $A$ is isotropic, $A\subset A^\perp$. Reductivity gives an \(H\)-equivariant splitting,
\[
 A^\perp=A\oplus N,
\]
where \(N\) maps isomorphically onto \(A^\perp/A\). The restriction of the pairing to \(N\) is nondegenerate. Therefore $E=N\oplus N^\perp$, and $A\subset N^\perp$ is Lagrangian. Choose an
$H$-stable complement $B_0$ to $A$ \emph{inside $N^\perp$}: $N^\perp=A\oplus B_0$. Then $E=N\oplus N^\perp=N\oplus A\oplus B_0=A^\perp\oplus B_0$ shows that $B_0$ is also a complement of $A^\perp$ in $E$ and is orthogonal to $N$. The pairing between \(A\) and \(B_0\) is perfect. Hence there is a
unique \(H\)-equivariant map
\[
\phi:B_0\longrightarrow A
\]
satisfying
\[
\omega\bigl(\phi(b),b'\bigr)=-\frac12\omega(b,b')
 \qquad
\text{for all }b,b'\in B_0.
\]
The graph
\[
 B=\{b+\phi(b)\ |\ b\in B_0\}
\]
is isotropic, and
\[
 E=A\oplus N\oplus B
\]
is an \(H\)-equivariant Witt decomposition. This gives, in particular, simultaneous \(G_x\)- and \(\mathsf S\)-stable orbit, moment-map, coisotropic, and symplectic-normal complements in the proof of the relative slice theorem. The factor \(\mathbb C_\chi\) just twists the dual summands and does not affect the argument.

We now identify the actual slice maps to which the above linear algebra consideration applies. This is an equivariant adaptation of \cite[Section~3.3, Lemmas~3.5-3.6 and Theorem~3.7]{BC}. Take $A=\mathfrak g\cdot x$ to be the tangent space of the $G$-orbit at $x$ and set $C=N\oplus B$.  Then $C$ is an $H$-stable coisotropic complement to the tangent space, with $C^\perp=B$ and $C/C^\perp=N$. Let $\mathfrak{g}_x\subset \mathfrak{g}$ be the Lie algebra of the stabilizer group $G_x\subset G$. Choose also an $H$-stable decomposition $\mathfrak g=\mathfrak g_x\oplus L$. Write $c=n+b\in C$ and define
\[
 \nu: C\to L^\vee, \,\,\nu(c)=\mu(x+c)|_L,\qquad
 \mu_N: N\to \mathfrak{g}_x^\vee,\,\, \mu_N(n)=\mu(x+c)|_{\mathfrak g_x}.
\]
The second expression is independent of $b$: elements of $\mathfrak g_x$ annihilate $x$, and $B$ is the radical of $\omega|_C$. It is the quadratic moment map of the symplectic-normal representation $N$. The zero locus used in the slice construction is therefore
\[
 Z_0=\{c=n+b\in C\ | \ \mu_N(n)=0,\ \nu(c)=0\}.
\]
Both equations are $H$-equivariant, with the coadjoint representations twisted by $\chi$. The two maps are explicitly
\[
 p_0:G\times_{G_x}Z_0\longrightarrow\mu^{-1}(0), \quad [g,c]\longmapsto g\cdot(x+c),
 \qquad
 q_0:Z_0\longrightarrow\mu_N^{-1}(0),\quad n+b\longmapsto n.
\]
For the compensated action on the target, the $\mathsf S$-action on the induced space is
$$\sigma\star[g,c]=[\rho(\sigma)^{-1}g\rho(\sigma),\sigma\star c];$$
this verifies equivariance of $p_0$ directly. Equivariance of $q_0$ follows from the $H$-stable splitting. The derivative of $\nu$ on $B$ is an isomorphism to $L^\vee\otimes\mathbb C_\chi$, because the pairing between $B$ and the tangent space $A$ is perfect. Thus projection along $B$ gives precisely the \'etale slice projection at the distinguished point. The quiver factor and the other factor in the  representation can be chosen using the polystable summands;
the resulting direct-sum decomposition is preserved by $\mathsf S$.

The remaining issue is the choice of saturated neighborhoods and semistable loci, which is necessary for the projective comparison. The loci on which these maps and their quotient squares are strongly \'etale are invariant under $\mathsf S$: the maps are equivariant and
the defining Cartesian property is preserved by automorphisms. More explicitly, on each source affine quotient, take the union of opens on which the natural comparison with the quotient fiber product is an isomorphism and the quotient map is \'etale. The lemmas of \cite{BC} show that this good locus contains the distinguished point, and equivariance makes it $\mathsf S$-stable. Choose an invariant principal open within this locus, as justified below, and pull it back to obtain a saturated open. Thus the shrinkings in \cite[Lemmas~3.5--3.6]{BC} can be made invariant. The restriction of the GIT character is unchanged, and its semi-invariant spaces are preserved by $\mathsf S$. The base-change identifications of those spaces used to compare semistability therefore remain equivariant. This gives the two Cartesian squares of
\cite[Theorem~3.7]{BC} in the equivariant category.

For completeness, we show that invariant shrinkings exist as follows. Let \(O\) be an \(H\)-stable open neighborhood of the fixed point in an affine \(H\)-scheme, and let \(I\) be the ideal of its complement. Choose \(f\in I\) with \(f(x)\neq0\). The Reynolds projection (namely taking invariants) gives
\[
f_0\in I^H,\qquad f_0(x)=f(x)\neq0.
\]
Hence the invariant principal open \(D(f_0)\) contains \(x\) and is contained in \(O\).

Use the same argument with $H=\mathsf S$ on an affine quotient when a saturated shrinking is required. Passing to completions preserves equivariance. Finally, the comparison of tautological bundles in \cite[Proposition~3.8]{BC} comes from associated bundles by restricting the defining gauge representations to $G_x$. The polystable decomposition is preserved by $\rho$ by assumption, so those associated-bundle identifications respect $\mathsf S$ as well. This proves each claimed compatibility.
\end{proof}

\subsection{Equivariant local model}\label{sec:equiv_local_model}

We apply Lemma~\ref{lem:equivariant-relative-slice} to $\pi_n$ with $\mathsf S=\sfK_w$. No finite covering of \(\sfK_w\) is needed in the present ADHM situation, because the compensating homomorphism can be written explicitly.

Consider the ADHM monad in (\ref{eqn:polystable_pt}). Let
\[
 e_{uv}=x^uy^v, \qquad  0\leq u<i,\quad 0\leq v<j,
\]
be the monomial basis of \(R_z\). For
\[
 k=(t_1,t_2;u_1,\ldots,u_r)\in \sfK_w,
\]
define
\begin{equation}\label{eqn:rho_R}
 \rho_R(k)e_{uv}= u_\alpha^{-1}t_1^ut_2^v e_{uv}.
\end{equation}
This is an algebraic homomorphism
\[
 \rho_R: \sfK_w\longrightarrow\GL(R_z).
\]
A direct calculation gives
\[
 \rho_R(k)B_1\rho_R(k)^{-1}=t_1B_1,\quad 
 \rho_R(k)B_2\rho_R(k)^{-1}=t_2B_2,\quad 
 \rho_R(k)I_\alpha=u_\alpha^{-1}I_\alpha.
\]
At the opposite corner,
\begin{equation}\label{eqn:rho_R2}
 \begin{aligned}
 J_\beta\rho_R(k)^{-1} =u_\alpha t_1^{1-i}t_2^{1-j}J_\beta =u_\beta t_1t_2J_\beta,
 \end{aligned}
\end{equation}
where the last equality is precisely the defining relation
\[
 \frac{u_\alpha}{u_\beta}=t_1^it_2^j
\]
on \(\sfK_w\). Equations (\ref{eqn:rho_R})--(\ref{eqn:rho_R2}) show that
\[
 k\cdot R_z=\rho_R(k)\cdot R_z.
\]
The zero unframed simple \(S\) is fixed without compensation. Thus
\[
 \rho(k)=\rho_R(k)\oplus\id_{S^{\oplus m}}
\]
is a compensating homomorphism for the polystable datum \(R_z\oplus S^{\oplus m}\), and it preserves the two summand types.

For the ADHM representation space, $\omega=\operatorname{tr}\bigl(dB_1\wedge dB_2+dI\wedge dJ\bigr)$, and $\mu=[B_1,B_2]+IJ$. They satisfy
\[
 k^*\omega=\hbar\,\omega, \qquad \mu(k\cdot x)=\hbar\,\mu(x), \qquad \hbar=t_1t_2.
\]
Thus Lemma~\ref{lem:equivariant-relative-slice} applies even though \(\sfK_w\) need not be contained in \(\ker(\hbar)\): the construction is valid when the symplectic form is preserved by the group action up to a character. We conclude that the common cartesian diagram (\ref{eqn:etale_charts})--(\ref{eqn:local_isom}) may be chosen \(\sfK_w\)-equivariantly.

We now compute the \(\sfK_w\)-equivariant transverse framing $W^{\text{sl}}$. The deformation complex of an ADHM datum \((B_1,B_2,I,J)\) is 
\begin{equation}\label{eqn:ADHM_object}
 \End(V) \xrightarrow{\,d_0\,}
 \begin{matrix}
 t_1\End(V)\oplus t_2\End(V)\\
 {}\oplus\Hom(W,V)\oplus
 \hbar\Hom(V,W)
 \end{matrix}
 \xrightarrow{\,d_1\,} \hbar\End(V),
\end{equation}
where
\[
 d_0(\xi) = \bigl([\xi,B_1],[\xi,B_2],\xi I,-J\xi\bigr),\quad
 d_1(b_1,b_2,\iota,\jmath) = [b_1,B_2]+[B_1,b_2]+\iota J+I\jmath.
\]
Its virtual \(K\)-theory character is
\begin{equation}\label{eqn:K-theory-character}
 \mathsf T(V,W) = W^*V+\hbar V^*W-CV^*V, \qquad C=(1-t_1)(1-t_2).
\end{equation}

Let \(\mathsf R\) denote the \(\sfK_w\)-character of the rectangular summand $R_z$. The compensating gauge action on the tautological fiber is the inverse of \(\rho_R\). Using the wall relation $\frac{u_\alpha}{u_\beta}=t_1^it_2^j$ in the definition of $\sfK_w$ in (\ref{eqn:torus_w}) and reindexing the monomial basis, we obtain
\begin{equation}\label{eqn:taut_char_R}
 \mathsf R =
 e^{a_\beta}\hbar\frac{1-t_1^i}{1-t_1} \frac{1-t_2^j}{1-t_2}
 \quad
 \text{in }K_{\sfK_w}(\mathrm{pt}).
\end{equation}
For \(V=\mathsf R\oplus V'\), define in \(K_{\sfK_w}(\mathrm{pt})\)
\begin{equation}\label{eqn:W_sl}
 W^{\mathrm{sl}}:= W-\hbar^{-1}C\mathsf R.
\end{equation}
Since
\[
 C^*=(1-t_1^{-1})(1-t_2^{-1})=\hbar^{-1}C,
\]
we have $(W^{\mathrm{sl}})^*=W^*-C\mathsf R^*$. Expanding (\ref{eqn:K-theory-character}) therefore gives
\begin{equation}\label{eqn:K-char1}
 \mathsf T(\mathsf R\oplus V',W)
 = \mathsf T(V',W^{\mathrm{sl}})
 + \bigl(W^*\mathsf R+\hbar\mathsf R^*W-C\mathsf R^*\mathsf R\bigr).
\end{equation}
The construction of $U$ implies that the second summand above is the character of the factor \(U\):
\begin{equation}\label{eqn:char_U}
 U = W^*\mathsf R+\hbar\mathsf R^*W-C\mathsf R^*\mathsf R.
\end{equation}
Although (\ref{eqn:char_U}) is written as a virtual character, it is the character of an actual vector space. Indeed, it is the middle cohomology of the ADHM deformation complex of \(R_z\); stability and costability imply that the degree-$0$ and degree-$2$ cohomology groups vanish.

The expression \(W^{\mathrm{sl}}\) in (\ref{eqn:W_sl}) is likewise an honest rank-\(r\) representation. Indeed, $\hbar^{-1}C\mathsf R = e^{a_\beta}(1-t_1^i)(1-t_2^j)$, and the wall relation reads $e^{a_\alpha}=e^{a_\beta}t_1^it_2^j$. Consequently
\begin{equation}\label{eqn:W_sl2}
 W^{\mathrm{sl}} = e^{a_\beta}t_1^i + e^{a_\beta}t_2^j + \sum_{\gamma\neq\alpha,\beta}e^{a_\gamma}.
\end{equation}
Thus the original bottom-left and top-right framing characters $e^{a_\beta}, e^{a_\alpha} = e^{a_\beta}t_1^it_2^j$ are replaced by the two opposite corners
\[
 e^{a_\beta}t_1^i, \quad e^{a_\beta}t_2^j.
\]
This is the corner flip\footnote{Here the \emph{corner flip} denotes the following change of equivariant framing. On the wall
\[
 e^{a_\alpha}=e^{a_\beta}t_1^it_2^j,
\]
the two original framing characters
\[
 e^{a_\beta}, \quad e^{a_\alpha}=e^{a_\beta}t_1^it_2^j
\]
occupy the southwest and northeast corners of the corresponding \(i\times j\) character rectangle.  Removing the rectangular ADHM summand replaces them in the transverse slice by the opposite two corner
characters
\[
 e^{a_\beta}t_1^i, \quad  e^{a_\beta}t_2^j.
\]
Equivalently,
\[
 \{a_\beta,a_\alpha\} \longmapsto \{a_\beta+i\epsilon_1,a_\beta+j\epsilon_2\}.
\]
This changes only the equivariant framing of the lower-charge ADHM factor, not its rank or its underlying nonequivariant variety.} associated with removal of the \(i\times j\) rectangle.

The Cartesian \'etale identification follows from \cite[Theorem~3.2]{BC}. Its symplectic normal representation is identified in (\ref{eqn:K-char1}), which also records the \(\sfK_w\)-equivariant structure. The equivariant form of (\ref{eqn:pi_n_local_model}) is 
\begin{equation}\label{eqn:equiv_pi_n_local_model}
 \M(r,n)
 \mathrel{\sim_{\mathrm{et},\,\sfK_w}}
 \M(W^{\mathrm{sl}},m)\times U
 \quad\text{over}\quad
 \M_0(r,n)
 \mathrel{\sim_{\mathrm{et},\,\sfK_w}}
 \M_0(W^{\mathrm{sl}},m)\times U,
\end{equation}
which is shorthand for the common Cartesian cover (\ref{eqn:etale_charts})--(\ref{eqn:local_isom}).

\subsection{Calculations on local model}

Every stable map to \(\M(r,n)\) is constant after composing with the Gieseker-Uhlenbeck map $\pi_n$ to the affine quotient \(\M_0(r,n)\).  
The same statement holds universally in families. Thus the stack of stable maps to \(\M(r,n)\) is canonically the stack of stable maps to
the fibers of \(\pi_n\), together with the underlying locus in the affine quotient. In view of (\ref{eqn:etale_charts})--(\ref{eqn:local_isom}), we obtain a \(\sfK_w\)-equivariant isomorphism
\begin{equation}\label{eqn:isom_stable_map_moduli1}
 \begin{aligned}
 \Mbar_{g,k}\bigl(\M(r,n),d\bigr)\times_{\M_0(r,n)}\hat{Z}
 \xrightarrow{\ \sim\ }
\Mbar_{g,k}\bigl(\M(W^{\mathrm{sl}},m)\times U,d\bigr)\times_{\M_0(W^{\mathrm{sl}},m)\times U}\hat{Z}.
 \end{aligned}
\end{equation}
Because every map from a connected complete curve to \(U\) is constant,
there is also a canonical identification
\begin{equation}\label{eqn:isom_stable_map_moduli2}
 \Mbar_{g,k}\bigl(\M(W^{\mathrm{sl}},m)\times U,d\bigr)
 \cong
 \Mbar_{g,k}\bigl(\M(W^{\mathrm{sl}},m),d\bigr)\times U.
\end{equation}

\subsubsection{Obstruction theory}
Next, we compare obstruction theories. Tangent bundles pull back under \'etale maps in (\ref{eqn:etale_charts}) and isomorphism \(\Psi\) identifies the pullbacks of their tangent bundles. It follows that (\ref{eqn:isom_stable_map_moduli1}) identifies the pullbacks of the obstruction theories, and the virtual classes $[\Mbar_{g,k}\bigl(\M(r,n),d\bigr)]^{\mathrm{vir}}$ and $[\Mbar_{g,k}\bigl(\M(W^{\mathrm{sl}},m)\times U,d\bigr)]^{\mathrm{vir}}$ agree on the common cover (\ref{eqn:isom_stable_map_moduli1}). The same assertion holds $\sfK_w$-equivariantly.

On the product target, the tangent bundle decomposes as
\[
 T_{\M(W^{\mathrm{sl}},m)\times U}= T_{\M(W^{\mathrm{sl}},m)} \oplus  \bigl(U\otimes\mathcal O\bigr).
\]
Consequently the obstruction complex on the common cover (\ref{eqn:isom_stable_map_moduli1}) splits \(\sfK_w\)-equivariantly as\footnote{$f$ and $f'$ are appropriate universal stable maps.}
\begin{equation}\label{eqn:obs_theory_split}
 R\pi_*f^*T_{\M(r,n)}
 \cong
 R\pi_*f'^*T_{\M(W^{\mathrm{sl}},m)}
 \oplus
 R\pi_*\mathcal O_{\mathcal C}\otimes U.
\end{equation}

\subsubsection{Tautological bundles}\label{sec:taut_bundle_compatibility}

By \cite[Proposition~3.8]{BC}, the tautological bundles are compatible with the common relative slice. On the common target in (\ref{eqn:local_isom}), there is a \(\sfK_w\)-equivariant decomposition
\begin{equation}\label{eqn:taut_bundle_decomp}
\mathcal V_n \cong \mathcal V_m^{\mathrm{sl}}\oplus\mathcal R_Z,
\end{equation}
where \(\mathcal V_m^{\mathrm{sl}}\) is the rank-\(m\) tautological bundle of the slice \(\M(W^{\mathrm{sl}},m)\) and \(\mathcal R_Z\) is a rank-\(ij\) bundle pulled back from the base \(\hat{Z}\). Its fiber at the distinguished point is the rectangular \(\sfK_w\)-representation
\(\mathsf R\) of (\ref{eqn:taut_char_R}).

After a further \(\sfK_w\)-invariant \'etale shrinking, we may trivialize \(\mathcal R_Z\) as a vector bundle compatibly with its character decomposition. Taking determinants in (\ref{eqn:taut_bundle_decomp}), with every class restricted to $\sfK_w$, the torus acting on the common cover, gives
\begin{equation}\label{eqn:taut_div_decomp}
 H_n|_{\sfK_w} = c_1^{\sfK_w}(\det\mathcal V_n)
 = c_1^{\sfK_w}(\det\mathcal V_m^{\mathrm{sl}}) + c_1^{\sfK_w}(\det\mathcal R_Z).
\end{equation}
The term $c_1^{\sfK_w}(\det\mathcal R_Z)$ is pulled back from the affine base. Hence it has degree $0$ on every curve contained in a fiber of the GIT morphism. Hence $\int_{d\beta} H_n= \int_{d\beta} H_m^{\mathrm{sl}}$, so the degree-\(d\) components on the two sides of (\ref{eqn:isom_stable_map_moduli1}) correspond.

At a fixed base point, \(c_1^{\sfK_w}(\det\mathcal R_Z)\) is an equivariant scalar \(c\). Thus the divisor $H_n$ restricts to
\begin{equation}\label{eqn:H_restriction}
 H_n\big|_{\mathrm{slice}} = H_m^{\mathrm{sl}}+c\,1.
\end{equation}

\subsubsection{The \(\sfK_w\)-fixed part of the factor $U$}\label{sec:K_w-fixed part leaf}

We next show that the factor \(U\) contains exactly one \(\sfK_w\)-fixed line. First consider the affine quotient $\M_0(W^{\mathrm{sl}},m)$ of the transverse slice. The framing characters in \(W^{\mathrm{sl}}\) are
\begin{equation}\label{eqn:chars_W_sl}
 \widetilde u_1=u_\beta t_1^i,
 \qquad
 \widetilde u_2=u_\beta t_2^j,
 \qquad
 u_\gamma, \quad \gamma\neq\alpha,\beta.
\end{equation}
As in Sections \ref{sec:coord_ring}-\ref{sec:K_w-fixed}, the coordinate ring of $\M_0(W^{\mathrm{sl}},m)$ is generated by nonempty trace words and framed paths. A nonempty trace word of bidegree \((p,q)\) has character $t_1^pt_2^q$ and cannot restrict trivially to \(\sfK_w\).

At the generic point of the wall $w=0$, a framed path involving one of the unchanged colors \(u_\gamma\), 
$\gamma\neq\alpha,\beta$, contains an independent framing character and therefore also has nontrivial \(\sfK_w\)-weight. It remains only to inspect paths between the two corner-flipped colors. A framed path of bidegree \((p,q)\), where \(p,q\geq0\), has one of the two characters
\begin{equation}\label{eqn:pq-char1}
 \frac{\widetilde u_1}{\widetilde u_2}
 t_1^{p+1}t_2^{q+1} = t_1^{p+i+1}t_2^{q-j+1},
\end{equation}
or
\begin{equation}\label{eqn:pq-char2}
 \frac{\widetilde u_2}{\widetilde u_1}
 t_1^{p+1}t_2^{q+1} = t_1^{p-i+1}t_2^{q+j+1}.
\end{equation}
The first character cannot be trivial because its \(t_1\)-exponent is strictly positive, and the second cannot be trivial because its \(t_2\)-exponent is strictly positive. Paths beginning and ending at the same flipped color have character $t_1^{p+1}t_2^{q+1}$, which is also nontrivial. Hence no nonconstant homogeneous generator of \(\mathcal O(\M_0(W^{\mathrm{sl}},m))\) has trivial \(\sfK_w\)-weight.

The scheme-theoretic fixed-locus description from Section \ref{sec:K_w-fixed} therefore gives
\begin{equation}\label{eqn:MW_sl_K_w-fixed}
 \M_0(W^{\mathrm{sl}},m)^{\sfK_w} = \{0\} \quad\text{scheme-theoretically}.
\end{equation}
It follows that
\begin{equation}\label{eqn:MW_sl_K_w-fixed2}
 \bigl(\M_0(W^{\mathrm{sl}},m)\times U\bigr)^{\sfK_w} = \{0\}\times U^{\sfK_w}.
\end{equation}
We compare this with the original affine base. By (\ref{eqn:K_w-fixed loci final}), $\M_0(r,n)^{\sfK_w}=\mathbb A^1_z$ scheme-theoretically near \(x_z\). Taking $\sfK_w$-fixed formal neighborhoods in (\ref{eqn:etale_charts}) gives
\begin{equation}\label{eqn:formal_loc_ring_isom}
 \widehat{\mathcal O}_{\mathbb A^1_z,x_z} \cong \widehat{\mathcal O}_{U^{\sfK_w},0}.
\end{equation}
The left-hand side is a one-variable regular complete local ring. Because \(U^{\sfK_w}\) is a vector space, (\ref{eqn:formal_loc_ring_isom}) implies
\begin{equation}\label{eqn:dim_fixed_subspace}
 \dim U^{\sfK_w}=1.
\end{equation}
Thus there is a \(\sfK_w\)-equivariant decomposition
\begin{equation}\label{eqn:equiv-decomp}
 U=L_w\oplus U_{\mathrm{mov}}, \qquad  L_w=U^{\sfK_w}\cong\mathbb C,
\end{equation}
where every weight of \(U_{\mathrm{mov}}\) is nontrivial on \(\sfK_w\). After localization at the prime \((w)\), the corresponding equivariant Euler classes are units.  The line \(L_w\) occurs with multiplicity $1$: this follows from the reduced $1$-dimensional fixed axis and not just from a tangent-character calculation.

Finally, we must distinguish \(\sfK_w\)-equivariance from equivariance for the residual torus. The quotient $\T/\sfK_w\cong\mathbb C^*$ does not fix the point \(x_z\) when \(z\neq0\), so the slice centered at \(x_z\) is not a linear \(\T/\sfK_w\)-representation. Globally, however, the axis coordinate transforms by
\[
 z\longmapsto\chi_w^{-1}z,
\]
and hence has additive \(\T/\sfK_w\)-weight
\begin{equation}\label{eqn:T/K_w-weight}
 a_\beta-a_\alpha+i\epsilon_1+j\epsilon_2=-w.
\end{equation}
Therefore the tangent bundle of the global axis is the equivariant line
\[
 T_{\mathbb A^1_z} \cong \mathcal O_{\mathbb A^1_z}\otimes\mathbb C_{-w},
\]
whose Euler class is \(-w\). The \(\sfK_w\)-equivariant slice proves that this is the unique fixed direction in $U$ and that it occurs once. The global \(\T\)-action on the axis gives its residual character \(-w\). All remaining weights in $U$ are units in the localization at \((w)\).

\subsubsection{Calculations}

Again put $m=n-ij$. To simplify notation, let $Y=\M(W^{\mathrm{sl}},m)$ denote the charge-$m$ instanton moduli appearing in the equivariant relative slice of Section \ref{sec:local_description_Uhlenbeck}. As an ordinary variety, \(Y\) is \(\M(r,m)\), but its equivariant framing representation is the
corner-flipped representation
\[
 W^{\mathrm{sl}}
 =  e^{a_\beta}t_1^i+ e^{a_\beta}t_2^j +\sum_{\gamma\neq\alpha,\beta}e^{a_\gamma}.
\]
Let $H_Y=c_1^\T(\det\mathcal V_Y)$ be its tautological divisor.

The calculations we do here are done coefficientwise in the Novikov variable and in the equivariant coefficient ring $\mathbb{Q}[\epsilon_1,\epsilon_2, a_1,...,a_r]$ localized at the prime \((w)\). Thus every character whose restriction to \(\sfK_w\) is nontrivial is a unit, while \(w\) itself is not inverted until the second localization with respect to \(\T/\sfK_w\).

By Section \ref{sec:K_w-fixed}, the \(\sfK_w\)-fixed locus in the affine base is the reduced axis $\mathbb A^1_z$, with residual \(\T/\sfK_w\)-weight \(-w\). In equivariant Chow group, we have
\begin{equation}\label{eqn:chow_class_equality}
 [0] = c_1^\T(\mathbb C_{-w})[\mathbb A^1_z] = -w[\mathbb A^1_z].
\end{equation}
Consequently, every contribution whose proper pushforward to \(\mathbb A^1_z\) is supported at \(z=0\) contains an additional factor of \(w\). Such a term is regular at \(w=0\), even if one subsequently divides by the Euler class \(-w\) of the axis.

It follows that the residue $\text{Res}_w$ may be computed entirely over the {\em punctured} 
axis \(\{z\neq 0\}\subset \mathbb{A}^1_z\), where the \(\sfK_w\)-equivariant relative slice theorem
applies. In particular, no information about the extension of the slice
across \(z=0\) can affect the residue.

Write
\[
\Mbar_{1,1}(Y,d)\xleftarrow{\ \pi\ }\mathcal C\xrightarrow{\ f_Y\ } Y
\]
for the universal family. Recall (\ref{eqn:isom_stable_map_moduli2}) and (\ref{eqn:obs_theory_split}). Let $\mathbb E=\pi_*\omega_\pi$ be the Hodge bundle. Since the fibers have arithmetic genus one, we have
\[
 R^0\pi_*\mathcal O_{\mathcal C}=\mathcal O,
 \qquad
 R^1\pi_*\mathcal O_{\mathcal C}=\mathbb E^\vee,
 \qquad
 c_1(\mathbb E)=\lambda_1.
\]
Thus, for a $1$-dimensional constant target direction \(\mathbb C_\nu\subset U\) of equivariant weight \(\nu\), its moving deformation is $R^0\pi_*\mathcal O_{\mathcal C}\otimes\mathbb C_\nu$, while its moving obstruction is $R^1\pi_*\mathcal O_{\mathcal C}\otimes\mathbb C_\nu = \mathbb E^\vee\otimes\mathbb C_\nu$. Its contribution to localization is
\begin{equation}\label{eqn:contribution1}
 \frac{e_\T(\mathbb E^\vee\otimes\mathbb C_\nu)}{e_\T(\mathbb C_\nu)}=\frac{\nu-\lambda_1}{\nu}.
\end{equation}

Recall the \(\sfK_w\)-equivariant decomposition (\ref{eqn:equiv-decomp}): $U=L_w\oplus U_{\mathrm{mov}}$, where $L_w=U^{\sfK_w}$ is $1$-dimensional and \(U_{\mathrm{mov}}\) has no \(\sfK_w\)-fixed weights. The line \(L_w\) is fixed in the first localization by \(\sfK_w\), but it is moving for the residual quotient torus $\T/\sfK_w\cong\mathbb C^*$. Its residual weight is the weight of the axis coordinate:
\[
 \operatorname{Wt}_{\T/\sfK_w}(L_w)=-w.
\]

Applying (\ref{eqn:contribution1}) with \(\nu=-w\), the obstruction contribution of this line is
\[
 e_\T(\mathbb E^\vee\otimes\mathbb C_{-w}) = -w-\lambda_1,
\]
while residual localization along the affine axis divides by
\[
 e_\T(\mathbb C_{-w})=-w.
\]
Hence the complete contribution of the distinguished line is
\begin{equation}\label{eqn:contribution2}
 \frac{-w-\lambda_1}{-w} = 1+\frac{\lambda_1}{w}.
\end{equation}
The sign is important: the axis has weight \(-w\), so the coefficient of \(w^{-1}\) is \(+\lambda_1\). The multiplicity of this factor is exactly one. Indeed, we have seen in (\ref{eqn:K_w-fixed loci final}) that the \(\sfK_w\)-fixed affine base is scheme-theoretically the reduced line $\mathbb{A}^1_z$, and in Section \ref{sec:K_w-fixed part leaf} we identify its tangent line with \(U^{\sfK_w}\). 

Write the \(\sfK_w\)-moving part of \(U\) as
\[
 U_{\mathrm{mov}} = \bigoplus_{\ell=1}^s\mathbb C_{\nu_\ell}.
\]
Every \(\nu_\ell\) above restricts nontrivially to \(\sfK_w\). Equivalently, $\nu_\ell\notin\mathbb Z w$, and therefore each \(\nu_\ell\) is a unit in the coefficient ring localized at \((w)\). By (\ref{eqn:contribution1}), the total contribution of these directions is
\begin{equation}\label{eqn:contribution3}
 \mathsf E_{\mathrm{mov}} = \prod_{\ell=1}^s \frac{\nu_\ell-\lambda_1}{\nu_\ell}.
\end{equation}
The Hodge class on the stack of stable maps is pulled back from \(\overline{\mathcal M}_{1,1}\), so $\lambda_1^2=0$. Hence,
\begin{equation}\label{eqn:contribution4}
 \mathsf E_{\mathrm{mov}} = 1-\lambda_1\sum_{\ell=1}^s\frac1{\nu_\ell},
\end{equation}
which is regular at \(w=0\). Multiplying (\ref{eqn:contribution2}) and (\ref{eqn:contribution4}) gives
\begin{equation}\label{eqn:contribution5}
 \begin{aligned}
 \left(1+\frac{\lambda_1}{w}\right)\mathsf E_{\mathrm{mov}}
 =
 \left(1+\frac{\lambda_1}{w}\right)\left(1-\lambda_1\sum_{\ell=1}^s\frac1{\nu_\ell}\right)
 =
 1+\frac{\lambda_1}{w}-\lambda_1\sum_{\ell=1}^s\frac1{\nu_\ell},
 \end{aligned}
\end{equation}
because every cross term contains \(\lambda_1^2\). Therefore
\begin{equation}\label{eqn:residue1}
 \operatorname{Res}_{w=0}
 \left[
 \left(1+\frac{\lambda_1}{w}\right)
 \mathsf E_{\mathrm{mov}}\,dw
 \right]
 =
 \lambda_1.
\end{equation}
Thus the moving directions in $U$ contribute regular factors but do not change either the coefficient or the sign of the residue $\text{Res}_w$.

The remaining factors in the localization formula are precisely the \(\sfK_w\)-localization factors for the stable maps to \(Y\). As in (\ref{eqn:MW_sl_K_w-fixed}), the affine quotient of \(Y\) has scheme-theoretic \(\sfK_w\)-fixed locus equal to $\{0\}$. Hence the \(\sfK_w\)-fixed locus of the stack of stable maps is proper: it lies over the fiber of the projective morphism $Y\to \M_0(W^{\mathrm{sl}},m)$ over $\{0\}$. It follows that its equivariant pushforward is regular at \((w)\). Equivalently, if one computes it by a further localization for \(\T/\sfK_w\), any apparent denominators cancel after summing all fixed components. Therefore, no additional \(w^{-1}\)-factor can arise from the charge-$m$ instanton moduli.

More explicitly, let \(F\) run over the \(\sfK_w\)-fixed components of \(\Mbar_{1,1}(Y,d)\). Virtual localization gives
\begin{equation}\label{eqn:localization_Y}
 [\Mbar_{1,1}(Y,d)]^{\mathrm{vir}}
 =
 \sum_F (\iota_F)_* \frac{[F]^{\mathrm{vir}}}{e_{\sfK_w}(N_F^{\mathrm{vir}})}.
\end{equation}
All denominators on the right-hand side are units in the localization at \((w)\). Since the axis residue contributes the common factor \(\lambda_1\) by (\ref{eqn:residue1}), summing the residual contributions over \(F\) reconstructs the full virtual localization on $Y=\M(W^{\text{sl}},m)$. This uses the complete \(\sfK_w\)-fixed stable-map obstruction theory; it does not require the image of a \(\sfK_w\)-fixed map to be pointwise fixed in \(Y\).

The compatibility of tautological bundles in the relative slice, discussed in Section  \ref{sec:taut_bundle_compatibility}, gives
\begin{equation}\label{eqn:taut_bundle_decomp2}
 \mathcal V_n \cong \mathcal V_Y\oplus \mathcal R,
\end{equation}
where \(\mathcal R\) is the rank-\(ij\) rectangular summand and \(\mathcal V_Y\) is the tautological bundle on the lower-charge instanton moduli $Y=\M(W^{\text{sl}},m)$. Taking determinants gives $\det\mathcal V_n  \cong \det\mathcal R\otimes\det\mathcal V_Y$ and hence
\begin{equation}\label{eqn:tut_div_decomp2}
 H_n\big|_{\mathrm{slice},\sfK_w}
 =H_Y\big|_{\sfK_w}+c_1^{\sfK_w}(\det\mathcal R).
\end{equation}
The rectangular bundle is pulled back from the \'etale slice base. Since every stable map lies in a fiber of the affine quotient, its pullback to the universal curve is constant. Thus $c_1^{\sfK_w}(\det\mathcal R)=c\,1$ on the stable-map fiber, for an equivariant scalar \(c\). In a residue formula we may choose any lift of $c$ to a $\T$-equivariant scalar: two lifts differ by a multiple of $w$, which cannot alter the residue. In particular, this summand has degree $0$ on every curve, so the slice identifies the degree-\(d\) class on \(\M(r,n)\) with the degree-\(d\) class on \(Y\).

Combining (\ref{eqn:residue1})--(\ref{eqn:tut_div_decomp2}), for every \(d>0\) one obtains
\begin{equation}\label{eqn:residue2}
 \begin{aligned}
 \operatorname{Res}_{w=0}
 \left([Q^d]\mathcal F^+_{r,n}(Q)\,dw\right)=
 (-1)^{rd}\int_{[\overline{\mathcal{M}}_{1,1}(Y,d)]^{\mathrm{vir}}}\lambda_1\,\operatorname{ev}_1^*(H_Y+c\,1).
 \end{aligned}
\end{equation}
The factor $(-1)^{rd}$ is required because $[Q^d]\mathcal F^+_{r,n}$ is the signed geometric invariant. $Y$ also parametrizes framed sheaves of rank $r$, so the sign convention for its Novikov variable is the same. All integrals in these residue identities are equivariant localization integrals. The term $c1$ vanishes by the string equation (the Hodge class pulls back under forgetting the insertion), giving
\begin{equation}\label{eqn:residue3}
 \operatorname{Res}_{w=0} \left([Q^d]\mathcal F^+_{r,n}(Q)\,dw\right) 
 = (-1)^{rd}\int_{[\overline{\mathcal{M}}_{1,1}(Y,d)]^{\mathrm{vir}}}\lambda_1\,\operatorname{ev}_1^*H_Y.
\end{equation}
The nonseparating boundary expression for $\lambda_1$ gives
\[
 \int_{[\Mbar_{1,1}(Y,d)]^{\mathrm{vir}}}^{\T}
 \lambda_1\operatorname{ev}_1^*H_Y
 =\frac1{24}\sum_\mu\langle H_Y,\phi_\mu,\phi^\mu\rangle^Y_{0,3,d},
\]
where $\{\phi_\mu\}$ and $\{\phi^\mu\}$ are dual bases of localizaed equivariant cohomology. Since
\[
 [Q^d]\tau_{W^{\mathrm{sl}},m}(Q)
 =(-1)^{rd}\sum_\mu\langle H_Y,\phi_\mu,\phi^\mu\rangle^Y_{0,3,d},
 \qquad d>0,
\]
the signed residue identity is 
\begin{equation}\label{eqn:residue4}
 \operatorname{Res}_{w=0} \left([Q^d]\mathcal F^+_{r,n}(Q)\,dw\right)
 =
 \frac1{24}[Q^d]\tau_{W^{\mathrm{sl}},m}(Q),
\end{equation}
where
\[
 \tau_{W^{\mathrm{sl}},m}(Q) := \operatorname{Tr} \bigl(H_Y\star_Q-H_Y\cup\bigr).
\]

A wall in the space $\text{Spec}\,\CC[\epsilon_1,\epsilon_2,a_1,...,a_r]$ of framing parameters requires \(\alpha\neq\beta\), so necessarily \(r\geq2\).
For \(r\geq2\), the complete positive-degree Maulik-Okounkov divisor
formula \cite[Section 15.1]{MOk} is
\begin{equation}\label{eqn:MO_div_formula}
 \bigl(H\star_Q-H\cup\bigr)^+
 =
 h\sum_{k\geq1}
 k\frac{Q^k}{1-Q^k}\,
 \beta_{-k}\beta_k.
\end{equation}
There is no additional scalar term in this formula. In particular, (\ref{eqn:MO_div_formula}) depends on \(r\), \(h\), \(Q\), and the diagonal Heisenberg operators, but it is independent of the individual framing parameters \(a_\gamma\). The corner flip changes only the equivariant framing characters. It does not change the rank, the charge, or the Heisenberg expression
(\ref{eqn:MO_div_formula}). Therefore
\begin{equation}\label{eqn:trace_equal}
 \tau_{W^{\mathrm{sl}},m}(Q) = \tau_{r,m}(Q).
\end{equation}
Substituting (\ref{eqn:trace_equal}) into (\ref{eqn:residue4}) and summing over \(d>0\) yields
\[
 \operatorname{Res}_{w=0}\mathcal F^+_{r,n}(Q)\,dw
 =
 \frac1{24}\tau_{r,n-ij}(Q),
\]
which is what we want to prove.

\begin{remark}[Extreme cases]
If \(ij=n\), then \(m=0\), the charge-$m$ instanton moduli is a point, which has no positive-degree contributions. Our convention $\tau_{r,0}=0$ therefore gives zero residue, as required. The rank one case has no wall of the form considered here because there are no distinct colors \(\alpha\neq\beta\).
\end{remark}

{\large\part{The polynomial part}}
In this Part, we study the polynomial part of $\langle H\rangle_{1,1,d}^{\M(r,n)}$. As in (\ref{eqn:parameter_scaling}), set $a_\alpha=tb_\alpha, \alpha=1,...,r$, and assume $b_1,...,b_r$ are algebraically independent over $\mathbb Q(\epsilon _1,\epsilon _2)$. We consider the behavior of structures in Gromov-Witten theory of $\M(r,n)$ in the specialization $$t^{-1}=:s=0.$$  

Define
\[
 R=K[\![s,Q]\!],\quad \text{where } K=\mathbb{Q} (\epsilon_1,\epsilon_2,b_1,\ldots,b_r).
\]
We introduce a few terminologies to be used in this Part. 

Given a finite free $K(\!(s)\!)[\![Q]\!]$-module $\mathscr{H}$, a {\em Rees lattice} in $\mathscr{H}$ is a free $R$-submodule $\mathscr H_R\subset \mathscr{H}$ whose $Q$-adically completed localization\footnote{In particular, $\widehat{K[\![s,Q]\!][1/s]}^{\,(Q)}\simeq K(\!(s)\!)[\![Q]\!]$. Ordinary localization allows only a uniform bound on negative powers of $s$, while this completion permits a bound depending on the degree in the $Q$ term. For example, $\sum_{d\geq0}s^{-d}Q^d$ belongs to the completed ring but not to $K[\![s,Q]\!][1/s]$.} is $\mathscr H$:
\[
 \widehat{\mathscr H_R[1/s]}^{\,(Q)}\simeq\mathscr H.
\]
The quotient $\mathscr H_R/s\mathscr H_R$ is called the {\em large framing limit}. An endomorphism of $\mathscr{H}$ is called \emph{Rees regular} if its matrix belongs to $\End_R(\mathscr H_R)$, so it has a well-defined $s=0$ specialization. A  \emph{Rees projector} is an idempotent in $\End_R(\mathscr H_R)$. A \emph{Rees estimate} is an $s$-adic order estimate.

We assume $r\geq2$ in Sections~\ref{sec:Frob}-\ref{sec:limit_genus1}, except in Lemma ~\ref{lem:gamma-initial}: the constant-map calculation applies for every $r\geq1$. The Gromov-Witten theory in the rank-one cases $\M(1,N_\alpha)=\text{Hilb}^{N_\alpha}(\mathbb{C}^2)$ enters with its geometric scalar normalization, recalled in Section~\ref{sec:rank-one-limit}; it is not an instance of the higher-rank divisor decomposition below.

\section{Frobenius structure}\label{sec:Frob}

We study the $s=0$ specialization of the Frobenius structure of $QH_\T^*(\M(r,n))_{\loc}$. Write $\eta=\eta_{\M(r,n)}$ for the localized equivariant Poincar\'e pairing, and set
\begin{equation}\label{eqn:exponent}
 k:=n(r-1),\qquad c:=2k=2n(r-1).
\end{equation}
Recall the effective framing torus $\mathsf A$ from (\ref{eqn:eff_frame_torus}) and the decomposition of its fixed loci into components
\begin{equation}\label{eqn:framed_torus_fixed_loci2}
 \M(r,n)^{\mathsf A}=\coprod_{|\mathbf N|=n}F_{\mathbf N},
 \qquad
 F_{\mathbf N}:=\prod_{\alpha=1}^r
       \operatorname{Hilb}^{N_\alpha}(\mathbb C^2), \quad \mathbf{N}=(N_1,...,N_r).
\end{equation}
So $c$ is the codimension of any component of $\M(r,n)^{\mathsf{A}}$ in $\M(r,n)$. Let $i_{\mathbf N}:F_{\mathbf N}\hookrightarrow\M(r,n)$ denote the inclusions. For a chamber $\mathfrak C$ of $\mathsf A$, let $N^-_{\mathfrak C,\mathbf N}$ be the repelling half of the normal bundle $N_{\mathbf N}$; it has rank $k$.

We fix the instanton-moduli polarization of \cite[Section~12.1.5, equation~(12.2)]{MOk}:
\[
 \varepsilon_{\mathbf N}(a)
 =\prod_{\alpha=1}^r\prod_{\beta\ne\alpha}(a_\beta-a_\alpha)^{N_\alpha}.
\]
This is the product of the nonzero framing weights in the dual of the half-tangent class used in the sign convention for Novikov variable. Indeed, $(t_1-1)\End(\mathcal V)$ cancels upon restriction to the framing torus, and the moving part of $\Hom(\mathcal V,W)$ has weights
$a_\beta-a_\alpha$, with multiplicity $N_\alpha$ for $\beta\ne\alpha$. For every chamber, we choose the sign $\sigma_{\mathfrak C,\mathbf N}\in\{\pm1\}$ by
\[
 \sigma_{\mathfrak C,\mathbf N}e_{\mathsf A}(N^-_{\mathfrak C,\mathbf N})
=\varepsilon_{\mathbf N}(a).
\]
Here $e_{\mathsf A}$ means the product of the framing weights, with coordinate characters and positive-degree classes on $F_{\mathbf N}$ set to zero. Thus the same polarization is used for opposite chambers.

The stable envelope with this polarization is denoted by $\operatorname{Stab}_{\mathfrak C}$. Its diagonal restriction is
\begin{equation}\label{eqn:signed-stable-normalization}
 i_{\mathbf N}^*\operatorname{Stab}_{\mathfrak C}\big|_{H_\T^*(F_{\mathbf N})}
=\sigma_{\mathfrak C,\mathbf N} e_\T(N^-_{\mathfrak C,\mathbf N})\cup(-).
\end{equation}
This is the signed normalization in \cite[Definition~3.3.1 and Theorem~3.3.4(ii)]{MOk}.
The support and degree conditions make the restriction matrix triangular, with every off-diagonal restriction of framing degree (i.e. degree in $a_\alpha$'s) strictly less than $k$ \cite[Theorem~3.3.4(i),(iii)]{MOk}. After $a_\alpha=tb_\alpha$, its actual signed diagonal has expansion
\begin{equation}\label{eqn:euler_class_expand}
 \sigma_{\mathfrak C,\mathbf N} e_\T(N^-_{\mathfrak C,\mathbf N})
 =t^k\widetilde d_{\mathbf N}+O(t^{k-1}),\qquad
 \widetilde d_{\mathbf N}
 =\prod_{\alpha=1}^r\prod_{\beta\ne\alpha} (b_\beta-b_\alpha)^{N_\alpha}\ne0.
\end{equation}
In particular, $\widetilde d_{\mathbf N}$ is independent of $Q$ and of partition shapes within $F_{\mathbf N}$.

Every stable-envelope/Yangian $R$-matrix\footnote{These chamber $R$-matrices are distinct from the $\mathsf R(z)$ arising from Givental's work, which we use in Section~\ref{sec:limit_genus1}.} and transfer matrix below uses this same polarization, including on products with an auxiliary framing factor. For two chambers, the common leading diagonal gives
$\operatorname{Stab}_{\mathfrak C'}^{-1}
 \operatorname{Stab}_{\mathfrak C}=I+O(t^{-1})$.
For a root $R$-matrix, this is precisely the normalization of 
\cite[Definition~4.1.1 and Proposition~4.1.3]{MOk}. 

For a color-charge vector $\mathbf N$, write
\begin{equation}\label{eqn:H_N}
 \mathcal H_{\mathbf N}
 =\bigotimes_{\alpha=1}^r
 H_\T^\bullet\bigl(\operatorname{Hilb}^{N_\alpha}(\mathbb C^2)\bigr)_{\loc}.
\end{equation}
Let $\eta_{\mathbf N}$ be the product Poincar\'e pairing on this vector space, as used in the tensor product of the Gromov--Witten theories of $\text{Hilb}^{N_\alpha}(\mathbb{C}^2)$.

\begin{proposition}\label{prop:stable-lattice}
For $r\geq2$, there is a $Q$-independent isomorphism
\begin{equation}\label{eqn:stab_isom}
 \Phi_t:\bigoplus_{|\mathbf N|=n}\mathcal H_{\mathbf N}
 \longrightarrow H_\T^\bullet(\M(r,n))_{\loc},\qquad
 \Phi_t=\operatorname{Stab}_{\mathfrak C}\circ\,t^{-k}\circ
          \operatorname{diag}(\widetilde d_{\mathbf N}^{-1}),
\end{equation}
over $\mathbb Q(t,\epsilon_1,\epsilon_2,b_1,\ldots,b_r)$ such that
\begin{equation}\label{eqn:Id_expand}
 \Phi_t^{-1}(\id)=\bigoplus_{\mathbf N}\id_{\mathbf N}+O(s),
\end{equation}
\begin{equation}\label{eqn:metric_expand}
 t^{2n(r-1)}\Phi_t^*\eta_{\M(r,n)}
 =\bigoplus_{\mathbf N}c_{\mathbf N}\eta_{\mathbf N}+O(s),
 \qquad c_{\mathbf N}=(-1)^k\widetilde d_{\mathbf N}^{-2}\ne0,
\end{equation}
and
\begin{equation}\label{eqn:divisor_expand}
 \Phi_t^{-1}(H\star_Q)\Phi_t
 =tB+A^{\rm bd}(Q)+A^{\rm off}(Q).
\end{equation}
Here
\begin{equation}\label{eqn:divisor_expand1}
 B|_{\mathcal H_{\mathbf N}}=\nu_{\mathbf N}\id,
 \qquad \nu_{\mathbf N}=\sum_\alpha b_\alpha N_\alpha,
\end{equation}
$A^{\rm off}=O(1)$ changes the color-charge vector $\mathbf N$, and
\begin{equation}\label{eqn:divisor_expand2}
 A^{\rm bd}|_{\mathcal H_{\mathbf N}}
 =\sum_{\alpha=1}^r\left(1\otimes\cdots\otimes
       A^{\rm osc}_{N_\alpha}(Q)\otimes\cdots\otimes1\right).
\end{equation}
The operator $A_N^{\rm osc}$ is the geometric divisor operator of
$\operatorname{Hilb}^N(\mathbb C^2)=\M(1,N)$ with the following scalar subtracted:
\begin{equation}\label{eqn:divisor_expand3}
 hN\frac{Q}{1-Q}\id.
\end{equation}
 All $c_{\mathbf N}$ are independent of $Q$.
\end{proposition}

\begin{proof}
Let $i^*=\bigoplus_{\mathbf N}i_{\mathbf N}^*$ be the restriction map. The degree bound and (\ref{eqn:euler_class_expand}) give
\[
 i^*\operatorname{Stab}_{\mathfrak C}
 =t^k\bigl(\operatorname{diag}(\widetilde d_{\mathbf N})+O(s)\bigr),
 \qquad i^*\Phi_t=I+O(s), \quad \text{where } k=n(r-1).
\]
The diagonal leading matrix is invertible. Since $i^*\id$ is the direct sum of the product units, inversion yields
\begin{equation}\label{eqn:stab_env_id}
 \operatorname{Stab}_{\mathfrak C}^{-1}(\id)
 =\bigoplus_{\mathbf N}t^{-k}\widetilde d_{\mathbf N}^{-1}
       \id_{\mathbf N}+O(t^{-k-1}),
\end{equation}
and proves (\ref{eqn:Id_expand}).

For the metric, we use equivariant localization directly. The normal weights occur in pairs $\chi,h-\chi$, with nonzero opposite framing parts. Therefore
\[
 t^{-2k}e_\T(N_{\mathbf N})
 =(-1)^k\widetilde d_{\mathbf N}^{\,2}+O(s).
\]
For classes $\gamma=(\gamma_{\mathbf N})$ and $\delta=(\delta_{\mathbf N})$ in $\bigoplus_{|\mathbf N|=n}\mathcal H_{\mathbf N}$, localization and $i^*\Phi_t=I+O(s)$ give
\begin{align*}
 t^{2k}\eta(\Phi_t\gamma,\Phi_t\delta)
 &=\sum_{\mathbf N}\int_{F_{\mathbf N}}^\T
    \frac{(i_{\mathbf N}^*\Phi_t\gamma)
          (i_{\mathbf N}^*\Phi_t\delta)}
         {t^{-2k}e_\T(N_{\mathbf N})}\\
 &=\sum_{\mathbf N}(-1)^k\widetilde d_{\mathbf N}^{-2}
       \eta_{\mathbf N}(\gamma_{\mathbf N},\delta_{\mathbf N})+O(s).
\end{align*}
This proves (\ref{eqn:metric_expand}) for the pairings specified above.

Finally, apply the complete divisor formula of \cite[Section~15.1]{MOk} in the tensor product chart $\Phi_t$ with the fixed polarization. Its part linear in framing parameters is $\sum_\alpha a_\alpha N_\alpha=tB$. The bounded terms preserving $\mathbf N$ are the sum of the same-color oscillator divisors $A_N^{\text{osc}}$. Every cross-color term changes $\mathbf N$ and is  bounded as $t\to\infty$. In higher rank, there is no additional scalar; in rank one the geometric divisor has the scalar (\ref{eqn:divisor_expand3}). This gives the stated decomposition before normalization. Conjugation by
$\operatorname{diag}(\widetilde d_{\mathbf N}^{-1})$ leaves every block-diagonal term unchanged and multiplies a cross-block term by a nonzero $t$- and $Q$-independent ratio of the corresponding $\widetilde d$'s. The common power $t^{-k}$ cancels. Thus the same
block formula and bounds hold in the chart $\Phi_t$.
\end{proof}

\begin{remark}
For $r=2,n=1$, choose the chamber $\delta:=a_1-a_2>0$. On the components of charges $(1,0)$ and $(0,1)$ the repelling Euler factors are $h-\delta$ and $-\delta$, respectively, but the polarization signs are $+1$ and $-1$. The signed diagonal is therefore $(h-\delta,\delta)$. For the opposite chamber it is $(-\delta,h+\delta)$, with the same leading diagonal $(-\delta,\delta)$. The ordinary normal Euler factors have leading term $-\delta^2$, in agreement with $c_{\mathbf N}=-(b_1-b_2)^{-2}$.
\end{remark}

Proposition~\ref{prop:stable-lattice} gives a Rees lattice on which the unit and the rescaled metric are regular, and the divisor has the controlled pole $s^{-1}B$. The limit metric differs from the product metric by the nonzero constants $c_{\mathbf N}$. Their logarithmic $Q$-derivatives vanish, which is the property needed in the subsequent genus-one specialization.

\section{Baxter algebra and canonical projectors at the large framing limit}\label{sec:projector}

It is known \cite[Section 15.2]{MOk} that the divisior operator $H\star_Q$ has simple spectrum at generic values of equivariant parameters. We want to show that the primitive projectors of these eigenlines extend over the Rees lattice of Proposition~\ref{prop:stable-lattice} and that their $s=0$ specialization are tensor products of the corresponding projectors for $\operatorname{Hilb}^{N_\alpha}(\mathbb C^2)$. To do this, we place the divisor operator inside an algebra of commuting operators using the Baxter algebra construction in \cite{MOk}.

Put
\begin{equation}\label{eqn:max_ideal_in_R}
 R=K[\![s,Q]\!],\qquad \mathfrak m=(s,Q),\qquad
 K=\mathbb Q(\epsilon_1,\epsilon_2,b_1,\ldots,b_r),
\end{equation}
and, with the notation of (\ref{eqn:H_N}), set
\begin{equation}\label{eqn:Rees_lattice}
 \mathscr H_R=
 \left(\bigoplus_{|\mathbf N|=n}\mathcal H_{\mathbf N}\right)\otimes_K R.
\end{equation}
All regular endomorphisms are written in the normalized stable-envelope chart $\Phi_t$. There are three coefficient rings in the argument. The closed point of $\text{Spec}\, R$ corresponding to the maximal ideal $\mathfrak{m}$ has residue field $R/\mathfrak m=K$, the $s=0$ specialization with its Novikov variable retained has ring $R/(s)=K[\![Q]\!]$, and comparison with the original theory takes place over
\[
 \widehat R_s:=\varprojlim_m R[1/s]/(Q^m)
       =K(\!(s)\!)[\![Q]\!].
\]
This completion allows the order of a pole in $s$ to depend on the Novikov degree. 
The quantum cohomology ring at the generic parameters is considered over $K(\!(s)\!)(\!(Q)\!)$, which contains $\operatorname{Frac}(R)$. 

The classical state space $\mathscr H_R/\mathfrak m\mathscr H_R$ has a basis labeled by
tuples of partitions $\boldsymbol\lambda=(\lambda^1,\ldots,\lambda^r)$ of total size $n$. When a finite separable coefficient extension $K'/K$ is needed, we write $R'=K'[\![s,Q]\!]$.

\subsection{Divisor operator and Baxter algebra}
Fro completeness of the presentation, we describe in some details the construction of Baxter algebra and its relation to the divisor operator $H\star_Q$.
\subsubsection{Complete Novikov series and spectral operations}

We describe the representation spaces, the auxiliary trace, and the coefficient rings entering the Baxter algebra. Throughout this subsection, the charge $n$ is fixed, the geometric Novikov variable is $q$, and the variable used in the transfer-matrix twist is $Q=(-1)^r q$, as in (\ref{eqn:signed_novikov}). The spectral parameter introduced below is independent of both $Q$ and the large framing parameter $t=s^{-1}$.

\medskip\noindent
\textit{The physical Fock spaces and the energy grading.}
Let $\mathsf T_0=(\mathbb C^*)^2$ be the coordinate torus, and put $k_0=\mathbb Q(\epsilon_1,\epsilon_2)=\text{Frac}(H_{\mathsf{T}_0}^*(\text{pt},\mathbb{Q}))$. The underlying vector space of a rank-one Fock representation is
\[
 \mathsf F=\bigoplus_{m\geq0}\mathsf F[m],\qquad
 \mathsf F[m]
 =H^*_{\mathsf T_0}\bigl(\operatorname{Hilb}^{m}(\mathbb C^2);\mathbb Q\bigr)
     \otimes_{\mathbb Q[\epsilon_1,\epsilon_2]}k_0.
\]
The vacuum is $1\in\mathsf F[0]=k_0$. The Nakajima creation operators identify this space with the bosonic Fock space: a monomial with creation indices $\lambda_1,\ldots,\lambda_\ell$ has energy $\lambda_1+\cdots+\lambda_\ell$. Thus energy is the number of points, rather than cohomological degree. In particular, $\dim_{k_0}\mathsf F[m]=p(m)$, the number of partitions of $m$. See \cite[Sections~12.4.1--12.4.2 and~13.1.5]{MOk} for the geometric Fock-space identification, and \cite[Section~12.4.7, Lemma~12.4.3]{MOk}
for the number operator.

For each framing color $\beta=1,\ldots,r$, let $\mathsf F_\beta$ be a copy of this Fock space, with framing weight $a_\beta$. The weight specifies the spectral translate of the Yangian representation; it does not change the underlying vector space with the energy grading. This is the evaluation-parameter construction of
\cite[Sections~5.2.3--5.2.4, equations~(5.1)--(5.2)]{MOk}. Write $\mathsf N_\beta$ for the operator acting by $m$ on $\mathsf F_\beta[m]$. After extending scalars from $k_0$ to 
$K=\mathbb Q(\epsilon_1,\epsilon_2,b_1,\ldots,b_r)$, the physical space of total energy $n$ is
\[
 \mathsf V_n
 :=\left(\bigotimes_{\beta=1}^r\mathsf F_\beta\right)[n]
 =\bigoplus_{N_1+\cdots+N_r=n}
       \bigotimes_{\beta=1}^r\mathsf F_\beta[N_\beta]
 =\bigoplus_{|\mathbf N|=n}\mathcal H_{\mathbf N}.
\]
This is the tensor product description of the cohomology of the $\mathsf{A}$-fixed locus $\M(r,n)^{\mathsf{A}}$ in (\ref{eqn:framed_torus_fixed_loci}) (see \cite[Section~12.1.4]{MOk}), followed by the localized stable-envelope identification with the cohomology of $\M(r,n)$. The integers $N_\beta$ label individual color charges, whereas $n$ is their sum. Operators on $\mathsf V_n$ may change $\mathbf N$ while preserving $n$.

\medskip\noindent
\textit{Auxiliary representations and geometric $R$-matrices.}
Choose an auxiliary representation
\[
 F_0=\bigoplus_{d\geq0}F_0[d],\qquad
 N_0|_{F_0[d]}=d\,\operatorname{id},\qquad
 \dim F_0[d]<\infty.
\]
The basic choice is another rank-one Fock representation. More generally, an admissible auxiliary representation means a geometric Yangian module with this grading and the compatible rational $R$-matrices of
\cite[Sections~5.2.2--5.2.6 and~6.5.1]{MOk}. Finite tensor products of the Fock representations are examples, with total energy as grading; their evaluation parameters are included in the coefficient field when needed. The index $0$ always denotes an auxiliary factor and is not an additional physical framing color.

For a rank-one auxiliary factor, the two-factor geometric $R$-matrix is the change between opposite stable-envelope charts for the decomposition of a rank-two framing into two lines:
\[
 \mathsf R_{0\beta}(u)
 =\operatorname{Stab}_{-}^{-1}\operatorname{Stab}_{+},
 \qquad u=a_0-a_\beta.
\]
Here the chambers are ordered by the sign of $u$, and the two stable envelopes use the same polarization fixed in
Section~\ref{sec:Frob}. This is \cite[Definition~4.1.1 and Sections~4.1.5-4.1.6]{MOk}; the signed
normalization of the stable envelopes is \cite[Definition~3.3.1 and Theorem~3.3.4]{MOk}. We use the convention in which $\mathsf R_{0\beta}$ is an endomorphism of $F_0\otimes\mathsf F_\beta$, without a permutation of these factors. On each finite total-energy sector, its matrix entries are rational in $u$ and the equivariant parameters, and
\[
 \mathsf R_{0\beta}(u) =\operatorname{id}+O(u^{-1}),\qquad u\rightarrow\infty,
\]
by \cite[Section~4.1.3 and Proposition~4.1.3]{MOk}. These statements apply to tensor-product auxiliary modules by the factorization construction.

The relevant conservation law is
\[
 [\mathsf R_{0\beta}(u),N_0+\mathsf N_\beta]=0.
\]
The separate energies need not be preserved. Consequently, an $R$-matrix can exchange energy between its two factors. The grading-zero condition and the finiteness of each graded piece are precisely the assumptions in \cite[Section~5.2.5]{MOk}. The Yang--Baxter relation
\cite[Section~4.1.9, equation~(4.4)]{MOk} makes the constructions for
three or more factors compatible. These are the geometric
stable-envelope/Yangian $R$-matrices.\footnote{The $\mathsf{R}$-matrix for the semi-simple Gromov-Witten theory used later in this Part is a different construction.}

\medskip\noindent
\textit{Monodromy and the formal partial trace.}
First give the auxiliary factor an independent framing parameter $u$. Define the ordered monodromy operator
\[
 \mathsf L_0(u)
 =\mathsf R_{0r}(u-a_r)\cdots\mathsf R_{01}(u-a_1)
 \quad\text{on}\quad
 F_0\otimes\mathsf F_1\otimes\cdots\otimes\mathsf F_r.
\]
Each factor acts as the identity on all unlisted tensor factors, and the rightmost factor acts first. This ordered product is \cite[Section~5.2.6, equation~(5.6)]{MOk}; it implements the tensor-product Yangian action of \cite[Section~5.2.12]{MOk}. The twisted transfer matrix is the partial trace
\[
 \widetilde{\mathsf T}^{F_0}(u;Q)
 =\operatorname{Tr}_{F_0}\bigl[(Q^{N_0}\otimes
             \operatorname{id})\mathsf L_0(u)\bigr].
\]
The twist $Q^{N_0}$ acts by $Q^d$ on $F_0[d]$. It is the one-vertex specialization of the Cartan-torus twist in
\cite[Section~6.5.1, equations~(6.23)--(6.24)]{MOk}. For two auxiliary modules $F_0,F_{0'}$, the same variable $Q$ is used in both twists. Energy conservation gives
\[
 [\mathsf R_{00'}(u-u'),Q^{N_0}\otimes Q^{N_{0'}}]=0,
\]
which is the twist-compatibility condition in the transfer-matrix construction of \cite[Section~1.2.2]{MOk}.

The trace has a purely formal meaning. Let $\iota_d:F_0[d]\to F_0$ and $p_d:F_0\to F_0[d]$ be inclusion and projection. On $\mathsf V_n$, set
\[
 \widetilde{\mathsf T}^{F_0}(u;Q)
 =\sum_{d\geq0}Q^d\operatorname{Tr}_{F_0[d]}
 \bigl[(p_d\otimes\operatorname{id})\mathsf L_0(u)(\iota_d\otimes\operatorname{id})\bigr].
\]
The operator inside the trace starts with auxiliary energy $d$ and ends with the same auxiliary energy. Conservation of total energy then forces its physical output to have energy $n$. During the ordered product, all
states lie in the finite-dimensional vector space
\[
 \bigoplus_{e+m_1+\cdots+m_r=n+d}\left(
       F_0[e]\otimes\bigotimes_{\beta=1}^r\mathsf F_\beta[m_\beta]\right).
\]
This proves that the coefficient of every $Q^d$ is a finite matrix of rational functions. The trace here is the ordinary algebraic partial trace of an endomorphism of a
finite-dimensional space; choosing a basis and its algebraic dual gives the same operator independently of that basis. No analytic convergence or choice of an adjoint pairing is involved. This explains explicitly why the completion in \cite[Section~6.5.1]{MOk} is applicable here.

For example, for a rank-one auxiliary Fock space,
\[
 \operatorname{Tr}_{F_0}Q^{N_0}
 =\sum_{d\geq0}p(d)Q^d
 =\prod_{m\geq1}(1-Q^m)^{-1}.
\]
This scalar series is the leading term of the transfer matrix as $u\to\infty$. In particular, the normalization
$\mathsf R(u)=\operatorname{id}+O(u^{-1})$ does not normalize the auxiliary trace to have leading scalar $1$; we use the unnormalized trace throughout.

\medskip\noindent
\textit{Centering and the normalized stable-envelope chart.}
Fix a physical color $\alpha$. Substitute $a_\beta=tb_\beta$ and center the auxiliary parameter at $a_\alpha$,
that is, set $u=tb_\alpha+v$. Then
\[
 \widetilde{\mathsf T}_{\alpha,t}^{F_0}(v;Q)
 =\operatorname{Tr}_{F_0}\!\left[
 Q^{N_0}\mathsf R_{0r}\bigl(t(b_\alpha-b_r)+v\bigr)\cdots
 \mathsf R_{01}\bigl(t(b_\alpha-b_1)+v\bigr)\right].
\]
Thus the color-$\alpha$ factor has argument $v$, while the other arguments have a nonzero term proportional to $t$. Centering leaves the twist and the order of the factors unchanged.

This formula acts in the tensor product chart defined by $\operatorname{Stab}_{\mathfrak C}$. To express it in the chart $\Phi_t$ of Proposition~\ref{prop:stable-lattice}, define the invertible block-scalar operator
\[
 \mathsf D|_{\mathcal H_{\mathbf N}} =\widetilde d_{\mathbf N}^{-1}\operatorname{id},
 \qquad
 \Phi_t=\operatorname{Stab}_{\mathfrak C}s^{n(r-1)}\mathsf D,
\]
where $\widetilde d_{\mathbf N}$ is the signed coefficient in (\ref{eqn:euler_class_expand}). An operator written in the first chart is therefore conjugated by $s^{n(r-1)}\mathsf D$ to obtain its matrix in the second chart. The common scalar power of $s$ cancels, giving 
\begin{equation}\label{eqn:transfer_matrix}
 \mathsf T_{\alpha,t}^{F_0}(v;Q)
 =\mathsf D^{-1}\widetilde{\mathsf T}_{\alpha,t}^{F_0}(v;Q)\mathsf D.
\end{equation}
Although $\mathsf D$ is constant on each color-charge block, it need not commute with a transfer matrix that changes color charges. The conjugation in (\ref{eqn:transfer_matrix}) is therefore part of the definition.

\medskip\noindent
\textit{Spectral operations and the completed Baxter algebra.}
We now specify the completion and the allowed operations. After $t=s^{-1}$, use the ambient algebra
\[
 \mathscr E_n
 :=\operatorname{End}_{K(\!(s)\!)}(\mathsf V_n\otimes_K K(\!(s)\!))[\![Q]\!]
 \cong\operatorname{Mat}_{\dim\mathsf V_n} \bigl(K(\!(s)\!)[\![Q]\!]\bigr).
\]
In this algebra, the auxiliary data and the specialization at a fixed centered spectral value are required to be defined over $K(\!(s)\!)$. To allow additional independent scalar parameters, first choose one common coefficient field $L\supset K$ for the entire family and replace $K$ by $L$ throughout the definition. The specializations of spectral parameters in the rank-one cases used in the next subsection are already defined over $K$. The topology here is the $Q$-adic topology: agreement modulo $Q^m$ means agreement of all coefficients in degrees less than $m$. The coefficient ring is the completed localization $\widehat R_s$ specified above.

There are two spectral operations. First, for each fixed auxiliary energy, expand the rational matrix entries at infinity in an independent spectral variable. For instance,
\[
 \mathsf T_{\alpha,t}^{F_0}(v;Q)
 \sim\sum_{j\geq0}C_j^{F_0,\alpha}(Q)v^{-j},\qquad
 C_j^{F_0,\alpha}(Q)
 =\sum_{d\geq0}Q^d[v^{-j}]
       \bigl([Q^d]\mathsf T_{\alpha,t}^{F_0}(v;Q)\bigr).
\]
Each $C_j^{F_0,\alpha}(Q)$ is an element of $\mathscr E_n$ and keeps its entire $Q$-series. This is spectral coefficient extraction as in \cite[Section~6.5.1, equation~(6.23)]{MOk}, up to the harmless scalar factor $1/h$ used there, where $h=\epsilon_1+\epsilon_2$.

Second, a finite spectral evaluation at $v=v_*$ means evaluating the \emph{rational} matrix coefficient of each $Q^d$ at $v_*$. It is admissible when all these evaluations are defined in the chosen coefficient field. The value is their complete $Q$-series. This does not mean substituting $v_*$ into an arbitrary formal series in $v^{-1}$. The next subsection constructs an evaluation that is simultaneously defined in every auxiliary-energy sector.

Let $\mathcal G$ be the set consisting of these complete spectral coefficients and admissible finite evaluations, for the auxiliary representations under consideration. Define
\[
 \mathscr B_Q
 :=\overline{K(\!(s)\!)[\![Q]\!][\operatorname{id},\mathcal G]}^{\,(Q)}
       \ \subset\ \mathscr E_n.
\]
Here the expression inside the bar denotes the unital associative subalgebra generated by $\mathcal G$, and the bar denotes closure in $\mathscr E_n$. Equivalently, an element belongs to this closure if, modulo every $Q^m$, it can be represented by an expression in the generators with scalar coefficients. This is the represented, completed version of the Baxter algebra of \cite[Section~6.5.1, immediately after equation~(6.24)]{MOk}, with
admissible finite evaluations included explicitly.

The distinction between spectral and Novikov coefficients is essential. For two generators $C(Q)=\sum_d C_dQ^d$ and $C'(Q)=\sum_e C'_eQ^e$, an identity $[C(Q),C'(Q)]=0$ gives
\[
 \sum_{d+e=m}[C_d,C'_e]=0\qquad\text{for every }m\geq0;
\]
it does not give the vanishing of each individual commutator $[C_d,C'_e]$. Accordingly, the separate endomorphisms $[Q^d]C_j(Q)$ are not adjoined as additional generators. By contrast, extracting coefficients in two independent spectral variables, or making admissible finite evaluations of the rational coefficients, preserves the transfer-matrix commutation identity of \cite[Section~1.2.2, equation~(1.5), and Section~6.5.1]{MOk}.

Finally, the definition of $\mathscr B_Q$ takes place after completed localization and does not assert that every generator preserves $\mathscr H_R=\mathsf V_n\otimes_K R$. Rees regularity requires the stronger condition that every matrix entry belong to $K[\![s,Q]\!]$. The next subsection establishes this condition for the centered finite evaluations used to construct the projectors in the $s=0$ limit. The following proposition records the commutativity and the divisor-operator input needed for that construction.

\begin{proposition}\label{prop:baxter-divisor}
The following statements hold.
\begin{enumerate}
\item The algebra $\mathscr B_Q$ is commutative.
\item With $Q=(-1)^r q$ as in (\ref{eqn:signed_novikov}), the completely normalized divisor operator $H\star_Q$ belongs to $\mathscr B_Q$. In particular every complete transfer coefficient or admissible finite evaluation $C$ satisfies
\begin{equation}\label{eqn:commuting_with_baxter}
 [C,H\star_Q]=0.
\end{equation}
Spectral shifts and scalar subtraction preserve this commutation relation.
\item For $r\ge2$, the positive-degree part of the divisor operator is
\begin{equation}\label{eqn:divisor_quantum_part}
 (H\star_Q-H\cup)^+
 =h\sum_{k\ge1}k\frac{Q^k}{1-Q^k}\,\beta_{-k}\beta_k,
\end{equation}
with no extra scalar. On the rank-one, charge-$N$ sector (i.e. $\text{Hilb}^N(\mathbb{C}^2)$), the geometric divisor operator equals $A_N^{\rm osc}$ plus 
\begin{equation}\label{eqn:scalar_divisor}
 hN\frac{Q}{1-Q}\,\id.
\end{equation}
\item In generic parameter values the localized small quantum ring is generated by $H$, and $H\star_Q$ has simple spectrum after a splitting field extension. Its spectral projectors are the primitive multiplication projectors.
\end{enumerate}
\end{proposition}

\begin{proof}
The transfer-matrix identity of \cite[Sections 6.5.1 and 1.2.2]{MOk} is
\[
 [\mathsf T^{F_0}(u;Q),\mathsf T^{F_{0'}}(v;Q)]=0
\]
with independent spectral paramters $u,v$ and the same formal twist $Q$. Extracting spectral coefficients gives commutation of complete $Q$-series. Finite evaluations preserve the identity wherever they are defined, including its mixed form with a spectral coefficient in the other variable. These facts also imply commutativity after $Q$-adic closure. 

The identification with divisor operators is established in \cite[Proposition~6.5.2 and Theorem~10.2.1]{MOk}. After the change of Novikov variable (\ref{eqn:signed_novikov}), it gives an element $B_H\in \mathscr{B}_Q$ of the form $B_H=H\star_Q+c_H\id$. Scalars are included in our coefficient ring, so $H\star_Q=B_H-c_H\id$ belongs to $\mathscr{B}_Q$. Its scalar normalization is fixed by $(H\star_Q)(\id)=H$. The explicit Jordan-quiver calculation \cite[Section~15.1]{MOk} gives (\ref{eqn:divisor_quantum_part})--(\ref{eqn:scalar_divisor}).

The generation and simple-spectrum statements are proved in \cite[Theorem~15.2.1 and its proof]{MOk}. To check the last assertion, suppose $H\star_Q$ have distinct eigenvalues $\rho_1,\ldots,\rho_D$. Its spectral projectors are
\[
 E_i=\prod_{j\ne i}\frac{(H\star_Q)-\rho_j\id}{\rho_i-\rho_j}, \quad i=1,...,D.
\]
They are polynomials in a multiplication operator, and hence are multiplications by $e_i=E_i(\id)$. The projector identities imply $e_i\star_Q e_j=\delta_{ij}e_i$ and $\sum_i e_i=\id$. Since their images are  $1$-dimensional, $e_i$ are primitive. 
\end{proof}

\subsubsection{Finite spectral evaluation and factorization in the limit}

Put $k_0=\mathbb Q(\epsilon_1,\epsilon_2)$ and fix 
$$v_*=b_1,$$ 
which is transcendental over $k_0$. We use auxiliary rank-one Fock representations with their common framing character removed, so the rank-one finite-energy $R$-matrices have entries in $k_0(v)$. Evaluation at $v_*$ is therefore defined on every such rational entry, simultaneously in all finite-energy sectors.

For $\beta\ne\alpha$, put $d_{\alpha\beta}=b_\alpha-b_\beta\in K^\times$. On each finite-energy sector, the rational asymptotic \cite[Proposition~4.1.3]{MOk} gives
\[
 \mathsf R_{0\beta}(u)=1+O(u^{-1}),\qquad
 \mathsf R_{0\beta}(s^{-1}d_{\alpha\beta}+v_*)=1+O(s)
       \quad\hbox{in }\operatorname{End}(F_0\otimes\mathsf F_\beta)[\![s]\!].
\]
For example, a factor $u^{-1}$ becomes $s/(d_{\alpha\beta}+sv_*)$, whose denominator has nonzero constant term. More generally, rational functions regular at $u=\infty$ admit this substitution into $K[\![s]\!]$. The factor of color $\alpha$ is $\mathsf R_{0\alpha}(v_*)$, independent of $s$ and defined over $K$. It follows, degree by degree in the auxiliary trace, that
\[
 C_\alpha^{F_0}(s,Q):=\mathsf T_{\alpha,t}^{F_0}(v_*;Q)
       \in\operatorname{End}_R(\mathscr H_R).
\]
This assertion concerns the complete series: each $Q^d$ coefficient is regular at $s=0$, so together they define an element of $K[\![s,Q]\!]$. At $s=0$, all off-color factors become the identity. The resulting rank-one transfer preserves each color-charge block, so it commutes with the block-scalar normalization $\mathsf D$. Thus in the normalized chart,
\begin{equation}\label{eqn:transfer_matrix_expand}
 C_\alpha^{F_0}(0,Q)\big|_{\mathcal H_{\mathbf N}}
 =1\otimes\cdots\otimes
   \mathsf T_{(1),N_\alpha}^{F_0}(v_*;Q)
   \otimes\cdots\otimes1.
\end{equation}
This is an identity of endomorphisms with complete $Q$-series. No Laurent coefficient at $v=\infty$ is extracted after centering: that operation does not in general commute with $s=0$ specialization.

\subsubsection{Separation in the classical limit}\label{sec:seperation_in_limit}

At $Q=0$, the auxiliary trace is its vacuum matrix coefficient. For the positive stability chamber, the vacuum-transfer calculation \cite[Sections~4.7.3 and 5.4.1--5.4.2]{MOk} identifies its spectral coefficients with characteristic-class multiplication and shows that they generate the tautological characteristic classes. In case of the one-vertex Jordan quiver, the auxiliary rank-one representation suffices: its negative normal bundle is the tautological bundle, and the coefficients of the associated Euler-class ratio recover its characteristic classes over $k_0$, where $h$ is invertible.

At the fixed point indexed by $\lambda\vdash N$, the tautological weights after removing the common framing weights are
\begin{equation}\label{eqn:taut_chern_roots}
 \{-(i-1)\epsilon_1-(j-1)\epsilon_2\mid (i,j)\in\lambda\}.
\end{equation}
The multiset determines $\lambda$, since algebraic independence of $\epsilon_1,\epsilon_2$ recovers $(i,j)$ from each weight. Let $\tau_\lambda(v)\in k_0(v)$ be the vacuum-transfer eigenvalue. If $\tau_\lambda(v)=\tau_\mu(v)$, equality of their Laurent expansions implies equality of all the tautological characteristic classes, hence of the multisets (\ref{eqn:taut_chern_roots}), and therefore $\lambda=\mu$. This use of Laurent coefficients is solely an argument about the already specialized rank-one rational functions. Because $v_*$ is transcendental over $k_0$, the evaluation map $k_0(v)\rightarrow K$, $v\mapsto v_*$, is injective. Consequently, if $\lambda\ne\mu$, then $\tau_\lambda(v_*)\ne\tau_\mu(v_*)$. Thus finite evaluations apply to different partitions remain different.

Choose this auxiliary representation for each color and write $C_\alpha=C_\alpha^{F_0}$, $1\le\alpha\le r$. Adjoin $C_0=s(H\star_Q)$, which is Rees regular by Proposition~\ref{prop:stable-lattice}. In the specialization $s=Q=0$, its eigenvalue on $\mathcal H_{\mathbf N}$ is $\nu_{\mathbf N}=\sum_\alpha b_\alpha N_\alpha$. Algebraic independence of the $b_\alpha$ distinguishes the color-charge vectors. Within a fixed color-charge vector, the $C_\alpha$ distinguish every partition in its own color. Hence the eigenvalues of the collection $C_0,C_1,\ldots,C_r$ of commuting operators define a map 
\begin{equation}\label{eqn:joint_char}
 \boldsymbol\lambda=(\lambda^1,...,\lambda^r)\longmapsto
 \bigl(\nu_{\mathbf N},\tau_{\lambda^1}(v_*),\ldots,
                    \tau_{\lambda^r}(v_*)\bigr),
 \qquad N_\alpha=|\lambda^\alpha|, \alpha=1,...,r,
\end{equation}
which is injective in the specialization $s=Q=0$. Commutativity follows from Proposition~\ref{prop:baxter-divisor} over the generic coefficient ring and therefore also over $R$, since the latter embeds in that ring.

\subsection{Regular primitive projectors and their limits}
We begin with an algebraic lemma.
\begin{lemma}[Complete local projector]\label{lem:complete_local_projector}
Let $(A,\mathfrak n,k)$ be a complete Noetherian local domain with a coefficient field $k$ of characteristic zero, and let $M$ be a free $A$-module of rank $D$. Given commuting operators  $C_0,\ldots,C_l\in\operatorname{End}_A(M)$. Suppose that, after a finite separable extension $k'/k$, $(M/\mathfrak nM)\otimes_k k'$ is a direct sum of $D$ joint eigenlines of $C_0,...,C_l$ so that the $(l+1)$-tuples of eigenvalues of $C_0,...,C_l$ are distinct. Put $A'=A\widehat\otimes_k k'$. Then the mod $\mathfrak{n}$ projectors to eigenlines have canonical lifts to orthogonal rank-one idempotents
\[
 P_1,\ldots,P_D\in A'[C_0,\ldots,C_l]
       \subset\operatorname{End}_{A'}(M\otimes_A A'),
 \qquad\sum_iP_i=\id.
\]
Their images are joint eigenlines. The lifts are unique among idempotents in the indicated commutative algebra with the given labels.
\end{lemma}

\begin{proof}
Write $c_j(i)\in k'$ for the eigenvalue of $C_j \bmod \mathfrak{n}$ on the $i$-th eigenline. Because the tuples of eigenvalues are distinct and $k'$ is infinite, there exist $u_0,\ldots,u_l\in k'$ so that
\[
 \sum_{j=0}^l u_j(c_j(i)-c_j(i'))\neq 0,\qquad 1\leq i\ne i'\leq D.
\]
Then $L=\sum_j u_j C_j$ has pairwise distinct mod $\mathfrak{n}$ eigenvalues $\bar\ell_1,\ldots,\bar\ell_D\in k'$. The ring $A'$ is complete local with residue field $k'$. Simple-root Hensel lifting gives unique $\ell_i\in A'$ with
\[
 \det(z-L)=\prod_{i=1}^D(z-\ell_i),\qquad
 \ell_i\equiv\bar\ell_i\bmod{\mathfrak nA'}, i=1,...,D.
\]
All differences $\ell_i-\ell_j$ are units in $A'$. Set
\begin{equation}\label{eqn:projector}
 P_i=\prod_{j\ne i}\frac{L-\ell_j\id}{\ell_i-\ell_j}.
\end{equation}
The Chinese remainder theorem and Cayley--Hamilton identity give $P_i^2=P_i$, $P_iP_j=0$ for $i\ne j$, and $\sum_iP_i=\id$. Each image is of finite rank and projective, and its reduction mod $\mathfrak{n}$ has dimension $1$; it is therefore free of rank one over the local ring $A'$. Every $C_j$ commutes with the $P_i$'s, hence acts by a scalar on each image. In particular
\[
 A'[C_0,\ldots,C_l]=\bigoplus_i A'P_i.
\]
This equality follows since the right-hand side contains every generator and the left-hand side contains every $P_i$.

The construction is independent of $L$. Indeed, two idempotents $e,f$ in this commutative algebra with equal reduction mod $\mathfrak{n}$ satisfy $(e-f)(1-e-f)=0$. The reduction of $1-e-f$ is $1-2\bar e$, whose square is $1$, so $1-e-f$ is an invertible endomorphism over the complete local ring. Thus $e=f$. This also proves uniqueness.
\end{proof}

\begin{theorem}[Joint Rees projectors]\label{thm:joint-rees}
After a finite separable extension $K'/K$, the primitive projectors of the generic localized small quantum ring have unique regular extensions
\begin{equation}\label{eqn:projector_lifted}
 P_{\boldsymbol\lambda}(s,Q)
       \in\operatorname{End}_{R'}(\mathscr H_R\otimes_RR'),
 \qquad R'=K'[\![s,Q]\!],
\end{equation}
labeled by tuples of partitions of total size $n$. These are mutually orthogonal rank-one idempotents and sum to the identity. Their complete $s=0$ specializations are
\begin{equation}\label{eqn:project_limit}
 P_{\boldsymbol\lambda}(0,Q)
       =\bigotimes_{\alpha=1}^rP^{(1)}_{\lambda^\alpha}(Q)
       \quad\hbox{on }\mathcal H_{\mathbf N},
 \qquad\mathbf N=(|\lambda^1|,\ldots,|\lambda^r|),
\end{equation}
and $0$ on other blocks. Here $P^{(1)}_\lambda$ is the primitive rank-one projector for the quantum cohomology of $\text{Hilb}^{|\lambda|}(\mathbb{C}^2)$ labeled by the partition $\lambda$. Subtracting the scalar in (\ref{eqn:scalar_divisor}) does not alter this projector. For charge $0$, $\M(1,0)=\text{pt}$ and the projector is the identity on the $1$-dimensional cohomology space.
\end{theorem}

\begin{proof}
The operators $C_0,\ldots,C_r$ constructed in Section \ref{sec:seperation_in_limit} are Rees regular and commuting, and (\ref{eqn:joint_char}) proves simple joint spectrum at
the closed point $(s,Q)=(0,0)$. Thus Lemma~\ref{lem:complete_local_projector} applies and produces (\ref{eqn:projector_lifted}). In fact the mod $\mathfrak{m}$ eigenvalues already
belong to $K$, so the projector construction itself does not need a field extension; we allow $K'$ for later normalization of idempotents.

Compare with the quantum cohomology ring over $K'(\!(s)\!)(\!(Q)\!)$. Every constructed projector commutes with $H\star_Q$, by (\ref{eqn:commuting_with_baxter}). The divisor operator $H\star_Q$ has simple spectrum  for generic values of parameters, so each $1$-dimensional image of the projectors is an eigenline of the divisor operator $H\star_Q$. Since the projectors are mutually orthogonal and exhaust the cohomology, they are exactly the spectral projectors of the divisor operator $H\star_Q$. Proposition~\ref{prop:baxter-divisor} identifies them with primitive multiplication projectors. The embedding
\[
 \operatorname{End}_{R'}(\mathscr H_R\otimes_RR')
 \hookrightarrow
 \operatorname{End}_{K'(\!(s)\!)(\!(Q)\!)}
  (\!(\mathscr H_R\otimes_RR')\otimes_{R'}K'(\!(s)\!)(\!(Q)\!))
\]
shows that their regular extensions are unique.

We next identify the entire $s=0$ limit. Apply Lemma \ref{lem:complete_local_projector} to the rank-one evaluated transfer operators over $K'[\![Q]\!]$. The partition separation proved in Section \ref{sec:seperation_in_limit} applies for every charge $N\le n$, giving regular rank-one projectors. By the
rank-one Baxter--divisor relation and generic simple spectrum, they are the primitive multiplication projectors $P^{(1)}_\lambda(Q)$. This argument also proves their regularity at $Q=0$ without requiring the classical divisor to have distinct eigenvalues.

After $s=0$, formula (\ref{eqn:transfer_matrix_expand}) makes the $C_\alpha$ factorwise rank-one operators, and $C_0$ is $B$. Consequently the tensor products on the right of
(\ref{eqn:project_limit}) decompose each block into joint eigenlines. 
A separating linear combination of this limit family still has pairwise distinct eigenvalues modulo $Q$, so its unique Hensel lifted projectors are both these tensor-product projectors and the specializations of
$P_{\boldsymbol\lambda}(s,Q)$. This proves (\ref{eqn:project_limit}) as an identity over $K'[\![Q]\!]$.

All constructions of the projectors already take place over $K$. If $K'$ is enlarged to choose square roots of their unit norms, the projectors, their traces, and the symmetric sums used below keep their descent to $K$.
\end{proof}

\subsection{The regular pointwise quantum product}

We absorb an allowed finite separable coefficient extension into the notation $K,R$.

\begin{corollary}[Pointwise quantum product extension]
\label{cor:product-extension}
Put
\[
 e_i=P_i(\id),\qquad
 \widehat\eta=t^c\eta_{\M(r,n)},\qquad
 g_i=\widehat\eta(e_i,e_i),\qquad c=2n(r-1).
\]
Every $g_i$ is a unit of $R$. The pointwise small quantum product extends over the Rees lattice by
\begin{equation}\label{eqn:quantum_product_extension}
 x\circ y=\sum_i
 \frac{\widehat\eta(x,e_i)\widehat\eta(y,e_i)}{g_i^2}\,e_i.
\end{equation}
At $s=0$, it is, on each color-charge block, the tensor product of the quantum products of the rank-one instanton moduli, and the product between distinct blocks is zero.
\end{corollary}

\begin{proof}
The unit and the projectors are regular, so $e_i\in\mathscr H_R$. Their reductions at $(s,Q)=(0,0)$ are the classical fixed-point idempotents in the indicated tensor products. They form a basis by Nakayama's lemma. The rescaled metric is regular by Proposition~\ref{prop:stable-lattice}; in the specialization $s=Q=0$, its value on $e_i$ is a nonzero block constant times a product of nonzero classical rank-one idempotent norms. Hence $g_i\neq 0\bmod\mathfrak m$, proving that $g_i$ is a unit.

Over the field $K(\!(s)\!)(\!(Q)\!)$, the $e_i$ are orthogonal primitive idempotents. Their metric orthogonality extends over $R$ by injectivity of $R\hookrightarrow K(\!(s)\!)(\!(Q)\!)$. Thus $x=\sum_i\widehat\eta(x,e_i)g_i^{-1}e_i$, and multiplication over the field $K(\!(s)\!)(\!(Q)\!)$ gives (\ref{eqn:quantum_product_extension}). Its coefficients are regular, so it defines the required extension; associativity and the unit identities also extend from the generic field. Finally, (\ref{eqn:Id_expand}), (\ref{eqn:metric_expand}), and (\ref{eqn:project_limit}) identify its limit with the asserted tensor-product algebra.
\end{proof}

\section{Four-point correlators}\label{sec:four-point-rees}

In Section \ref{sec:Frob}, we construct the Rees lattice in $QH_\T^*(\M(r,n))_\loc$ using stable envelopes, and in Section \ref{sec:projector} we extend its idempotents to $s=0$. The purpose of this section is to study the order in $s$ of one flat derivative of the quantum product, represented by a genus-$0$, $4$-point correlator. We need this to analyze $\mathsf{R}$-matrix coefficients in Section \ref{sec:limit_genus1}.

Recall the notations 
$a_\alpha=t b_\alpha, s=t^{-1}, K=\mathbb Q(\epsilon_1,\epsilon_2,b_1,\ldots,b_r),R=K[\![s,Q]\!], c=2n(r-1)$, and the Rees lattice $\mathscr{H}_R$ in (\ref{eqn:Rees_lattice}). Using the normalized stable-envelope map \(\Phi_t\) in (\ref{eqn:stab_isom}), define
\[
 \eta_R:=\Phi_t^*\eta_{\M(r,n)}.
\]
By Proposition~\ref{prop:stable-lattice}, 
\(\widehat\eta:=t^c\eta_R\) is a nondegenerate \(R\)-valued pairing. Equivalently,
\begin{equation}\label{eqn:pairing_rescaled}
 \eta_R^{-1}=t^c\widehat\eta^{-1}.
\end{equation}

For \(Q\)-independent vectors \(\gamma_1,\ldots,\gamma_4\in\mathscr H_R\), define the genus $0$, $4$-point correlator 
\begin{equation}\label{eqn:4-pt-series}
 \mathsf C_4(\gamma_1,\gamma_2,\gamma_3,\gamma_4)
 :=\sum_{d\geq0}(-1)^{rd}Q^d
 \int_{[\Mbar_{0,4}(\M(r,n),d)]^{\mathrm{vir}}}\prod_{k=1}^4\operatorname{ev}_k^*(\Phi_t\gamma_k).
\end{equation}
The main result of this section, Proposition \ref{prop:fourpoint}, is the following identification of the leading order of $\mathsf{C}_4$:
\begin{equation}\label{eqn:4pt-series-regular}
 t^c\mathsf C_4(\gamma_1,\gamma_2,\gamma_3,\gamma_4)\in R.
\end{equation}
After raising one index with \(\eta_R^{-1}\), the factor \(t^c\) in (\ref{eqn:pairing_rescaled}) cancels the factor \(t^{-c}\) in (\ref{eqn:4pt-series-regular}), giving a Rees regular endomorphism. We prove this by separating the deformation complex of the map from the deformation complex of its marked domain.

\subsection{The moving index of the map deformation complex}
Consider the effective framing torus $\mathsf{A}$ in (\ref{eqn:eff_frame_torus}). The fixed locus $\M(r,n)^{\mathsf{A}}$ given by (\ref{eqn:framed_torus_fixed_loci}) has pure dimension \(2n\). Since \(\dim \M(r,n)=2rn\), the codimension of every fixed component is $c=2rn-2n=2n(r-1)$. The following index lemma explains why the same number occurs for every fixed genus-$0$ stable map, independently of its graph.

\begin{lemma}[Symplectic moving-index]
\label{lem:rooted-index}
Let \(X\) be a smooth symplectic variety equipped with an algebraic $\mathbb{C}^*$-action preserving the symplectic form. For a connected $\mathbb{C}^*$-fixed genus-$0$ stable map $u:C\rightarrow X$, suppose that all components of \(X^{\mathbb C^*}\) meeting the image of $u$ have the same dimension \(\mathsf{f}\). Then the virtual ranks satisfy
\[
 \operatorname{vrk}\bigl(R\Gamma(C,u^*T_X)^{\mathrm{fix}}\bigr)=\mathsf{f} 
 \quad \text{and hence }
 \operatorname{vrk}\bigl(R\Gamma(C,u^*T_X)^{\mathrm{mov}}\bigr)=\dim_\CC X-\mathsf{f}.
\]
\end{lemma}

\begin{proof}
Consider the normalization of $C$:
\begin{equation}\label{eqn:normalization_curve}
 \nu:\widetilde C=\coprod_v C_v\longrightarrow C,
\end{equation}
 let \(q_e\) be the nodes. The normalization exact sequence gives an equality 
\begin{equation}\label{eqn:class_RGamma}
 [R\Gamma(C,u^*T_X)] =\sum_v[R\Gamma(C_v,u_v^*T_X)] -\sum_e[T_{u(q_e)}X]
\end{equation}
in equivariant \(K\)-theory. We compute the virtual rank of the $\CC^*$-fixed part of each term.

If the torus acts trivially on \(C_v\), equivariance forces the image of \(u_v\) to lie in a fixed component \(F\subset X^{\mathbb C^*}\). The weight-zero subbundle of \(u_v^*T_X\) is \(u_v^*T_F\).  Since the torus preserves the symplectic form, \(F\) is symplectic; in particular,
\(\operatorname{rk}T_F=\mathsf{f}\) and \(c_1(T_F)=0\). Riemann-Roch on \(C_v\cong\mathbb P^1\) therefore gives
\[
 \operatorname{vrk}R\Gamma(C_v,u_v^*T_X)^{\mathrm{fix}} =\chi(C_v,u_v^*T_F)=\mathsf{f}.
\]

Suppose that the torus acts nontrivially on \(C_v\cong\mathbb P^1\). Its two fixed points map to points \(x_0,x_\infty\in X^{\mathbb C^*}\). We decompose \(u_v^*T_X\), in a way compatible with the symplectic form, into a sum of pairs of line bundles
\[
 L\oplus L^\vee
\]
in equivariant \(K\)-theory. Let \(\kappa\neq0\) be the tangent character of \(C_v\) at \(0\), and let \(a,b\) be the fiber characters of \(L\) at \(0,\infty\), respectively. Then localization gives
\[
 a-b=d\kappa
\]
where \(d=\deg L\). Denote by \(\iota_\kappa(a,b)\) the multiplicity of the trivial character in
\(\chi(\mathbb P^1,L)\). Explicitly,
\begin{equation}\label{eqn:char_multiplicity}
 \iota_\kappa(a,b)=
 \begin{cases}
  1,&d\geq0\text{ and }0\in\{a,a-\kappa,\ldots,b\},\\
 -1,&d\leq-2\text{ and }0\in\{a+\kappa,\ldots,b-\kappa\},\\
  0,&\text{otherwise}.
 \end{cases}
\end{equation}
The first line lists the \(H^0\)-weights; the second lists the interior \(H^1\)-weights with their negative sign. Comparing these for \(L\) and \(L^\vee\) gives 
\begin{equation}\label{eqn:char_multiplicity2}
 \iota_\kappa(a,b)+\iota_\kappa(-a,-b) =\delta_{a,0}+\delta_{b,0}.
\end{equation}
Summing over all pairs of line bundles gives
\begin{equation}\label{eqn:fixed_vrk}
 \operatorname{vrk}R\Gamma(C_v,u_v^*T_X)^{\mathrm{fix}}
 =\frac12\bigl(\dim T_{x_0}X^{\mathbb C^*} +\dim T_{x_\infty}X^{\mathbb C^*}\bigr)
 =\mathsf{f}.
\end{equation}

If \(C\) has \(N\) irreducible components, then its dual graph is a tree with \(N-1\) nodes. For every component term in (\ref{eqn:class_RGamma}), its $\mathbb{C}^*$-fixed part has virtual rank \(\mathsf{f}\), and the fixed part of the tangent space at a node has dimension \(\mathsf{f}\). Thus
\[
 \operatorname{vrk}R\Gamma(C,u^*T_X)^{\mathrm{fix}} =N\mathsf{f}-(N-1)\mathsf{f}=\mathsf{f}.
\]
Finally, as \(C\) is of genus $0$ and \(c_1(T_X)=0\), Riemann-Roch gives $\operatorname{vrk}R\Gamma(C,u^*T_X)=\dim X.$
\end{proof}

\begin{corollary}[Framing moving rank]
\label{cor:framing-index}
For every full-\(\T\)-fixed genus-$0$ stable map to \(\M(r,n)\), the virtual rank of the part of
\[
 R\Gamma(C,u^*T_{\M(r,n)})
\]
with nonzero framing-torus character (i.e. the moving part) is
\[
 c=2n(r-1).
\]
\end{corollary}

\begin{proof}
A generic one-parameter subgroup \(\lambda:\mathbb C^*\to \mathsf{A}\) satisfies
\[
 \M(r,n)^\lambda=\M(r,n)^{\mathsf{A}}.
\]
Choose \(\lambda\) additionally so that no nonzero \(\mathsf{A}\)-character occurring in the fixed-map deformation complex restricts trivially to \(\lambda\). This is possible because only finitely many characters occur on each fixed component of the stack of stable maps. The torus \(\lambda(\mathbb C^*)\) preserves the symplectic form, and all components of its fixed locus have dimension \(2n\). Lemma~\ref{lem:rooted-index} therefore gives moving virtual rank
\[
 2rn-2n=2n(r-1).
\]
By the choice of \(\lambda\), this moving part is exactly the sum of subspaces with nonzero framing characters.
\end{proof}

\subsection{The correction from the four-pointed domain}

Corollary \ref{cor:framing-index} concerns deformations of a map with fixed domain. The absolute virtual tangent complex also contains
\[
 \operatorname{Ext}^1(\Omega_C(D),\mathcal O_C)
 -\operatorname{Ext}^0(\Omega_C(D),\mathcal O_C),
\]
where \(D\) is the sum of the four markings. We show that the virtual rank of its moving part is nonnegative.

\begin{lemma}[Nonnegative domain correction]
\label{lem:domain-correction}
Let \(u:(C,D)\to X\) be a connected $4$-pointed genus-$0$ stable map fixed by a $1$-dimensional torus. Then
\[
 d_C:=\operatorname{vrk}
 \left(\operatorname{Ext}^1(\Omega_C(D),\mathcal O_C)
 -\operatorname{Ext}^0(\Omega_C(D),\mathcal O_C)\right)^{\mathrm{mov}}
 \geq0.
\]
\end{lemma}

\begin{proof}
Consider the normalization of \(C\) as in (\ref{eqn:normalization_curve}). If \(D_v\) denotes the markings and preimages of nodes on an irreducible component \(C_v\), the standard local-to-global sequence for the deformation theory of a pointed nodal curve gives the following identity in equivariant \(K\)-theory:
\[
 \begin{aligned}
 &[\operatorname{Ext}^1(\Omega_C(D),\mathcal O_C)
   -\operatorname{Ext}^0(\Omega_C(D),\mathcal O_C)]\\
 &\qquad =\sum_v
   [H^1(C_v,T_{C_v}(-D_v))-H^0(C_v,T_{C_v}(-D_v))]
   +\sum_e[T_{q_e'}C_{v'}\otimes T_{q_e''}C_{v''}],
 \end{aligned}
\]
where the last sum consists of the node-smoothing contributions. Thus it suffices to count moving infinitesimal automorphisms and moving node smoothings.

A component on which the compensated torus action is trivial contributes only weight-zero component deformations and infinitesimal automorphisms. A component with nontrivial action is a \(\mathbb P^1\) on which every marking and every branch node is one of the two torus-fixed points. If it carries \(s_v\) special points, then \(s_v\leq2\), so
\[
 H^1\bigl(C_v,T_{C_v}(-D_v)\bigr)=0,
\]
and the moving part of \(H^0(C_v,T_{C_v}(-D_v))\) has dimension \(2-s_v\).

The components with nontrivial actions form a disjoint union of chains. Fix one such chain \(\Gamma\). Let
\[
 N_\Gamma=\#\{\text{components of }\Gamma\},\qquad
 m_\Gamma=\#\{\text{markings on }\Gamma\},
\]
and let \(e_\Gamma\) be the number of nodes joining \(\Gamma\) to components with trivial torus action. The chain has \(N_\Gamma-1\) internal nodes, so summing the numbers of special points over its components gives
\[
 \sum_{v\in\Gamma}s_v =m_\Gamma+e_\Gamma+2(N_\Gamma-1).
\]
Consequently, the moving infinitesimal-automorphism dimension contributed by \(\Gamma\) is
\begin{equation}\label{eqn:moving_auto_rk_Gamma}
 \sum_{v\in\Gamma}(2-s_v)=2N_\Gamma-\sum_{v\in\Gamma}s_v=2-m_\Gamma-e_\Gamma.
\end{equation}

Every attaching node counted by \(e_\Gamma\) has nonzero smoothing weight: such a node is the intersection of two branches one of which has a nontrivial torus action and the other has the trivial torus action. Let \(e_{m,\Gamma}\geq0\) be the number of internal nodes of \(\Gamma\) whose smoothing character is nonzero. The total virtual rank of the moving part contributed by \(\Gamma\) is therefore
\begin{equation}\label{eqn:moving_vrk_Gamma}
 e_\Gamma+e_{m,\Gamma}-(2-m_\Gamma-e_\Gamma) =m_\Gamma+2e_\Gamma+e_{m,\Gamma}-2.
\end{equation}
It remains to show that (\ref{eqn:moving_vrk_Gamma}) is nonnegative. A moving chain has only two free fixed points, one at each end. Hence it can meet the rest of the curve only at its ends, so \(e_\Gamma\leq2\). The case \(e_\Gamma=0\) cannot occur: connectedness would make \(\Gamma\) the whole curve, but then all four markings would have to lie on a moving chain, which has room for at most two distinct marked points. Thus \(e_\Gamma=1\) or \(2\). If \(e_\Gamma=1\), expression (\ref{eqn:moving_vrk_Gamma}) is \(m_\Gamma+e_{m,\Gamma}\geq0\); if \(e_\Gamma=2\), it is
\(m_\Gamma+e_{m,\Gamma}+2\geq0\). Summing over the moving chains proves \(d_C\geq0\). If there are no moving components, then \(d_C=0\).
\end{proof}

\subsection{The raised-index correlator}

We now combine the two virtual rank calculations. The relative map deformation complex has nonzero framing-character virtual rank \(c\), while the marked-domain complex adds a nonnegative rank \(d_C\). Characters with zero framing part depend only on \(\epsilon_1,\epsilon_2\) and do not affect the order in \(t\).

\begin{proposition}[Four-point Rees regularity]
\label{prop:fourpoint}
Let \(\gamma_1,\ldots,\gamma_4\in\mathscr H_R\) be independent of \(Q\). Then, coefficientwise in \(Q\),
\begin{equation}\label{eqn:4pt-series-regular2}
 t^c\mathsf C_4(\gamma_1,\gamma_2,\gamma_3,\gamma_4) \in R.
\end{equation}
For \(Q\)-independent \(\gamma_1,\gamma_2\in\mathscr H_R\), define \(\mathcal D_{\gamma_1,\gamma_2}\) (over the generic point of $\text{Spec}R$) by requiring
\begin{equation}\label{eqn:operator_D}
 \eta_R\bigl(\mathcal D_{\gamma_1,\gamma_2}x,y\bigr) =\mathsf C_4(\gamma_1,\gamma_2,x,y)
\end{equation}
for \(Q\)-independent \(x,y\) in a basis of \(\mathscr H_R\), and then extend \(R\)-linearly and \(Q\)-adically. Then
\[
 \mathcal D_{\gamma_1,\gamma_2}\in\operatorname{End}_R(\mathscr H_R).
\]
\end{proposition}

\begin{proof}
Fix a coefficient \(Q^d\) of (\ref{eqn:4-pt-series}). Apply the virtual localization with respect to \(\T\). Let \(F\) be a connected component of the fixed locus of the stack of stable maps. Choose a generic one-parameter subgroup of the framing torus as in Corollary~\ref{cor:framing-index}. The nonzero framing-character part of the relative map deformation complex has virtual rank \(c\). By
Lemma~\ref{lem:domain-correction}, the marked-domain deformation complex contributes an additional \(d_F\geq0\) to the virtual rank of the moving part. Therefore the part of the absolute virtual normal bundle carrying nonzero framing character has virtual rank $c+d_F$.

After the substitution \(a_\alpha=t b_\alpha\), every nonzero framing weight has the form
\[
 t\ell(b_1,\ldots,b_r)+\nu(\epsilon_1,\epsilon_2) =t\ell(b_1,...,b_r)\bigl(1+O(s)\bigr),
\]
where \(\ell(b_1,...,b_r)\in K^\times\). If the moving part of a virtual representation has virtual rank \(c+d_F\), the inverse of its Euler class is therefore contained in
\[
 t^{-(c+d_F)}\cdot K[\![s]\!]^\times =t^{-c}s^{d_F}\cdot K[\![s]\!]^\times.
\]
The remaining full moving \(\T\)-weights have zero framing character and are independent of \(t\). The restrictions of \(\Phi_t\gamma_i\) to the full-torus fixed loci are regular in \(s\): this is exactly the triangular degree estimate used to define the normalized stable-envelope lattice in Section \ref{sec:Frob}. Hence the contribution of \(F\) to $\mathsf{C}_4$ lies in
\[
 t^{-c}s^{d_F}K[\![s]\!]\subseteq t^{-c}K[\![s]\!].
\]
There are finitely many fixed components in each degree, and their equivariant integrations preserve this $t$-order. Summing them proves (\ref{eqn:4pt-series-regular2}), coefficientwise in \(Q\).

Write the correlator in (\ref{eqn:operator_D}) as
\[
 \mathsf C_4(\gamma_1,\gamma_2,-,-) =t^{-c}\widehat{\mathsf C}_{\gamma_1,\gamma_2},
\]
where \(\widehat{\mathsf C}_{\gamma_1,\gamma_2}\) is an \(R\)-valued bilinear form by (\ref{eqn:4pt-series-regular2}). Using (\ref{eqn:pairing_rescaled}), raising the first index gives
\[
 \mathcal D_{\gamma_1,\gamma_2}
 =\eta_R^{-1}\circ \mathsf C_4(\gamma_1,\gamma_2,-,-)
 =t^c\widehat\eta^{-1}\circ t^{-c}\widehat{\mathsf C}_{\gamma_1,\gamma_2}
 =\widehat\eta^{-1}\circ \widehat{\mathsf C}_{\gamma_1,\gamma_2}.
\]
Both factors on the right are defined over \(R\), so \(\mathcal D_{\gamma_1,\gamma_2}\) preserves \(\mathscr H_R\).
\end{proof}

\begin{corollary}[Regularity of the flat derivative of multiplication]
\label{cor:fourpoint-derivative}
Let \(\gamma\in\mathscr H_R\) be \(Q\)-independent and let \(B_\gamma(\mathbf t)\) denote operator of multiplication by \(\gamma\) in the localized big quantum product, with \(\mathbf t\) denoting the flat coordinates. On the locus $\mathbf{t}=0$ corresponding to small quantum cohomology, we have
\begin{equation}\label{eqn:flat_product_regular}
 \left.(\partial_\gamma B_\gamma)\right|_{\mathbf t=0}
 =\mathcal D_{\gamma,\gamma}\in\operatorname{End}_R(\mathscr H_R).
\end{equation}
\end{corollary}

\begin{proof}
For constant flat vectors \(x,y\), differentiating the primary genus-$0$ potential four times gives
\[
 \eta_R\!\left(\left.(\partial_\gamma B_\gamma)\right|_{\mathbf t=0}x,y\right)
 =\sum_{d\geq0}(-1)^{rd}Q^d\left\langle\gamma,\gamma,x,y\right\rangle^{\M(r,n)}_{0,4,d}
 =\mathsf C_4(\gamma,\gamma,x,y).
\]
By the defining identity (\ref{eqn:operator_D}), the derivative is \(\mathcal D_{\gamma,\gamma}\), which is Rees regular by Proposition~\ref{prop:fourpoint}.
\end{proof}

\section{Specialization of Givental's formula}\label{sec:limit_genus1}

This section contains two results. We first prove that $\mathsf{R}_1$ for the Gromov-Witten theory of $\M(r,n)$ is Rees regular. Then we prove that the two terms in Givental's genus-one formula specialize, at $s=0$, to the corresponding terms of the tensor product of the Gromov-Witten theories of the rank-one cases $\text{Hilb}^{N_\alpha}(\mathbb{C}^2)$. The two assertions use different inputs: the four-point estimate is essential for $\mathsf{R}_1$, whereas the first term uses only the projectors and the rescaled metric.

Throughout this section, we work after enlarging, if necessary, the finite separable extension allowed in Theorem~\ref{thm:joint-rees} and omit that from the notations. All matrix coefficients and generating series are coefficientwise formal in $Q$. Let
\[
 V=Q\frac{\partial}{\partial Q}.
\]
$V$ acts trivially on $s$, $t=s^{-1}$, and $K$. Let $\partial_H$ denote the constant flat divisor vector field on the big quantum cohomology. The divisor equation identifies\footnote{See \cite[Section 2.3]{pt} for a detailed argument in the case of $\text{Hilb}^n(\mathbb{C}^2)$. The same argument applies here.} these two derivations on ingredients of the Frobenius structure. 

\subsection{$\mathsf{R}$-matrix construction in the Rees frame}

Let $P_i\in\operatorname{End}_R(\mathscr H_R)$ be the primitive projectors of Theorem~\ref{thm:joint-rees}, and put $e_i=P_i(\id)$. By Corollary~\ref{cor:product-extension}, the $e_i$ are Rees regular primitive idempotents and their norms for the rescaled metric $\widehat\eta=t^c\eta$, are units:
\[
 \widehat g_i=\widehat\eta(e_i,e_i)\in R^\times.
\]
After extending the coefficient field, choose square roots $\widehat g_i^{1/2}\in R^\times$ and introduce the normalized idempotents:
\[
 f_i=\widehat g_i^{-1/2}e_i.
\]
The vectors $f_i$ form an $\widehat\eta$-orthonormal basis of the Rees module $\mathscr{H}_R$. Let $\Psi$ denote the matrix whose columns are the $f_i$ in a chosen flat basis of $\mathscr H_R$. Then
\[
 \Psi\in\operatorname{GL}(\mathscr H_R, R),\qquad \Psi^T\widehat\eta\Psi=\id.
\]

The frame normalized for the original metric $\eta=t^{-c}\widehat\eta$ is $t^{c/2}\Psi$. This differs from $\Psi$ by one common scalar, independent of all Frobenius-manifold variables and of $Q$. Consequently, the matrices
\[
 \Psi^{-1}d\Psi, \qquad \Psi^{-1}(\gamma\star)\Psi, \qquad \mathsf{R}_1
\]
are unchanged if we work with the regular $\widehat\eta$-normalized frame instead of the original $\eta$-normalized frame. We shall therefore use $\Psi$ throughout.

Virtual localization defines the localized primary genus-$0$ potential and hence the formal big quantum product $\star_{\mathbf t}$ in constant flat coordinates $\mathbf t$. Only its germ at the locus $\mathbf t=0$ is used. For a constant flat class $\gamma$, write
\[
 B_\gamma(\mathbf t)=\gamma\star_{\mathbf t}(-), \qquad B_\gamma=B_\gamma(\mathbf{t}=0)=\gamma\star.
\]
The fourth derivative of the primary genus $0$ potential gives, for constant flat vectors $x,y$,
\begin{equation}\label{eqn:operator_derivative}
 \eta\!\left(\left.(\partial_\gamma B_\gamma)\right|_{\mathbf t=0}x,y\right)
 =\sum_{d\ge0}(-1)^{rd}Q^d\left\langle\gamma,\gamma,x,y\right\rangle^{\M(r,n)}_{0,4,d}.
\end{equation}
By Corollary~\ref{cor:fourpoint-derivative}, the endomorphism on the left is Rees regular whenever $\gamma$ is $Q$-independent.

The small quantum product is generically semisimple. Let $u_i$ be local canonical coordinates on the big germ, and let $d$ denote its differential in the flat variables. Put
\[
 U=\operatorname{diag}(u_i),\qquad W=\Psi^{-1}d\Psi,
 \qquad \rho_i=du_i(\partial_H)\big|_{\mathbf t=0}.
\]
Thus $\rho_i$ is the eigenvalue of $H\star_Q$. After absorbing the flat divisor coordinate into the Novikov variable, the canonical one-forms along the logarithmic divisor trajectory are
\begin{equation}\label{eqn:logarithmic-canonical-forms}
 du_i=\rho_i(Q;s)\frac{dQ}{Q}.
\end{equation}
This denotes the divisor trajectory, not ordinary restriction of a chosen $Q$-power-series coordinate to $\mathbf t=0$. If
$\rho_i(Q;s)=\rho_i(0;s)+\sum_{d\geq1}\rho_{i,d}(s)Q^d$, a primitive is
\[
 u_i^{\log}(Q;s)=\rho_i(0;s)\log Q
       +\sum_{d\geq1}\frac{\rho_{i,d}(s)}{d}Q^d+\text{constant}.
\]
In particular, $V(u_i^{\log})=\rho_i$ includes the generally nonzero classical divisor restriction $\rho_i(0;s)$. Only the one-forms, and never a choice of these logarithmic primitives, enter the argument.

Consider the $\mathsf{R}$-matrix of the $\T$-equivariant Gromov-Witten theory of $\M(r,n)$:
\[
 \mathsf{R}(z)=1+\mathsf{R}_1z+\mathsf{R}_2z^2+\cdots.
\]
By \cite[Proposition~1.1]{g1}, $\mathsf{R}_i$ are determined recursively: the first two recursions are
\begin{equation}\label{eqn:R-recursion}
 W=[dU,\mathsf{R}_1],
 \qquad
 d\mathsf{R}_1+W\mathsf{R}_1=[dU,\mathsf{R}_2].
\end{equation}
Since the frame is orthonormal, $W^T=-W$. The symplectic condition 
$\mathsf{R}(z)\mathsf{R}^T(-z)=1$ gives $\mathsf{R}_1^T=\mathsf{R}_1$. In particular, for $i\ne j$,
\[
 W_{ij}=d(u_i-u_j)(\mathsf{R}_1)_{ij}.
\]

For a flat vector field $X$, Givental's genus-one formula \cite{g1} gives the following expression for the generating function $\omega_1(X)$ of genus $1$ Gromov-Witten invariants with insertion $X$:
\begin{equation}\label{eqn:Givental_formula_general}
 \omega_1(X)=\frac1{48}X\left(\log\prod_i\Delta_i\right)+\frac12\sum_i(\mathsf{R}_1)_{ii}X(u_i),
 \qquad
 \Delta_i=\eta(e_i,e_i)^{-1}.
\end{equation}
Here $X(\log f)$ means $X(f)/f$; no choice of a formal logarithm is involved. Taking $X=\partial_H$ gives the one-divisor series $\F_{r,n}$. By divisor equation, its metric derivative can be written using $V$, while $du_i(\partial_H)=\rho_i$. Its positive-$Q$-degree part is $\F^+_{r,n}$.

The diagonal recursion for $\mathsf{R}_1$ will be integrated in the $Q$ variable. Its initial value at $Q=0$ is fixed by degree-$0$, genus $1$ invariants as the Lemma below, see the discussion around \cite[Theorem 2.1]{g1}.

\begin{lemma}[Constant-map initial condition]
\label{lem:gamma-initial}
For every $r\geq1$ and $n\geq1$, at $Q=0$ in the normalized idempotent basis corresponding to $\T$-fixed points, we have
\begin{equation}\label{eq:gamma-initial}
 (\mathsf{R}_1(Q=0))_{\mathfrak p\mathfrak p}
 =-\frac1{12}\sum_{\chi\in\operatorname{Wt}(T_{\mathfrak p}\M(r,n))}\frac1\chi, \quad \mathfrak{p}\in \M(r,n)^\T.
\end{equation}
\end{lemma}

The initial value in Lemma~\ref{lem:gamma-initial} is Rees regular. Indeed, same-color tangent weights are independent of the large framing parameter, whereas every cross-color weight has the form
\[
 t(b_\beta-b_\alpha)+O(1),\qquad \alpha\ne\beta,
\]
so its reciprocal is $O(s)$.

\subsection{Regularity of $\mathsf{R}$-matrix coefficients}

The divisor operator $H\star_Q$ is not suitable for solving the first recursion near $s=0$ because its eigenvalues may collide inside a color-charge block $\mathcal{H}_{\mathbf{N}}$. We instead choose a constant flat class whose multiplication eigenvalues are already distinct at $s=0$.

\begin{lemma}[A separating flat multiplication operator]
There exists a $Q$-independent flat class $\gamma\in\mathscr H_R$ such that, if
\begin{equation}\label{eqn:diagonalize_B_gamma}
 B_\gamma=\gamma\star =\Psi\operatorname{diag}(\beta_1,\ldots,\beta_D)\Psi^{-1},
\end{equation}
then
\begin{equation}\label{eqn:distinct_eigenvalues}
 \beta_i-\beta_j\in R^\times, \quad i\ne j.
\end{equation}
\end{lemma}

\begin{proof}
At the closed point $(s,Q)=(0,0)$, the Frobenius algebra is the direct product of the $1$-dimensional subalgebras generated by the reductions $\bar e_i$. Choose pairwise distinct constants $c_i\in K$ and put
\begin{equation}\label{eqn:combination_idempotent}
 \bar\gamma=\sum_i c_i\bar e_i.
\end{equation}
The state space $\mathscr{H}_R$ is a constant free module in flat cohomology coordinates, so the reduction map 
\(
 \mathscr H_R\to\mathscr H_R/\mathfrak m\mathscr H_R
\)
admits a $K$-linear section. Choose a lift $\gamma$ with coefficients in $K$; in particular, it is independent of both $Q$ and the flat coordinates.

Corollary~\ref{cor:product-extension} implies that multiplication by $\gamma$ is an element of $\operatorname{End}_R(\mathscr H_R)$. Since $P_i$ is a multiplication projector, $B_\gamma$ preserves the free rank-one module $\operatorname{Im}P_i$ and acts there by a scalar $\beta_i\in R$. Reduction
modulo $\mathfrak m$ gives $\beta_i\equiv c_i$.  Therefore
\[
 \beta_i-\beta_j\equiv c_i-c_j\not=0 \bmod{\mathfrak m},
\]
which proves that every difference is a unit of the local ring $R$.
\end{proof}

\begin{theorem}[Regularity of Gromov-Witten $\mathsf{R}_1$]
\label{thm:R1-regular}
After the finite separable extension in Theorem~\ref{thm:joint-rees}, we have 
\begin{equation}\label{eqn:regular_R1}
 \mathsf{R}_1\in\operatorname{End}_R(\mathscr H_R).
\end{equation}
\end{theorem}

\begin{proof}
We prove regularity for off-diagonal and diagonal entries separately.

\noindent{\bf Step 1:} {\it control the off-diagonal entries using the four-point correlator.} 

Choose $\gamma$ as in the preceding lemma and let $\partial_\gamma$ be the corresponding constant flat vector field on the formal germ of the big quantum cohomology Frobenius structure. Set
\begin{equation}\label{eqn:W-matrix}
 W_\gamma=\Psi^{-1}\partial_\gamma\Psi.
\end{equation}
The eigenvalue of multiplication by $\gamma$ on the $i$-th eigenline is $du_i(\partial_\gamma)=\beta_i$. Evaluating the first recursion in (\ref{eqn:R-recursion}) on $\partial_\gamma$ therefore gives, for $i\ne j$,
\begin{equation}\label{eqn:R-recursion_gamma}
 (W_\gamma)_{ij} =(\beta_i-\beta_j)(\mathsf{R}_1)_{ij}.
\end{equation}

Differentiating the diagonalization (\ref{eqn:diagonalize_B_gamma}) of $B_\gamma$ and then conjugating by $\Psi$, we get
\[
 \Psi^{-1}(\partial_\gamma B_\gamma)\Psi
 =\partial_\gamma\operatorname{diag}(\beta_i)+[W_\gamma,\operatorname{diag}(\beta_i)].
\]
Taking $(i,j)$-entries give
\begin{equation}\label{eqn:R-recursion_ij}
 \begin{aligned}
 \bigl(\Psi^{-1}(\partial_\gamma B_\gamma)\Psi\bigr)_{ij}
 =(\beta_j-\beta_i)(W_\gamma)_{ij}=-(\beta_i-\beta_j)^2(\mathsf{R}_1)_{ij}.
 \end{aligned}
\end{equation}
The derivative $\partial_\gamma\operatorname{diag}(\beta_i)$ is diagonal and therefore does not enter this identity.

By (\ref{eqn:operator_derivative}) and Corollary~\ref{cor:fourpoint-derivative},
$\partial_\gamma B_\gamma$ preserves $\mathscr H_R$. The frame $\Psi$ and its inverse are Rees regular, and every difference $\beta_i-\beta_j$ is a unit by (\ref{eqn:distinct_eigenvalues}). Equation (\ref{eqn:R-recursion_ij}) therefore proves that all off-diagonal entries of $\mathsf{R}_1$ are in $R$.

\noindent{\bf Step 2:} {\it derive the diagonal recursion.}

Let
\[
 W_V=\Psi^{-1}V\Psi, \qquad
 \rho_i=du_i(\partial_H)\big|_{\mathbf t=0}.
\]
By divisor equation, $W_V=\Psi^{-1}\partial_H\Psi$ on the locus $\mathbf{t}=0$. Contracting the first recursion with the flat divisor field $\partial_H$ therefore gives
\begin{equation}\label{eqn:R-recursion_V}
 (W_V)_{ij}=(\rho_i-\rho_j)(\mathsf{R}_1)_{ij}, \quad i\ne j.
\end{equation}
Contract the second recursion in (\ref{eqn:R-recursion}) with $\partial_H$ and take its $(i,i)$-entry. The commutator $[\operatorname{diag}(\rho_i),\mathsf R_2]$ has zero diagonal. Divisor equation implies that
$\partial_H(\mathsf R_1)_{ii}=V(\mathsf R_1)_{ii}$. Hence we obtain
\[
 V(\mathsf{R}_1)_{ii}=-(W_V\mathsf{R}_1)_{ii}.
\]
Using $W_V^T=-W_V$, $\mathsf{R}_1^T=\mathsf{R}_1$, and (\ref{eqn:R-recursion_V}), we obtain
\begin{equation}\label{eqn:R-diagonal1}
 V(\mathsf{R}_1)_{ii}=-\sum_{j\ne i}(\rho_i-\rho_j)(\mathsf{R}_1)_{ij}^{\,2},
\end{equation}
see also the proof of \cite[Proposition 1.1]{g1}. This fixes both the sign and the absence of a factor of $2$.

\noindent{\bf Step 3:} {\it control the cross-block terms in the diagonal recursion.}

The operator $s(H\star_Q)$ is Rees regular and acts on the free rank-one module $\operatorname{Im}P_i$ by the scalar $s\rho_i$. If the label $i$ belongs to the color-charge block $\mathcal{H}_{\mathbf N}$, then Proposition~\ref{prop:stable-lattice} gives
\[
 s\rho_i\equiv\nu_{\mathbf N}\bmod{s}.
\]
Hence there is a uniquely determined $\widetilde\rho_i\in R$ such that
\begin{equation}\label{eqn:eigenvalue_expand}
 \rho_i=s^{-1}\nu_{\mathbf N}+\widetilde\rho_i.
\end{equation}

If $i$ and $j$ lie in the same color-charge block, the leading terms in (\ref{eqn:eigenvalue_expand}) cancel, so
\[
 \rho_i-\rho_j\in R.
\]
Together with the regularity of off-diagonal entries $(\mathsf{R}_1)_{ij}$ already proved, this makes the corresponding summand on the right of (\ref{eqn:R-diagonal1}) regular. 

Suppose instead that $i$ and $j$ lie in different color-charge blocks $\mathcal{H}_{\mathbf N}$ and $\mathcal{H}_{\mathbf M}$. Algebraic independence of the $b_\alpha's$ gives $\nu_{\mathbf N}-\nu_{\mathbf M}\in K^\times$, and
(\ref{eqn:eigenvalue_expand}) yields
\begin{equation}\label{eqn:eigenvalue_difference_expand}
 \rho_i-\rho_j
 =s^{-1}(\nu_{\mathbf N}-\nu_{\mathbf M})+O(1)
 \in s^{-1}R^\times.
\end{equation}
The leading scalar $s^{-1}B$ is independent of $Q$, so by (\ref{eqn:divisor_expand})
\[
 V(H\star_Q)\in\operatorname{End}_R(\mathscr H_R).
\]
Applying $V$ to the eigenvalue equation
\[
 (H\star_Q)P_i=\rho_iP_i
\]
and multiplying on the left by $P_j$ and on the right by $P_i$, we obtain, using $P_jP_i=0$ and $P_i^2=P_i$,
\begin{equation*}
P_jV(H\star_Q)P_i+P_j(H\star_Q)V(P_i)P_i=\rho_iP_jV(P_i)P_i.    
\end{equation*}
Multiplying the eigenvalue equation $(H\star_Q)P_j=\rho_jP_j$ by $V(P_i)P_i$ from the right, we obtain
\begin{equation*}
(H\star_Q)P_jV(P_i)P_i=\rho_jP_jV(P_j)P_j.    
\end{equation*}
Since the operators $H\star_Q$ and $P_j$ commute, taking the difference of the above two equations gives
\begin{equation}\label{eqn:derivative_projector_expand}
 P_jV(P_i)P_i =\frac{P_jV(H\star_Q)P_i}{\rho_i-\rho_j}=O(s).
\end{equation}
Apply both sides to the normalized idempotent $f_i$. Since $\Psi W_V=V\Psi$, taking the $i$-th column gives the identity
\[
 V(f_i)=\sum_k f_k(W_V)_{ki},
\]
which identifies the left-hand side of (\ref{eqn:derivative_projector_expand}) with $f_j(W_V)_{ji}$. Thus
\[
 (W_V)_{ij}=O(s).
\]
Combining this with (\ref{eqn:R-recursion_V}) and (\ref{eqn:eigenvalue_difference_expand}), we obtain
\begin{equation}\label{eqn:object_expand}
 (\mathsf{R}_1)_{ij}=O(s^2).
\end{equation}
Consequently, a cross-block summand in (\ref{eqn:R-diagonal1}) has order
\[
 (\rho_i-\rho_j)(\mathsf{R}_1)_{ij}^2=O(s^3),
\]
and is Rees regular.

\noindent{\bf Step 4:} {\it integrate coefficientwise in $Q$.}

We have proved that the right-hand side of (\ref{eqn:R-diagonal1}) belongs to $R$. The initial value $(\mathsf{R}_1)_{ii}|_{Q=0}$ is regular by Lemma~\ref{lem:gamma-initial} and the observation following its proof. Write
\[
 (\mathsf{R}_1)_{ii}=\sum_{d\ge0}r_{i,d}(s)Q^d,
 \qquad
 -\sum_{j\ne i}(\rho_i-\rho_j)(\mathsf{R}_1)_{ij}^2 =\sum_{d\ge0}F_{i,d}(s)Q^d.
\]
Equation (\ref{eqn:R-diagonal1}) says
\[
 d\,r_{i,d}(s)=F_{i,d}(s),\qquad d\ge1.
\]
Since the coefficient field has characteristic zero, division by the positive integer $d$ preserves $K[\![s]\!]$. Hence every $r_{i,d}(s)$ is Rees regular. This proves regularity of the diagonal entries and hence (\ref{eqn:regular_R1}).
\end{proof}

\begin{remark}[Why the auxiliary flat direction is essential]
The difference $\rho_i-\rho_j$ may vanish inside a color-charge block, so the first recursion cannot be solved by dividing by it. The auxiliary operator $B_\gamma$ was chosen so that $\beta_i-\beta_j$ is a unit for every $i\neq j$. Equation (\ref{eqn:R-recursion_ij}) contains the square of this difference: one factor comes from the first Givental recursion, and the second from differentiating the diagonalization of $B_\gamma$. The four-point correlator is therefore exactly the first-order estimate of an ingredient of the Frobenius structure missing from the pointwise projector construction.
\end{remark}

\subsection{Specialization of the first term}

The first term in Givental's formula depends only on norms of the idempotents. The original norm has a common power of $s$, but $V=Q\partial_Q$ annihilates that power. The remaining unit norm specializes by the product formulas of Sections \ref{sec:Frob} and \ref{sec:projector}.

\begin{lemma}[Specialization of the $\Delta$-term]
\label{lem:delta-specialization}
Fix a color-charge block $\mathcal{H}_{\mathbf{N}}$ with $\mathbf N=(N_1,\ldots,N_r)$ and write
\[
 \boldsymbol\lambda=(\lambda^1,\ldots,\lambda^r), \qquad \lambda^\alpha\vdash N_\alpha.
\]
Let
\[
 e_{\boldsymbol\lambda}=P_{\boldsymbol\lambda}(\id),
 \qquad
 g_{\boldsymbol\lambda} =\eta(e_{\boldsymbol\lambda},e_{\boldsymbol\lambda}),
 \qquad
 \Delta_{\boldsymbol\lambda}=g_{\boldsymbol\lambda}^{-1}.
\]
Then
\begin{equation}\label{eq:delta-individual-limit}
 \left.(V(\log\Delta_{\boldsymbol\lambda}))\right|_{s=0}
=\sum_{\alpha=1}^r V(\log\Delta^{(1)}_{\lambda^\alpha}),
\end{equation}
where $\Delta^{(1)}_{\lambda^\alpha}$ is the inverse of the norm of the idempotent (indexed by the partition $\lambda^\alpha$) in the Gromov-Witten theory of the rank-one case $\operatorname{Hilb}^{N_\alpha}(\mathbb C^2)$.

Let $I_N$ be the set of labels used in Gromov-Witten theory of $\text{Hilb}^N(\mathbb{C}^2)$, the rank-one, charge-$N$ case, with\footnote{The charge-$0$ case has norm $1$ and trivial $\mathsf R$-matrix.} $I_0=\{\emptyset\}$ and $d_0=1$. Put
\[
 d_N=|I_N|,\quad I_{\mathbf N}=\prod_{\alpha=1}^r I_{N_\alpha}.
\]
Then the contribution of the block to the first term satisfies
\begin{equation}\label{eq:delta-block-limit}
 \left.
 \frac1{48}V\left(\log \left(\prod_{\boldsymbol\lambda\in I_{\mathbf N}}\Delta_{\boldsymbol\lambda}\right)\right)\right|_{s=0}
=\sum_{\alpha=1}^r\left(\prod_{\beta\ne\alpha}d_{N_\beta}\right)
 \frac1{48}V\left(\log\left(\prod_{\lambda\in I_{N_\alpha}}\Delta^{(1)}_\lambda\right)\right).
\end{equation}
\end{lemma}

\begin{proof}
Put
\[
 \widehat g_{\boldsymbol\lambda}=\widehat\eta(e_{\boldsymbol\lambda},e_{\boldsymbol\lambda}).
\]
Since $\eta=t^{-c}\widehat\eta$,
\begin{equation}\label{eq:delta-rescaled-norm}
 g_{\boldsymbol\lambda}=t^{-c}\widehat g_{\boldsymbol\lambda},
 \qquad
 \Delta_{\boldsymbol\lambda}
 =t^c\widehat g_{\boldsymbol\lambda}^{-1}
 =s^{-c}\widehat g_{\boldsymbol\lambda}^{-1}.
\end{equation}
The element $\widehat g_{\boldsymbol\lambda}$ is a unit by Corollary~\ref{cor:product-extension}.  Since $V(s)=0$,
\begin{equation}\label{eq:delta-log-rescaled}
 V\left(\log\Delta_{\boldsymbol\lambda}\right)
 =-\frac{V(\widehat g_{\boldsymbol\lambda})}{\widehat g_{\boldsymbol\lambda}}.
\end{equation}
Thus the derivative is Rees regular even though $\Delta_{\boldsymbol\lambda}$ itself generally has a pole at $s=0$.

We use (\ref{eqn:project_limit}) in Theorem~\ref{thm:joint-rees} and the fact that the unit specializes blockwise to the product unit to deduce
\begin{equation}\label{eq:idempotent-product-limit}
 e_{\boldsymbol\lambda}(0,Q)=\bigotimes_{\alpha=1}^r e^{(1)}_{\lambda^\alpha}(Q),
 \qquad
 e^{(1)}_{\lambda^\alpha}=P^{(1)}_{\lambda^\alpha}(\id^{(1)}_{N_\alpha}).
\end{equation}
By Proposition~\ref{prop:stable-lattice}, the rescaled metric on the block has $s=0$ specialization
\[
 \left.\widehat\eta\right|_{s=0,\mathcal H_{\mathbf N}}
 =c_{\mathbf N}\bigotimes_{\alpha=1}^r\eta^{(1)}_{N_\alpha},
 \quad c_{\mathbf N}\in K^\times, \quad \text{where } \eta_{N_\alpha}^{(1)}=\eta_{\text{Hilb}^{N_\alpha}(\mathbb{C}^2)}.
\]
The factor $c_{\mathbf N}$ is independent of $Q$. Therefore, with
\[
 g^{(1)}_{\lambda^\alpha}
 =\eta^{(1)}_{N_\alpha}(e^{(1)}_{\lambda^\alpha},e^{(1)}_{\lambda^\alpha}),
 \quad
 \Delta^{(1)}_{\lambda^\alpha}=(g^{(1)}_{\lambda^\alpha})^{-1},
\]
we have
\begin{equation}\label{eq:norm-product-limit}
 \widehat g_{\boldsymbol\lambda}(0,Q)
 =c_{\mathbf N}\prod_{\alpha=1}^r g^{(1)}_{\lambda^\alpha}(Q),
 \qquad
 \widehat g_{\boldsymbol\lambda}(0,Q)^{-1}
 =c_{\mathbf N}^{-1}\prod_{\alpha=1}^r\Delta^{(1)}_{\lambda^\alpha}(Q).
\end{equation}
Apply $V(\log(-))$ and use $V(c_{\mathbf N})=0$. Combining \eqref{eq:delta-log-rescaled} and \eqref{eq:norm-product-limit}, we obtain \eqref{eq:delta-individual-limit}.

Finally, for fixed $\alpha$, every label in $I_{N_\alpha}$ occurs once for each choice of labels in the other factors, hence with multiplicity
\(
 \prod_{\beta\ne\alpha}d_{N_\beta}
\).
Summing \eqref{eq:delta-individual-limit} over $I_{\mathbf N}$ gives \eqref{eq:delta-block-limit}.
\end{proof}

\subsection{Tensor product description of the specialization}

Fix a color-charge block $\mathcal{H}_{\mathbf N}$. Let
\[
 \mathsf{R}^{\otimes}(z) =\bigotimes_{\alpha=1}^r \mathsf{R}^{(1)}_{N_\alpha}(z), \quad \mathsf{R}^{(1)}_{N_\alpha}=\mathsf{R}^{\text{Hilb}^{N_\alpha}(\mathbb{C}^2)}
\]
be the tensor product of the $\mathsf{R}$-matrices of the Gromov-Witten theories of $\text{Hilb}^{N_\alpha}(\mathbb{C}^2)$. Its $z^1$-coefficient is
\[
 \mathsf{R}_1^{\otimes}
 =\sum_{\alpha=1}^r 1\otimes\cdots\otimes \mathsf{R}^{(1)}_{1,N_\alpha}\otimes\cdots\otimes1.
\]
The projectors, product, and rescaled metric specialize to this tensor product theory by Theorem~\ref{thm:joint-rees} and Corollary~\ref{cor:product-extension}. We now compare $\mathsf{R}_1|_{\mathcal{H}_{\mathbf N}}$ with $\mathsf{R}^{\otimes}_1$.

\begin{proposition}[Large framing limit of $\mathsf{R}_1$]
\label{prop:boundary-calibration}
At $s=0$, the diagonal entries of $\mathsf{R}_1$ and the first term in (\ref{eqn:Givental_formula_general}) agree, block by block, with those of the tensor product of the Gromov-Witten theories of the rank-one cases $\text{Hilb}^{N_\alpha}(\mathbb{C}^2)$.
\end{proposition}

\begin{proof}
Fix a block $\mathcal{H}_{\mathbf N}$.

\noindent{\bf Step 1:} {\it specialize the normalized frame and the divisor direction.}

By \eqref{eq:idempotent-product-limit} and \eqref{eq:norm-product-limit}, the $\widehat\eta$-normalized idempotent frame specializes to the tensor-product normalized frame, multiplied on the whole block by the common scalar $c_{\mathbf N}^{-1/2}$. This scalar is independent of $Q$. Hence
\begin{equation}\label{eqn:W_V_limit}
 \left.W_V\right|_{s=0}=W_V^{\otimes},
\end{equation}
where $W_V^{\otimes}$ is defined analogously to (\ref{eqn:W-matrix}) using the tensor product Gromov-Witten theory.

The leading-order block scalar $s^{-1}\nu_{\mathbf N}$ in $H\star_Q$ is also independent of $Q$ and is therefore annihilated by $V$. The divisor operator specializes to the sum of the scalar-subtracted divisor operators of the rank-one cases. Thus the quantities
\[
 g_{ij}(Q) :=\left.(\rho_i-\rho_j)\right|_{s=0}
\]
are precisely the differences of the eigenvalues $du_i(\partial_H)-du_j(\partial_H)$ of the divisor operator of the tensor-product theory. The scalar shift is common to the block and cancels in these differences; these are not differences $u_i-u_j$ of canonical coordinates.

\noindent{\bf Step 2:} {\it compare the off-diagonal data whenever the difference is nonzero.}

Specializing the first recursion to $s=0$ gives
\begin{equation}\label{eqn:R-recursion2}
 W^{\otimes}_{V,ij}=g_{ij}\cdot(\mathsf{R}_1|_{s=0})_{ij}.
\end{equation}
The tensor product $\mathsf{R}$-matrix satisfies the same recursion:
\[
W^{\otimes}_{V,ij}=g_{ij}\cdot(\mathsf{R}_1^{\otimes})_{ij}.
\]
If $g_{ij}\ne0$, the coefficient ring after specialization is a domain, so one may cancel the nonzero element $g_{ij}$ even when it is not a unit. Hence
\begin{equation}\label{eqn:R_off_diag_limit}
 (\mathsf{R}_1|_{s=0})_{ij}=(\mathsf{R}_1^{\otimes})_{ij}
 \qquad\text{whenever }g_{ij}\ne0.
\end{equation}
If $g_{ij}=0$ identically, no division is made; that pair contributes $0$ to the diagonal recursion in both theories.

\noindent{\bf Step 3:} {\it compare the diagonal recursions and their initial values.}

By (\ref{eqn:R-diagonal1}), the diagonal recursion depends on an off-diagonal pair only through
\[
 g_{ij}\cdot(\mathsf{R}_1)_{ij}^{\,2}.
\]
For $g_{ij}\ne0$, equality follows from (\ref{eqn:R_off_diag_limit}); for $g_{ij}=0$, both
contributions vanish. Hence the $s=0$ specialization of $\mathsf{R}_1$ and $\mathsf{R}_1^{\otimes}$ have the same diagonal recursions.

Their initial values at $Q=0$ also agree. Lemma~\ref{lem:gamma-initial} expresses the initial value of $\mathsf{R}_1$ as a sum of reciprocals of tangent weights. At the large framing limit $s=0$, every cross-color reciprocal tangent weight is $O(s)$ and vanishes, while the same-color weights are exactly the tangent weights of the rank-one Hilbert-scheme factors. The surviving sum is therefore the sum of the initial values of the $\text{Hilb}^{N_\alpha}(\mathbb{C}^2)$ theories, which is the diagonal initial value of $\mathsf{R}_1^{\otimes}$. Coefficientwise integration in $Q$ now gives
\begin{equation}\label{eqn:R_diag_limit}
 (\mathsf{R}_1|_{s=0})_{ii}=(\mathsf{R}_1^{\otimes})_{ii}.
\end{equation}

The specialization of the first term is already obtained in Lemma~\ref{lem:delta-specialization}. The proof is now complete.
\end{proof}

\section{Rank-one theories, tensor products, and the large framing limit}\label{sec:rank-one-limit}

We have identified the large framing limits of the ingredients in Givental's formula. Three pieces of
bookkeeping remain before we can apply them:

\begin{enumerate}
\item the geometric rank-one divisor differs from the oscillator divisor by a scalar operator;
\item the genus-one formula for a tensor product contains cross terms; 
\item the full divisor eigenvalues contain a blockwise scalar of order $s^{-1}$, whose positive-degree contribution must be shown to vanish.
\end{enumerate}

We now treat these points. Throughout this section, the variable $z$ keeps track of the charges (i.e. number of points) of Hilbert schemes appearing in the tensor products, while $Q$ is the (signed) Novikov variable. All identities are formal in $z,Q$ and coefficientwise in $Q$.

\subsection{Rank-one input and the scalar shift}

Put
\begin{equation}\label{eqn:partition_gen_func}
 P(z)=\prod_{k\geq1}(1-z^k)^{-1}, \qquad L(z)=\log P(z).
\end{equation}
$P(z)$ is the generating function for the numbers of partitions. For the Hilbert scheme theory, define
\[
 F^{\mathrm{geom}}(z,Q)=\sum_{N\geq0}z^N\F^+_{1,N}(Q),
 \qquad
 T^{\mathrm{geom}}(z,Q)=\sum_{N\geq0}z^N\tau_{1,N}(Q).
\]
The positive-degree specialization of the result of \cite{ILPT} for Hilbert scheme cases, after the variable change $q=-Q$, is
\begin{equation}\label{eqn:rank_1_formula}
 F^{\mathrm{geom}}(z,Q) =-\frac h{24p}\bigl(1+L(z)\bigr)T^{\mathrm{geom}}(z,Q).
\end{equation}
The subtraction of the degree-$0$ term and the sign change are both included here.

The tensor-product theory of Section \ref{sec:limit_genus1} is expressed instead in the \emph{oscillator normalization}. On the rank-one, charge-$N$ sector, let $A_N^{\mathrm{osc}}$ be the scalar-subtracted divisor operator of (\ref{eqn:divisor_expand2})--(\ref{eqn:divisor_expand3}), and write
\[
 T_N(Q)=\Tr\bigl((A_N^{\mathrm{osc}}-A_N^{\mathrm{osc}}|_{Q=0})\bigr),
\]
namely $T_N$ is the trace of the positive $Q$-degree part the oscillator divisor operator. Let $F_N(Q)$ be the one form corresponding to genus one positive-degree Gromov-Witten theory of $\text{Hilb}^N(\mathbb{C}^2)$ evaluated on the oscillator divisor vector field $H-hN\frac Q{1-Q}\id$, and set
\[
 T(z,Q)=\sum_{N\geq0}T_N(Q)z^N,
 \qquad
 F(z,Q)=\sum_{N\geq0}F_N(Q)z^N.
\]
Finally, let $\mathsf{R}_{1,N}$ denote the initial value of $\mathsf{R}_1$ of the rank-one, charge-$N$ case and put
\[
 S_0(z)=\sum_{N\geq0}\Tr(\mathsf{R}_{1,N})z^N.
\]

\begin{lemma}[Constant-map trace in rank-one]
\label{lem:rankone-R1-trace}
The series $S_0(z)$ is independent of $Q$ and satisfies
\[
 S_0(z)=-\frac h{12p}P(z)L(z).
\]
\end{lemma}

\begin{proof}
Let $d$ be the exterior differential on the Frobenius manifold for the Gromov-Witten theory in the rank-one case. The second Givental recursion is
\[
 d\mathsf{R}_1+W\mathsf{R}_1=[dU,\mathsf{R}_2].
\]
Taking the trace gives
\[
 d\,\Tr(\mathsf{R}_1)=-\Tr(W\mathsf{R}_1).
\]
In the normalized idempotent frame, $W^T=-W$ and $\mathsf{R}_1^T=\mathsf{R}_1$. Hence
\[
 \Tr(W\mathsf{R}_1)=\Tr((W\mathsf{R}_1)^T)
 =\Tr(\mathsf{R}_1^TW^T)=-\Tr(\mathsf{R}_1W)=-\Tr(W\mathsf{R}_1),
\]
so $\Tr(W\mathsf{R}_1)=0$. Thus $\Tr(\mathsf R_1)$ is constant in the flat variables. By divisor equation, $V(\Tr(\mathsf R_1))=0$, so its formal $Q$-series equals its initial value $\Tr(\mathsf R_{1,N})$.

At $Q=0$, Lemma~\ref{lem:gamma-initial} gives
\[
 \Tr(\mathsf{R}_{1,N})
 =-\frac1{12}
 \sum_{\lambda\vdash N}
 \sum_{\chi\in\operatorname{Wt}(T_\lambda\operatorname{Hilb}^{N}(\mathbb C^2))}
 \frac1\chi.
\]
By the calculations of \cite[Section 1]{ILPT}, 
\[
 \sum_{N\geq1}z^N
 \sum_{\lambda\vdash N}
 \sum_{\chi\in\operatorname{Wt}(T_\lambda\operatorname{Hilb}^{N}(\mathbb C^2))}
 \frac1\chi
 =\frac hpP(z)L(z).
\]
The charge-zero term is zero. Multiplication by $-1/12$ proves the desired formula.
\end{proof}

\begin{lemma}[Rank-one scalar shift]\label{lem:rankone-shift}
\begin{equation}\label{eqn:rankone-shift}
 F(z,Q)=-\frac h{24p}\bigl(1+L(z)\bigr)T(z,Q).
\end{equation}
\end{lemma}

\begin{proof}
By the complete rank-one normalization (\ref{eqn:scalar_divisor}), the operators given by geometric divisor and the oscillator divisor differ on charge $N$ by
\[
 c_N(Q)\id,
 \qquad
 c_N(Q)=hN\frac Q{1-Q}.
\]
Consequently,
\[
 T_N^{\mathrm{geom}}-T_N=c_N(Q)\,\dim H_\T^\bullet(\operatorname{Hilb}^{N}(\mathbb C^2))_{\loc}.
\]
Since $\dim H_\T^\bullet(\operatorname{Hilb}^{N}(\mathbb C^2))_{\loc}=p(N)$, the number of partitions of $N$, summing over $N$ yields
\begin{equation}\label{eqn:trace_difference}
 T^{\mathrm{geom}}-T
 =h\frac Q{1-Q}\sum_{N\geq0}Np(N)z^N
 =h\frac Q{1-Q}z\frac{d}{dz}P(z).
\end{equation}

Let $e$ be the unit of the quantum product $\star$. In canonical coordinates,
\[
 e(u_i)=1,
 \qquad
 e(\Delta_i)=0.
\]
Givental's genus-one formula therefore gives
\[
 \omega_1(e)=\frac12\sum_i(\mathsf{R}_1)_{ii}=\frac12\Tr (\mathsf{R}_1).
\]
Replacing the geometric divisor vector field by the oscillator vector field subtracts $c_Ne$. Hence, on charge $N$,
\begin{equation}\label{eqn:g1_funct_difference}
 F_N-F_N^{\mathrm{geom}}=-\frac{c_N(Q)}2\Tr(\mathsf{R}_{1,N}).
\end{equation}
After summing over $N$ and using Lemma~\ref{lem:rankone-R1-trace},
\[
 \begin{aligned}
 F-F^{\mathrm{geom}}
 &=-\frac h2\frac Q{1-Q}\,z\frac{d}{dz}S_0(z)\\
 &=\frac{h^2}{24p}\frac Q{1-Q}
   z\frac{d}{dz}\bigl(P(z)L(z)\bigr)\\
 &=\frac{h^2}{24p}\frac Q{1-Q}
   \bigl(1+L(z)\bigr)z\frac{d}{dz}P(z)\\
 &=\frac h{24p}\bigl(1+L(z)\bigr)
   \bigl(T^{\mathrm{geom}}-T\bigr),
 \end{aligned}
\]
where we used $\frac{d}{dz}P=P\frac{d}{dz}L$. Combining this identity with (\ref{eqn:rank_1_formula}) proves (\ref{eqn:rankone-shift}).
\end{proof}

\subsection{Tensor product theory in genus one}

We next compute the genus-one expression for the $r$-fold tensor product of the oscillator-normalized rank-one theories. The argument is most transparent first at fixed charges.

For a rank-one, charge-$N$ factor, write
\[
 d_N=p(N),
 \qquad
 S_N=\Tr(\mathsf{R}_{1,N}),
\]
and keep the notation $T_N,F_N$ above. For a color-charge vector $\mathbf N=(N_1,\ldots,N_r)$, let $T_{\mathbf N}^{\otimes}$ and $F_{\mathbf N}^{\otimes}$ denote respectively the positive-degree divisor trace and the positive-degree genus-one expression of the tensor-product theory, evaluated on the sum of the factorwise oscillator divisor vector fields.

\begin{lemma}[Tensor product]\label{lem:tensor-product-fixed}
For every color-charge vector $\mathbf N$, we have
\begin{equation}\label{eqn:tensor_product_trace}
 T_{\mathbf N}^{\otimes} =\sum_{a=1}^r T_{N_a}\prod_{b\ne a}d_{N_b},
\end{equation}
and
\begin{equation}\label{eqn:tensor_product_g1_formula}
 F_{\mathbf N}^{\otimes}
=\sum_{a=1}^rF_{N_a}\prod_{b\ne a}d_{N_b}+\frac12\sum_{a\ne b}T_{N_a}S_{N_b}\prod_{c\ne a,b}d_{N_c}.
\end{equation}
\end{lemma}

\begin{proof}
For two semisimple calibrated\footnote{Following \cite{G}, a calibration of a Frobenius structure is a choice of a fundamental solution of the flat connection.} Frobenius structures $\mathfrak A$ and $\mathfrak B$, the product calibration is the tensor product of fundamental solutions. Hence
\begin{equation}\label{eqn:tensor_product_R1}
\mathsf{R}_{\mathfrak A\otimes \mathfrak B}(z)=\mathsf{R}_{\mathfrak A}(z)\otimes \mathsf{R}_{\mathfrak B}(z),
\quad\text{in particular }
\mathsf{R}_{1,\mathfrak A\otimes \mathfrak B}=\mathsf{R}_{1,\mathfrak A}\otimes1+1\otimes \mathsf{R}_{1,\mathfrak B}.
\end{equation}
Primitive idempotents are given by tensor products, their norms multiply, and their canonical coordinates add. Therefore, the first term in Givental's formula is the sum of the first terms of the two factors, each repeated by the number of primitive idempotents in the other factor.

For the $\mathsf{R}_1$-term, let $A_{\mathfrak A}^+$ and $A_{\mathfrak B}^+$ denote the positive-degree parts of the two divisor multiplication operators. Expanding (\ref{eqn:tensor_product_R1}) and the sum of divisor vector fields gives, in addition to the two factorwise contributions, the cross term
\begin{equation}\label{eqn:R_term_mixed}
 \frac12\Tr(\mathsf{R}_{1,\mathfrak A})\Tr(A_{\mathfrak B}^+)
 +\frac12\Tr(A_{\mathfrak A}^+)\Tr(\mathsf{R}_{1,\mathfrak B}).
\end{equation}
Iterating this two-factor calculation proves (\ref{eqn:tensor_product_trace})--(\ref{eqn:tensor_product_g1_formula}). The ordered sum in (\ref{eqn:tensor_product_g1_formula}) records the two orientations in (\ref{eqn:R_term_mixed}); the factor $1/2$ is exactly the coefficient in Givental's formula.
\end{proof}

Define the series for the tensor product theory,
\[
 T^{(r)}(z,Q)=\sum_{N\geq0}z^N\sum_{|\mathbf N|=N}T_{\mathbf N}^{\otimes}(Q),\quad
 F^{(r)}(z,Q)=\sum_{N\geq0}z^N\sum_{|\mathbf N|=N}F_{\mathbf N}^{\otimes}(Q).
\]
Summing Lemma~\ref{lem:tensor-product-fixed} over all color-charge vectors $\mathbf{N}$ gives
\begin{equation}\label{eqn:sum_all_colors}
 T^{(r)}=rP^{r-1}T,
 \qquad
 F^{(r)}=rP^{r-1}F+\binom r2P^{r-2}S_0T.
\end{equation}
Indeed, the first identity chooses the unique active factor. In the second, the first term again chooses one active genus-one factor, while the second chooses one factor contributing $T$, one different factor contributing $S_0$, and lets all remaining factors contribute their state-space series $P$.

\begin{lemma}[Tensor-product genus-one series]\label{lem:tensor-product}
\begin{equation}\label{eqn:tensor-products}
 F^{(r)}(z,Q) =-\frac h{24p}\bigl(1+rL(z)\bigr)T^{(r)}(z,Q).
\end{equation}
\end{lemma}

\begin{proof}
Substitute Lemma~\ref{lem:rankone-shift} and Lemma~\ref{lem:rankone-R1-trace} into (\ref{eqn:sum_all_colors}):
\[
 \begin{aligned}
 F^{(r)}
 &=-\frac h{24p}rP^{r-1}(1+L)T -\frac h{12p}\binom r2P^{r-1}LT\\
 &=-\frac h{24p}rP^{r-1}\bigl(1+rL\bigr)T\\
 &=-\frac h{24p}\bigl(1+rL\bigr)T^{(r)}.
 \end{aligned}
\]
\end{proof}

\subsection{Cancellation of the divergent color-charge scalar}

Fix a color-charge block $C$ with color-charge vector $\mathbf N$. Let $I_C$ be the set of primitive labels in that block and put
\[
 \mathcal S_C(Q)=\sum_{i\in I_C}(\mathsf{R}_1)_{ii}(Q).
\]
The eigenvalue of the divisor operator on every eigenline in $C$ has the form
\[
 \rho_i=s^{-1}\nu_C+\widetilde\rho_i,
 \qquad
 \nu_C=\sum_{\alpha=1}^rb_\alpha N_\alpha,
\]
by (\ref{eqn:eigenvalue_expand}). The term $s^{-1}\nu_C$ is not part of the tensor product, and its contribution to the $\mathsf{R}_1$-part of Givental's formula is
\[
 \frac12s^{-1}\nu_C\,\mathcal S_C(Q).
\]
Only its positive $Q$-degree part is relevant for us.

\begin{lemma}[Block-trace cancellation]\label{lem:block-trace}
For every $d>0$,
\begin{equation}\label{eqn:R1_positive_degree_block}
 [Q^d]\bigl(\mathcal S_C(Q)-\mathcal S_C(0)\bigr)\in s^3K[\![s]\!].
\end{equation}
Consequently
\[
\frac12s^{-1}\nu_C [Q^{>0}]\mathcal S_C(Q)=O(s^2).
\]
\end{lemma}

\begin{proof}
Sum the diagonal recursion (\ref{eqn:R-diagonal1}) over $i\in I_C$:
\begin{equation}\label{eqn:R1_positive_degree_block1}
 V(\mathcal S_C)
=-\sum_{i\in I_C}\sum_{j\ne i}(\rho_i-\rho_j)(\mathsf{R}_1)_{ij}^2,
\quad V=Q\partial_Q.
\end{equation}
If $i,j\in I_C$, the term indexed by $(j,i)$ is the negative of the term indexed by $(i,j)$: the $\rho_i-\rho_j$ changes sign, while $\mathsf{R}_1$ is symmetric. Hence all terms with $i,j\in I_C$ cancel in pairs, leaving
\begin{equation}\label{eqn:R1_positive_degree_block2}
 V(\mathcal S_C) =-\sum_{\substack{i\in I_C\\ j\notin I_C}}(\rho_i-\rho_j)(\mathsf{R}_1)_{ij}^2.
\end{equation}
For a cross-block pair, (\ref{eqn:object_expand}) gives
\[
 \rho_i-\rho_j\in s^{-1}R^\times,
 \qquad
 (\mathsf{R}_1)_{ij}=O(s^2).
\]
Every summand in (\ref{eqn:R1_positive_degree_block2}) is therefore $O(s^3)$, and so $V(\mathcal S_C)=O(s^3)$ coefficientwise in $Q$. If
\[
 \mathcal S_C(Q)-\mathcal S_C(0)=\sum_{d\geq1}s_{C,d}(s)Q^d,
\]
then (\ref{eqn:R1_positive_degree_block2}) gives $d\,s_{C,d}(s)=O(s^3)$. Since the coefficient field has characteristic $0$, division by $d$ proves (\ref{eqn:R1_positive_degree_block}). Multiplication by $s^{-1}\nu_C/2$ gives the claimed Rees estimate.
\end{proof}

\subsection{The fixed-charge large framing limit}

The following elementary identity about $L(z)=\log P(z)$ is needed in the coefficient extraction:
\begin{equation}\label{eqn:log_of_partitions}
 L(z)=\log P(z)=\sum_{m\geq1}\sigma_{-1}(m)z^m.
\end{equation}

\begin{theorem}[Large framing limit]\label{thm:coulomb-limit}
For every $r\geq2$ and $n\geq1$,
\begin{equation}\label{eqn:coulomb-limit}
 \lim_{t\to\infty}\F^+_{r,n}(tb;Q)
=-\frac h{24p}\left(\tau_{r,n}(Q)+r\sum_{m=1}^{n-1} \sigma_{-1}(m)\tau_{r,n-m}(Q)\right).
\end{equation}
The limit is coefficientwise in $Q$.
\end{theorem}

\begin{proof}
Proposition~\ref{prop:boundary-calibration} identifies each color-charge block of the large framing limit with the corresponding tensor product of oscillator-normalized rank-one theories. Lemma~\ref{lem:block-trace} shows that the omitted scalar $s^{-1}\nu_C$ contributes $0$ in positive degree at $s=0$.

It remains to identify the trace of the divisor operator in the tensor product theory. In the stable-envelope chart $\Phi_t$, the diagonal part of the divisor is the sum of the rank-one oscillator divisors. Every cross-color term changes the color-charge vector, so its trace on the direct sum of blocks is zero. Therefore
\begin{equation}\label{eqn:trace_gen_func}
 T^{(r)}(z,Q)=\sum_{N\geq1}\tau_{r,N}(Q)z^N.
\end{equation}
This is also immediate from the complete higher-rank divisor formula (\ref{eqn:divisor_quantum_part}), which is independent of the individual framing parameters.

Apply Lemma~\ref{lem:tensor-product} and extract the coefficient of $z^n$. The term without $L$ contributes $\tau_{r,n}$. Using (\ref{eqn:log_of_partitions}), the term $rLT^{(r)}$ contributes
\[
 r\sum_{m=1}^{n}\sigma_{-1}(m)\tau_{r,n-m}.
\]
The endpoint $m=n$ vanishes because $\tau_{r,0}=0$. This gives (\ref{eqn:coulomb-limit}).
\end{proof}

\section{Proof of Theorem \ref{thm:main}}\label{sec:pf_main_thm}

Fix a $Q^d$. Every expression in this Section is then a rational function of the framing parameters over
\[
 k_0=\mathbb Q(\epsilon_1,\epsilon_2).
\]
Theorem~\ref{prop:polar_part} determines all finite polar hyperplanes and their residues. Theorem~\ref{thm:coulomb-limit} determines the remaining polynomial part.  What remains is therefore an algebraic uniqueness argument.

For $r=1$, Theorem~\ref{thm:main} is exactly (\ref{eqn:rank_1_formula}), already proven in \cite{ILPT}.

Assume now that $r\geq2$. The complete higher-rank divisor formula (\ref{eqn:divisor_quantum_part}) shows that every $\tau_{r,N}(Q)$ is independent of the framing parameters $a_\alpha's$. Recall the coefficients $b_{r,m}$ from (\ref{eqn:coefficients}). Define
\begin{equation}\label{eqn:trace_expression}
 G_{r,n}(Q) =-\frac h{24p}\sum_{m=0}^{n-1} b_{r,m}\tau_{r,n-m}(Q).
\end{equation}

\subsection{Matching the finite polar parts}

\begin{lemma}[Matching polar parts]\label{lem:matching-parts}
Coefficientwise in $Q$, the functions $\F^+_{r,n}$ and $G_{r,n}$ have the same polar part at every finite hyperplane in the space $\text{Spec}\,\CC[\epsilon_1,\epsilon_2,a_1,...,a_r]$ of framing parameters.
\end{lemma}

\begin{proof}
Let
\[
 x=i\epsilon_1+j\epsilon_2, \quad i,j\geq1,
\]
and consider the wall (recall that $a_{\alpha\beta}=a_\alpha-a_\beta$)
\[
 a_{\alpha\beta}=x.
\]
Theorem~\ref{prop:polar_part} says that a pole occurs at this wall only when $ij\leq n$, that the pole is simple, and that its residue is
\begin{equation}\label{eqn:residue5}
 \frac1{24}\tau_{r,n-ij}(Q).
\end{equation}
If $ij=n$, this residue is zero because $\tau_{r,0}=0$, and the sum (\ref{eqn:trace_expression}) contains no term with $m=n$.

Assume $ij<n$.  In the definition (\ref{eqn:coefficients}) of $b_{r,m}$, the divisor $\delta=i$ gives
\[
 x_{\delta,ij}=i\epsilon_1+j\epsilon_2=x.
\]
Its contribution to $G_{r,n}$ is
\[
 -\frac h{24p}\frac{2p}{h}
 \frac{x}{x^2-a_{\alpha\beta}^2}
 \tau_{r,n-ij}
 =\frac1{24}\left(
 \frac1{a_{\alpha\beta}-x}
 +\frac1{-a_{\alpha\beta}-x}\right)
 \tau_{r,n-ij}.
\]
The first fraction has the residue (\ref{eqn:residue5}) at $a_{\alpha\beta}=x$.  The second fraction gives the  residue at the wall $a_{\beta\alpha}=x$ because $a_{\beta\alpha}=-a_{\alpha\beta}$.

Because $\epsilon_1$ and $\epsilon_2$ are algebraically independent, distinct ordered pairs $(i,j)$ define distinct linear forms $i\epsilon_1+j\epsilon_2$. Hence no other pair $(m,\delta)$ contributes a
pole at the generic point of the chosen wall. Theorem \ref{prop:polar_part} excludes every other finite pole in the space $\text{Spec}\,\CC[\epsilon_1,\epsilon_2,a_1,...,a_r]$ of framing parameters. Thus all finite polar parts agree.
\end{proof}

\begin{lemma}[No finite poles implies polynomiality]
\label{lem:no-poles-polynomial}
If $f\in k_0(a_1,\ldots,a_r)$ has no pole at any codimension-$1$ point of the space $\text{Spec}\,\CC[\epsilon_1,\epsilon_2,a_1,...,a_r]$ of framing parameters, then $f\in k_0[a_1,\ldots,a_r]$.
\end{lemma}

\begin{proof}
Write $f=P/Q$ with coprime polynomials $P,Q\in k_0[a_1,\ldots,a_r]$. If $Q$ were nonconstant, any irreducible factor of $Q$ would define a codimension-$1$ prime along which $f$ has a pole, because that factor does not divide $P$. Thus $Q\in k_0^\times$, and $f$ is a polynomial.
\end{proof}

By Lemma~\ref{lem:matching-parts}, Theorem \ref{prop:polar_part}, and the simplicity of all
finite poles, the difference
\begin{equation}\label{eqn:F-minus-G}
 D_{r,n}(a;Q)=\F^+_{r,n}(Q)-G_{r,n}(Q)
\end{equation}
has no codimension-$1$ poles. Lemma~\ref{lem:no-poles-polynomial} therefore shows, coefficientwise in $Q$, that $D_{r,n}$ is a polynomial in the framing parameters $a_\alpha's$.

\subsection{Detection of the polynomial part at infinity}

The limit $t\to \infty$ is understood algebraically as $s=0$ where $s=t^{-1}$. Thus $f(tb)=O(t^{-1})$ means that, after substitution $a_\alpha=tb_\alpha$, the resulting Laurent series belongs to
$t^{-1}k_0(b_1,\ldots,b_r)[\![t^{-1}]\!]$.

\begin{lemma}[Detection on a generic ray in the space of framing parameters]
\label{lem:generic-ray-polynomial}
Let $b_1,\ldots,b_r$ be algebraically independent over the field $k_0$. If $P\in k_0[a_1,\ldots,a_r]$ and
\[
 P(tb_1,\ldots,tb_r)=O(t^{-1}),  \qquad t\to\infty,
\]
then $P=0$.
\end{lemma}
For $m>0$, the explicit formula for $b_{r,m}$ gives, along the generic ray $a_\alpha=tb_\alpha$,
\begin{equation}\label{eqn:coeff_b_expand}
 \begin{aligned}
 b_{r,m}(tb)
 =r\sigma_{-1}(m)
 +\frac{2p}{h}\sum_{\alpha<\beta}\sum_{\delta\mid m}
 \frac{x_{\delta,m}}{x_{\delta,m}^2-t^2(b_\alpha-b_\beta)^2}
 =r\sigma_{-1}(m)+O(t^{-2}).
 \end{aligned}
\end{equation}
Since $b_{r,0}=1$, equations (\ref{eqn:trace_expression}) and (\ref{eqn:coeff_b_expand}) give
\begin{equation}\label{eqn:G_expand}
 \begin{aligned}
 G_{r,n}(tb;Q)
=-\frac h{24p}\left(\tau_{r,n}(Q)+r\sum_{m=1}^{n-1}\sigma_{-1}(m)\tau_{r,n-m}(Q)\right)+O(t^{-2}).
 \end{aligned}
\end{equation}
The constant term is exactly the large framing limit in (\ref{eqn:coulomb-limit}). Therefore the difference (\ref{eqn:F-minus-G}) satisfies 
\begin{equation}\label{eqn:F-minus-G_expand}
 D_{r,n}(tb;Q)=O(t^{-1})
\end{equation}
coefficientwise in $Q$.

\begin{proof}[Proof of Theorem~\ref{thm:main} for $r\geq2$]
Fix $d\geq1$. By (\ref{eqn:F-minus-G}), the coefficient
\[
 [Q^d]D_{r,n}
\]
is a polynomial in the framing parameters $a_\alpha's$ over $k=\mathbb Q(\epsilon_1,\epsilon_2)$. Equation (\ref{eqn:F-minus-G_expand}) says that its restriction to the generic ray in the space $\text{Spec}\,\CC[\epsilon_1,\epsilon_2,a_1,...,a_r]$ of framing parameters is $O(t^{-1})$. Lemma \ref{lem:generic-ray-polynomial} implies
\[
 [Q^d]D_{r,n}=0.
\]
This holds for every $d>0$, hence
\[
 \F^+_{r,n}(Q)=G_{r,n}(Q).
\]
Substituting (\ref{eqn:trace_expression}) gives the desired formula
\end{proof}

\subsection{Evaluation of the divisor trace}

Theorem \ref{thm:main} reduces the remaining calculation to the trace
\[
 \tau_{r,N}(Q)=\Tr\bigl(H\star_Q-H\cup\bigr).
\]
We now evaluate this trace directly using the formula for $H\star_Q$ written in the Fock space. The calculation is formal and coefficientwise in both the charge variable $z$ and the Novikov variable $Q$.

\subsubsection{The trace of one Heisenberg number operator}

Let
\[
 |\boldsymbol\lambda\rangle =|\lambda^1,\ldots,\lambda^r\rangle
\]
be the standard colored-partition basis of the $r$-boson Fock space. For a partition $\lambda$, write $m_k(\lambda)$ for the multiplicity of the part $k$. In the Heisenberg convention of the Jordan-quiver divisor formula, the diagonal part of $\beta_{-k}\beta_k$ is
\begin{equation}\label{eqn:diag_1}
\left.\beta_{-k}\beta_k\right|_{|\boldsymbol\lambda\rangle}=-k\sum_{\alpha=1}^r m_k(\lambda^\alpha).
\end{equation}
The cross-color terms change the color-charge vector and contribute $0$ to the trace. Formula (\ref{eqn:diag_1}) is therefore enough for the trace.

For one color,
\[
\sum_{\lambda}m_k(\lambda)z^{|\lambda|}
=\left(\sum_{m\geq0}mz^{km}\right)\prod_{\ell\ne k}(1-z^\ell)^{-1}.
\]
Using the identities,
\[
 \sum_{m\geq0}mz^{km}=\frac{z^k}{(1-z^k)^2},
 \qquad
 \prod_{\ell\ne k}(1-z^\ell)^{-1}=P(z)(1-z^k),
\]
we obtain
\begin{equation}\label{eqn:gen_func_multiplicity}
 \sum_{\lambda}m_k(\lambda)z^{|\lambda|}=P(z)\frac{z^k}{1-z^k}.
\end{equation}
Choosing the color that contributes the number operator and letting all other colors contribute their partition functions gives
\begin{equation}\label{eqn:trace_diag}
 \Tr\left(z^{\mathrm{energy}}\beta_{-k}\beta_k\right)=-rk\,P(z)^r\frac{z^k}{1-z^k}.
\end{equation}
This proves the sign, the factor $k$, and the multiplicity $r$ in the chosen Heisenberg convention.

\subsubsection{Higher rank}

For $r\geq2$, the complete positive-degree divisor formula has no scalar correction and is
\begin{equation}\label{eqn:divisor_formula_repeat}
 \bigl(H\star_Q-H\cup\bigr)^+
 =h\sum_{k\geq1}k\frac{Q^k}{1-Q^k}\,
 \beta_{-k}\beta_k.
\end{equation}
At fixed charge, the state space is finite dimensional, so trace and coefficient extraction may be interchanged. Combining (\ref{eqn:divisor_formula_repeat}) and (\ref{eqn:trace_diag}), we obtain
\begin{equation}\label{eqn:trace_gen_func1}
 \sum_{N\geq1}\tau_{r,N}(Q)z^N =-rhP(z)^r\sum_{k\geq1}k^2\frac{Q^k}{1-Q^k}\frac{z^k}{1-z^k},
 \qquad r\geq2.
\end{equation}

Since
\[
 \frac{z^k}{1-z^k}=\sum_{s\geq1}z^{sk},
 \qquad
 P(z)^r=\sum_{M\geq0}p_r(M)z^M,
\]
the coefficient of $z^N$ in (\ref{eqn:trace_gen_func1}) is
\begin{equation}\label{eqn:trace_answer}
 \tau_{r,N}(Q)
 =-rh\sum_{k=1}^{N}k^2\frac{Q^k}{1-Q^k}\sum_{s=1}^{\lfloor N/k\rfloor}p_r(N-sk).
\end{equation}
Here $s$ is the positive multiplicity of the distinguished part $k$; after removing its energy $sk$, the remaining colored partition has charge $N-sk$.

{\large\part{Degree zero invariants}}
In this part, we compute the degree-zero invariants $\langle H\rangle_{1,1,0}^{\M(r,n)}$. This part can be read independently. The rank $r=1$ case was already proven in \cite[Section 1]{ILPT}.

\section{The formulae}
Recall from (\ref{eqn:partition_gen_func}) the two generating functions $P(z)$ and $L(z)$. Define further
\begin{equation}
E_2(z)=\sum_{k\geq1}k\frac{z^k}{1-z^k}= z\frac{dL}{dz},\quad
 E_3(z)=\sum_{k\geq1}k^2\frac{z^k}{1-z^k}.
\end{equation}
For $\alpha<\beta$, define
\begin{equation}\label{eq:cross-series-intro}
 \mathcal C_{\alpha\beta}(z)
 =\sum_{i,j\geq1}z^{ij}
 \left(\frac1{i\epsilon _1+j\epsilon _2- (a_{\alpha}-a_{\beta})}
 +\frac1{i\epsilon _1+j\epsilon _2 + (a_{\alpha}-a_{\beta})}\right).
\end{equation}
For $r\geq1$, consider the generating function
\begin{equation}
\F^{(0)}_{r,H}(z) := \sum_{n\geq0}z^n\left\langle H\right\rangle^{\M(r,n)}_{1,1,0},
\end{equation}
where we put the $n=0$ term as $0$ as convention. The main result of this part is an explicit formula of the generating function.
\begin{theorem}\label{thm:closed-main}
\begin{equation}\label{eq:closed-main}
\begin{split}
\F^{(0)}_{r,H}(z) = - \frac{P(z)^r}{24} & \left[\left(\left(\frac{rh}{2}+\sum_{\alpha=1}^r a_\alpha\right)E_2(z) - \frac{rh}{2}E_3(z)\right) \left(\frac{h}{p} + \frac{rh}{p}L(z)+\sum_{1\leq \alpha<\beta \leq r}\mathcal C_{\alpha\beta}(z)\right) \right. \\
& \left. + \frac{1}{2}\sum_{1\leq \alpha<\beta \leq r} \left((a_\alpha+a_\beta+h) z\frac{d\mathcal C_{\alpha\beta}(z)}{dz}\right)\right],
\end{split}
\end{equation}
where $h=\epsilon_1+\epsilon_2$ and $p=\epsilon_1\epsilon_2$ from (\ref{eqn:parameter_expressions}). 
\end{theorem}
\begin{example}[Charge one]\label{cor:n-one}
For $n=1$, the $z^1$-coefficient of (\ref{eq:closed-main}) gives
\begin{equation}\label{eq:n-one}
\left\langle H\right\rangle^{\M(r,1)}_{1,1,0}
 =-\frac1{24}\left[\frac hp\sum_{\alpha=1}^r a_\alpha
 +\sum_{1\leq \alpha<\beta \leq r}\frac{h^2+(a_\alpha+a_\beta)h}{h^2-(a_\alpha-a_\beta)^2}\right].
\end{equation}
\end{example}
\begin{example}[Rank one]\label{ex:rank1}
For $r=1$, the framing action is unnecessary, so we set $a_1=0$. Then, (\ref{eq:closed-main}) recovers the formula for degree-zero invariants in \cite[Section~1, especially equation~(1.3)]{ILPT}
\begin{equation}\label{eq:rank-one-effective}
\F^{(0)}_{1,H}(z)=-\frac{h^2}{48p}P(z)(1+L(z))(E_2(z)-E_3(z)).
\end{equation}
\end{example}

\section{Fixed point localization and sum over partitions}
By (the torus-equivariant version of) \cite[Proposition 2.7]{Getzler}, we have
\begin{equation}\label{eqn:degree_0_inv}
\langle H\rangle_{1,1,0}^{\M(r,n)}=-\frac{1}{24}\int^\T_{\M(r,n)}H\cup c_{\text{top}-1}(T_{\M(r,n)}),
\end{equation}
We compute the $\T$-equivariant integral on the RHS via fixed point localization. Recall from Section \ref{sec:fixed-pts} that the $\T$-fixed points in $\M(r,n)$ are indexed by $r$-colored partitions of $n$. Hence, the fixed point localization formula gives
\begin{equation} \label{eq:fixed-main}
\F^{(0)}_{r,H}(z)
 = -\frac{1}{24} \sum_{\boldsymbol\lambda =(\lambda^1,\ldots,\lambda^r) \in \mathcal{P}_r} z^{|\boldsymbol\lambda|} \cdot H|_{\boldsymbol\lambda} \cdot \Theta|_{\boldsymbol\lambda},
\end{equation}
where $\mathcal{P}_r$ denotes the set of $r$-colored partitions, $H|_{\boldsymbol\lambda}$ is the restriction of $H$ at the fixed point indexed by $\boldsymbol\lambda$, and $\Theta|_{\boldsymbol\lambda}$ is the sum of reciprocals of the tangent weights there. In Section \ref{sec:fixed-pts}, we computed
\begin{equation} \label{eq:H-restriction}
H|_{\boldsymbol\lambda}
=\sum_{\alpha=1}^r K_{a_\alpha}(\lambda^\alpha),
\end{equation}
where for a partition $\lambda$ and a $\T$-weight $a$, we denote
\begin{equation}\label{eq:K-a}
K_a(\lambda) := a|\lambda|-\epsilon _1n(\lambda)-\epsilon _2n(\lambda')
\end{equation}
Moreover, by (\ref{eqn:full_tangent_weights}),
\begin{equation} \label{equation: Theta restriction}
\Theta|_{\boldsymbol\lambda} = \sum_{\alpha=1}^r S(\lambda^\alpha) + \sum_{1 \leq\alpha<\beta\leq r} \Phi_{a_\alpha-a_\beta}(\lambda^\alpha,\lambda^\beta),
\end{equation}
where for a partition $\lambda$, we denote
\begin{equation}\label{eq:s-lambda}
S(\lambda) := \sum_{s\in\lambda}
 \left(
 \frac1{-L_\lambda(s)\epsilon _1+(A_\lambda(s)+1)\epsilon _2}
 +\frac1{(L_\lambda(s)+1)\epsilon _1-A_\lambda(s)\epsilon _2}
 \right),
\end{equation}
and for partitions $\lambda,\mu$ and $\T$-weight $a$, we denote
\begin{equation}\label{eq:Phi-delta}
\begin{split}
 \Phi_a(\lambda,\mu) := & 
 \sum_{s\in\lambda}
 \left(
 \frac{1}{-a-L_\mu(s)\epsilon _1+(A_\lambda(s)+1)\epsilon _2}
 +\frac{1}{a+(L_\mu(s)+1)\epsilon _1-A_\lambda(s)\epsilon _2}
 \right)\\
 + & \sum_{s\in\mu}
 \left(
 \frac{1}{-a+(L_\lambda(s)+1)\epsilon _1-A_\mu(s)\epsilon _2}
 +\frac{1}{a-L_\lambda(s)\epsilon _1+(A_\mu(s)+1)\epsilon _2}
 \right).   
 \end{split}
\end{equation}
In this way, the proof of Theorem \ref{thm:closed-main} boils down to computing the above sum over $r$-colored partitions. 

Given a function $f=f(\boldsymbol\lambda)$ on $\mathcal{P}_r$, the space of $r$-colored partitions, define
\begin{equation}\label{eq:q-average}
\langle f \rangle_r := P(z)^{-r} \cdot \sum_{\boldsymbol\lambda \in \mathcal{P}_r}
 z^{|\boldsymbol\lambda|} \cdot f(\boldsymbol\lambda),
\end{equation}
the norm of $f$ on $\mathcal{P}_r$ under the measure $\mathfrak{m}_r(\boldsymbol\lambda) := P(z)^{-r} z^{|\boldsymbol\lambda|}$. Note that $\mathfrak{m}_r$ is the product of $r$-copies of the measure $\mathfrak{m}_1(\lambda) := P(z)^{-1}z^{|\lambda|}$ on $\mathcal{P}_1$, the space of ($1$-colored) partitions. Hence, substituting (\ref{eq:H-restriction}) and (\ref{equation: Theta restriction}) into (\ref{eq:fixed-main}), we have
\begin{equation}
\F^{(0)}_{r,H}(z) = - \frac{P(z)^r}{24} \langle H \cdot \Theta \rangle_r,
\end{equation}
where
\begin{equation}
\begin{split}
\langle H \cdot \Theta \rangle_r = & \sum_{\alpha=1}^r \langle S(\lambda) \cdot K_{a_\alpha}(\lambda) \rangle_1 + (r-1) \cdot \langle S(\lambda) \rangle_1 \cdot \sum_{\alpha=1}^r \langle K_{a_\alpha}(\lambda)\rangle_1 \\
 + & \sum_{1\leq \alpha<\beta \leq r} \left(\langle \Phi_{a_\alpha-a_\beta}(\lambda,\mu) \rangle_2 \cdot \sum_{\substack{1\leq \gamma\leq r\\ \gamma \neq \alpha,\beta}} \langle K_{a_\gamma}(\lambda)\rangle_1 \right) \\
 + & \sum_{1\leq \alpha<\beta \leq r} \langle \Phi_{a_\alpha-a_\beta}(\lambda,\mu) \cdot (K_{a_\alpha}(\lambda)+K_{a_\beta}(\mu)) \rangle_2.
\end{split}
\end{equation}
Each of the norms appearing above, which are over either $\mathcal{P}_1$ or $\mathcal{P}_2$, can be computed explicitly. We write the results in the five lemmas below, from which it is not hard to check that Theorem \ref{thm:closed-main} follows.
\begin{lemma}
\begin{equation}\label{eq:mean-s}
\langle S(\lambda)\rangle_1=\frac{h}{p}L(z).
\end{equation}
\end{lemma}
\begin{proof}
\cite[Lemma~1]{RZ} gives
\begin{equation}
\langle M_{A,L}(\lambda) \rangle_1
 =\frac{z^{A+L+1}}{1-z^{A+L+1}}
\end{equation}
for 
\begin{equation}
M_{A,L}(\lambda) := \#\{s\in\lambda \ |\ A_\lambda(s) = A, L_\lambda(s) = L\}.
\end{equation}
Hence,
\begin{equation}
\begin{split}
\langle S(\lambda)\rangle_1 = & \sum_{A,L\geq 0} \left( \langle M_{A,L}(\lambda) \rangle_1 \cdot \left( \frac{1}{-L\epsilon_1+(A+1)\epsilon_2} + \frac{1}{(L+1)\epsilon_1-A\epsilon_2} \right) \right) \\
= & \sum_{B \geq 0} \left(\frac{z^{B+1}}{1-z^{B+1}} \cdot \sum_{\substack{A+L=B \\ A,L\geq 0}} \left( \frac{1}{-L\epsilon_1+(A+1)\epsilon_2} + \frac{1}{(L+1)\epsilon_1-A\epsilon_2} \right)\right) \\
= & \sum_{B \geq 0} \left(\frac{z^{B+1}}{1-z^{B+1}} \cdot \left( \frac{1}{(B+1)\epsilon_1} + \frac1{(B+1)\epsilon_2} \right)\right) = \frac{h}{p} L(z).
\end{split}
\end{equation}
\end{proof}
\begin{lemma} \label{lemma:mean-K}
For any $\T$-weight $a$,
\begin{equation}\label{eq:mean-K}
 \langle K_a(\lambda)\rangle_1 = \left(a+\frac{h}{2}\right)E_2(z)-\frac{h}{2}E_3(z). 
\end{equation}
\end{lemma}
\begin{proof}
For $\ell\geq 1$, denote by $m_{\ell}(\lambda)$ the function on the space of partitions $\mathcal{P}_1$ recording the number of times $\ell$ appears in the partition $\lambda$. Then, it is not hard to compute
\begin{equation}
\begin{split}
\langle m_{\ell}(\lambda)\rangle_1 = & P(z)^{-1} \cdot \sum_{\lambda} \left( z^{\sum_{k\geq 1} k \cdot m_k(\lambda)} \cdot m_{\ell}(\lambda) \right) \\
= & P(z)^{-1} \cdot \sum_{m_{\ell}\geq 0} \left(m_{\ell}z^{\ell \cdot m_{\ell}}\right) \cdot \prod_{k\geq 1,k\neq \ell} \left( \sum_{m_k\geq 0} z^{ k \cdot m_k}\right) \\
= & P(z)^{-1} \cdot \left(\frac{1}{\ell} \cdot z\frac{d}{dz}\frac{1}{1-z^\ell}\right)\cdot \prod_{k\geq 1,k\neq \ell} \left(\frac{1}{1-z^k} \right) = \frac{z^\ell}{1-z^\ell}.
\end{split}
\end{equation}
Then, since for $\lambda = (\lambda_1\geq \lambda_2 \geq \cdots)$, $|\lambda| = \sum_i\lambda_i = \sum_{\ell\geq 1} \left(m_{\ell}(\lambda) \cdot \ell \right)$ and $n(\lambda') = \sum_i \binom{\lambda_i}{2} = \sum_{\ell\geq 1} \left(m_{\ell}(\lambda) \cdot \binom{\ell}{2}\right)$, we have
\begin{equation}
\begin{split}
\langle K_a(\lambda)\rangle_1 = & \ a \langle |\lambda|\rangle_1 - \epsilon_1 \langle n(\lambda)\rangle_1 - \epsilon_2 \langle n(\lambda')\rangle_1 \\
= & \ a \langle |\lambda|\rangle_1 - (\epsilon_1+\epsilon_2)\langle n(\lambda')\rangle_1 \\
= &\left(a+\frac{h}{2}\right)E_2(z)-\frac{h}{2}E_3(z),
\end{split}
\end{equation}
where the second equality holds as the measure $\mathfrak{m}_1$ is invariant under transposition of partitions.
\end{proof}
\begin{lemma} \label{lem:same-covariance}
For any $\T$-weight $a$,
\begin{equation}
\langle K_{a}(\lambda) \cdot S(\lambda) \rangle_1 = \frac{h}{p} (1+L(z)) \left(\left(a+\frac{h}{2}\right)E_2(z)-\frac{h}{2}E_3(z)\right).
\end{equation}
\end{lemma}
\begin{proof}
By definition, 
\begin{equation}
\langle K_{a}(\lambda) \cdot S(\lambda) \rangle_1 = a \langle |\lambda| \cdot S(\lambda) \rangle_1 + \langle K_{0}(\lambda) \cdot S(\lambda) \rangle_1.
\end{equation}
The first term on the RHS can be computed from (\ref{eq:mean-s}) since 
\begin{equation}
P(z) \cdot \langle |\lambda| \cdot S(\lambda) \rangle_1 = z\frac{d}{dz}\left(P(z) \cdot \langle S(\lambda) \rangle_1\right).
\end{equation}
The second term on the RHS is now independent from the $\T$-action on the framing space, and is in fact exactly the fixed point localization of \cite[Equation (1.3)]{ILPT}. We re-write it in our notation as follows
\begin{equation}
\langle K_{0}(\lambda) \cdot S(\lambda) \rangle_1 = \frac{h^2}{2p}(E_2(z)-E_3(z))(1+L(z)).
\end{equation}
The lemma then follows simply from combining the above two equations.
\end{proof}
\begin{lemma} \label{lemma: norm Phi a}
For any $\T$-weight $a$,
\begin{equation}\label{eq:mean-Phi}
 \langle\Phi_a(\lambda,\mu)\rangle_2
 =\mathcal{C}_{a}(z):=
 \sum_{i,j\geq1}z^{ij}
 \left(\frac{1}{i\epsilon_1+j\epsilon_2-a}+\frac{1}{i\epsilon_1+j\epsilon_2+a}\right).
\end{equation}
\end{lemma}
\begin{lemma} \label{lem:pointed-two-partition}
For any $\T$-weights $a$ and $b$,
\begin{equation} \label{eq:cross-covariance}
\langle \Phi_{a-b}(\lambda,\mu) \cdot (K_a(\lambda)+K_b(\mu))\rangle_2
= \frac{a+b+h}{2} \cdot z\frac{d\mathcal{C}_{a-b}(z)}{dz} + \left( (a+b+h)E_2(z) - hE_3(z)\right) \cdot \mathcal{C}_{a-b}(z),
\end{equation}
where the function $\mathcal{C}$ is defined in Lemma \ref{lemma: norm Phi a}.
\end{lemma}
The former three lemmas involve only one color of the partition, and thus reduce essentially to the $r=1$ case, namely the case of Hilbert scheme $\M(1,n) = \mathrm{Hilb}^n(\CC^2)$, up to easy terms involving the equivariant parameter $a$. The latter two lemmas, on the other hand, are across different colors and are thus genuinely new for the $r>1$ case. The rest of this part is dedicated to proving these two lemmas.

\section{Proof of Lemma \ref{lemma: norm Phi a}}
\subsection{K-theoretic model}
We start by studying a K-theoretic version of $\Phi_a(\lambda,\mu)$ which is more readily understood in the literature. Denote respectively by $t_1, t_2, A$ the $\T$-characters with weights $\epsilon_1,\epsilon_2, a$. In other words, $c_1(t_1)=\epsilon_1, c_1(t_2)=\epsilon_2, c_1(A) = a$ on $B\T$. Define
\begin{equation}\label{eqn:char_inv}
R: K(B\T)  \longrightarrow  \mathrm{Frac}(H^*(B\T))
\end{equation}
as the map linearly extending $R(L) = \frac{1}{c_1(L)}$ for any line bundle $L$. Then,
\begin{equation}
\Phi_a(\lambda,\mu) = R(A^{-1}E_{\lambda\mu}+AE_{\mu\lambda})
\end{equation}
for 
\begin{equation}
E_{\lambda\mu} := \sum_{s\in\lambda} t_1^{-L_\mu(s)} t_2^{A_\lambda(s)+1} + \sum_{s\in\mu} t_1^{L_\lambda(s)+1} t_2^{-A_\mu(s)}
\end{equation}
which also satisfies
\begin{equation} \label{equation: E in terms of V}
E_{\lambda\mu} = -(1-t_1)(1-t_2)V_\lambda^\vee V_\mu + t_1t_2 V_{\lambda}^\vee  + V_\mu
\end{equation}
for
\begin{equation}
V_\lambda := \sum_{s=(i,j)\in\lambda} t_1^{1-i}t_2^{1-j}.
\end{equation}
(\ref{equation: E in terms of V}) is exactly \cite[Lemma 6]{co} (where $z_1,z_2$ are respectively $t_2,t_1$ in our notation), and has already been used implicitly in Section \ref{sec:tangent-weights} when we describe the tangent weights. Generally speaking, $\T$-representations naturally appearing in $\M(r,n)$ usually have appeared in $\M(1,n)= \mathrm{Hilb}^n(\CC^2)$, once the action on the framing space is dropped. Under such correspondence, cross-color computation of the former boil down to cross-fixed-point computation of the latter. The appearance of $E_{\lambda\mu}$ is nothing more than another example of this observation.

As preparation for computing $\langle\Phi_a(\lambda,\mu)\rangle_2$, we first compute $\langle V_\lambda \rangle_1$ and $\langle E_{\lambda\mu} \rangle_2$, where $V_\lambda$ and $E_{\lambda\mu}$ are regarded respectively as functions on $\mathcal{P}_1$ and $\mathcal{P}_2$ in the most natural way. 
\subsection{Norm of $V_\lambda$}
Consider the function $I_{\lambda_r=k}$ on $\mathcal{P}_1$ defined on $\lambda = (\lambda_1\geq \lambda_2\geq \cdots)$ by
\begin{equation}
I_{\lambda_r=k}(\lambda) = 
\begin{cases}
1 & \quad \lambda_r=k \\
0 & \quad \lambda_r\neq k \text{ or } \lambda_r \text{ does not exist}.
\end{cases}
\end{equation}
Given any $\lambda$ such that $I_{\lambda_r=k}(\lambda) = 1$, the Young diagram of $\lambda$ splits exactly into two pieces after the corner $r\times k$ rectangle is removed. Denote by $\lambda^1$ and $\lambda^2$ the two partitions associated to the two pieces. Then, the length of $\lambda^1$ is at most $r-1$, and numbers in $\lambda^2$ are at most $k$. Conversely, given any such $\lambda^1$ and $\lambda^2$, we can obtain an $\lambda$ with $I_{\lambda_r=k}(\lambda) = 1$ by patching together the two Young diagrams and the rectangle. 

For simplicity of notation, we denote from now
\begin{equation}
(x;z)_d := \prod_{l=1}^d (1-xz^{d-1}),
\end{equation}
for $d\in \mathbb{Z}_{\geq 0}$. Then,
\begin{equation}
P(z)\langle I_{\lambda_r=k}(\lambda) \rangle_1  = \sum_{\lambda^1\text{ as above}} \sum_{\lambda^2\text{ as above}} z^{kr + |\lambda^1| + |\lambda^2|} = \frac{z^{kr}}{(z;z)_{r-1}(z;z)_k}. 
\end{equation}
By construction
\begin{equation}
V_\lambda = \sum_{r\geq 1}\sum_{k\geq 1} \left( I_{\lambda_r=k}(\lambda) \cdot t_1^{1-r} \cdot (1+\cdots+t_2^{1-k}) \right),
\end{equation}
so
\begin{equation} \label{equation: norm of Vlambda middle step 1}
\begin{split}
\langle V_\lambda \rangle_1 = & P(z)^{-1}\sum_{r\geq 1}\sum_{k\geq 1} \frac{z^{kr} \cdot t_1^{1-r} \cdot (1+\cdots+t_2^{1-k})}{(z;z)_{r-1}(z;z)_k} \\
= & P(z)^{-1}\sum_{r\geq 1} \frac{t_1^{1-r}}{1-t_2^{-1}}\left(\frac{1}{(z;z)_{r-1}(z^r;z)_{\infty}} - \frac{1}{(z;z)_{r-1}(t_2^{-1}z^r;z)_{\infty}}\right) \\
= & \frac{1}{(1-t_2^{-1})}\sum_{r\geq 1}t_1^{1-r} - \frac{P(z)^{-1}}{(1-t_2^{-1}) (t_2^{-1}z;z)_{\infty}} \sum_{r\geq 1} \left(\frac{(t_2^{-1}z;z)_{r-1}}{(z;z)_{r-1}}t_1^{1-r}\right) \\
= & \frac{1}{(1-t_1^{-1})(1-t_2^{-1})} - \frac{(z;z)_{\infty}}{(1-t_2^{-1}) (t_2^{-1}z;z)_{\infty}} \frac{(t_1^{-1}t_2^{-1}z;z)_{\infty}}{(t_1^{-1}z;z)_{\infty}} \\
= & \frac{1}{(1-t_1^{-1})(1-t_2^{-1})} \left(1-\frac{(t_1^{-1}t_2^{-1}z;z)_{\infty}(z;z)_{\infty}}{(t_1^{-1}z;z)_{\infty}(t_2^{-1}z;z)_{\infty}}\right)
\quad \in \CC[\![t_1^{-1},t_2^{-1},z]\!],
\end{split}
\end{equation}
where for the second equality we use the identity 
\begin{equation}
\sum_{k\geq 0} \frac{\lambda^k}{(z;z)_k} = \frac{1}{(\lambda;z)_\infty} \quad \in \CC[\lambda][\![z]\!]
\end{equation}
and for the fourth equality we use \cite[Equation (1.3.2)]{GR}
\begin{equation}
\sum_{n\geq 0} \left(\frac{(a;q)_{n}}{(q;q)_{n}}z^{n}\right) = \frac{(az;q)_\infty}{(z;q)_\infty} \quad \in \CC[\![q,z]\!]
\end{equation}
by setting $n=r-1, a=t_2^{-1}z, q = z, z = t_1^{-1}$. The latter holds for analytic functions for $|z|<1$ and $|q|<1$, and thus for power series in $q$ and $z$ as well. (\ref{equation: norm of Vlambda middle step 1}) implies in particular that for a fixed power of $z$, the coefficient of the RHS, which a priori is a power series in $t_1^{-1}$ and $t_2^{-1}$, has only finite terms because of the same property of the LHS by construction. That is to say, $1-t_1^{-1}$ and $1-t_2^{-1}$ appear as divisors of the numerator, and (\ref{equation: norm of Vlambda middle step 1}) holds in $\CC[t_1^{\pm 1},t_2^{\pm 1}][\![z]\!]$. We write 
\begin{equation}
\langle V_\lambda \rangle_1 = \frac{1}{(1-t_1^{-1})(1-t_2^{-1})} \left(1-\frac{(t_1^{-1}t_2^{-1}z;z)_{\infty}(z;z)_{\infty}}{(t_1^{-1}z;z)_{\infty}(t_2^{-1}z;z)_{\infty}}\right)
\in \CC[t_1^{\pm 1},t_2^{\pm 1}][\![z]\!].
\end{equation}
By setting $t_1\mapsto t_1^{-1}$ and $t_2\mapsto t_2^{-1}$ in $\langle V_\lambda \rangle_1$, we also have
\begin{equation}
\langle V_\lambda^\vee \rangle_1 = \frac{t_1^{-1}t_2^{-1}}{(1-t_1^{-1})(1-t_2^{-1})} \left(1-\frac{(t_1t_2z;z)_{\infty}(z;z)_{\infty}}{(t_1z;z)_{\infty}(t_2z;z)_{\infty}}\right)
\in \CC[t_1^{\pm 1},t_2^{\pm 1}][\![z]\!].
\end{equation}

\subsection{Norm of $E_{\lambda\mu}$}
By (\ref{equation: E in terms of V}),
\begin{equation}
\begin{split}
\langle E_{\lambda\mu} \rangle_2 = &  -(1-t_1)(1-t_2) \cdot \langle V_\lambda^\vee \rangle_1 \langle V_\mu \rangle_1 + t_1t_2 \langle V_\lambda^\vee \rangle_1 + \langle V_\mu \rangle_1 \\
= & \frac{1}{(1-t_1^{-1})(1-t_2^{-1})}\left(1 - \frac{(t_1t_2z;z)_{\infty}(z;z)_{\infty}}{(t_1z;z)_{\infty}(t_2z;z)_{\infty}}\frac{(t_1^{-1}t_2^{-1}z;z)_{\infty}(z;z)_{\infty}}{(t_1^{-1}z;z)_{\infty}(t_2^{-1}z;z)_{\infty}}\right) \, \in \CC[t_1^{\pm 1},t_2^{\pm 1}][\![z]\!].
\end{split}
\end{equation}
More cancellation will happen once we apply $R$ in (\ref{eqn:char_inv}). Indeed, for any $\T$-character $\lambda$, $R(\lambda+\lambda^{-1}) = 0$ by construction, and thus in our final result, we will only need 
\begin{equation} \label{equation: Elambdamu middle step 1}
\begin{split}
& \langle E_{\lambda\mu} \rangle_2 - \langle E_{\lambda\mu} \rangle_2|_{t_1\mapsto t_1^{-1}, t_2\mapsto t_2^{-1}} \\
= & \frac{1-t_1^{-1}t_2^{-1}}{(1-t_1^{-1})(1-t_2^{-1})} \cdot \left( 1 - \frac{(t_1t_2z;z)_{\infty}(z;z)_{\infty}}{(t_1z;z)_{\infty}(t_2z;z)_{\infty}}\frac{(t_1^{-1}t_2^{-1}z;z)_{\infty}(z;z)_{\infty}}{(t_1^{-1}z;z)_{\infty}(t_2^{-1}z;z)_{\infty}} \right) \\
= & \frac{1-t_1^{-1}t_2^{-1}}{(1-t_1^{-1})(1-t_2^{-1})}  - \sum_{i\geq 0} \frac{t_1^{-i}}{1-t_2^{-1}z^{i}} - \sum_{i<0} \frac{t_1^{-i}t_2z^{-i}}{1-t_2z^{-i}} \\
= & \sum_{i,j\geq 1} z^{ij} (t_1^it_2^j-t_1^{-i}t_2^{-j}) \, \in \CC(\!(t_1^{-1},t_2^{-1})\!)[\![z]\!],
\end{split}
\end{equation}
where for the second step we use Ramanujan's ${}_1\psi_1$ summation, \cite[Equation~(5.2.1)]{GR}
\begin{equation}
\frac{(q;q)_\infty(b/a;q)_\infty(az;q)_\infty(q/(az);q)_\infty}{(b;q)_\infty(q/a;q)_\infty(z;q)_\infty(b/(az);q)_\infty} = \sum_{i\geq 0} \frac{(a;q)_i}{(b;q)_i}z^i + \sum_{i<0} \frac{(q/b;q)_{-i}}{(q/a;q)_{-i}} \left(\frac{b}{az}\right)^{-i},
\end{equation}
by setting $a = t_2^{-1}, b=t_2^{-1}z, q=z, z=t_1^{-1}$,
which holds for analytic functions in $|z|<|t_1^{-1}|<1$ and thus for elements in $\CC(\!(t_2^{-1})\!)[\![t_1z,z,t_1^{-1}]\!]\subset\CC(\!(t_1^{-1},t_2^{-1})\!)[\![z]\!]$. Again, (\ref{equation: Elambdamu middle step 1}) holds a priori in $\CC(\!(t_1^{-1},t_2^{-1})\!)[\![z]\!]$ but in fact in $\CC[t_1^{\pm 1}, t_2^{\pm 1}][\![z]\!]$ as well, even though the middle steps involve genuine Laurent series.

\subsection{Image under $R$} It is not hard to see that $R$ commutes with taking norms, so
\begin{equation}
\begin{split}
\langle\Phi_a(\lambda,\mu)\rangle_2 = & R((A^{-1}+A)\langle E_{\lambda\mu}\rangle_2) \\
= & \frac{1}{2}R\left((A^{-1}+A)\left(\langle E_{\lambda\mu} \rangle_2 - \langle E_{\lambda\mu} \rangle_2|_{t_1\mapsto t_1^{-1}, t_2\mapsto t_2^{-1}}\right)\right) = \mathcal{C}_a(z).
\end{split}
\end{equation}
Here in the second step, we use again the fact that $R(\lambda+\lambda^{-1}) = 0$ for any $\T$-character $\lambda$.

\section{Proof of Lemma \ref{lem:pointed-two-partition}}
The quantity $\langle \Phi_{a-b}(\lambda,\mu) \cdot (K_a(\lambda)+K_b(\mu))\rangle_2$ is a power series in $z$ with coefficients being rational functions in $a,b,\epsilon_1,\epsilon_2$. Moreover, by (\ref{eq:Phi-delta}), the only possible denominators are of the form $$a-b- u\epsilon_1-v\epsilon_2, \text{ for } u,v\in\mathbb{Z}.$$ 
Under the change of variables 
\begin{equation} \label{equation: change of var ab to cdelta}
\delta = a-b, \quad c=a+b,
\end{equation}
they give at most simple poles at $\delta = u\epsilon_1+v\epsilon_2$. Our idea is to understand the residues at these simple poles, which then determine the rational functions given asymptotic behavior at $\delta \rightarrow \infty$.

\subsection{Boundary of Young diagram}
Recall that for a partition $\lambda=(\lambda_1\ge\lambda_2\ge\cdots)$, we place its Young diagram $\mathsf{D}_\lambda$ in the first quadrant by setting
\begin{equation}
\mathsf{D}_\lambda=\{(i,j) \ |\ i\ge1, 1\le j\le\lambda_i\}
\end{equation}
where $(i,j)$ is understood as the box with vertices $(i,j), (i-1,j), (i,j-1), (i-1,j-1)$. 

We define $\Gamma_\lambda$ to be the boundary of the first quadrant with $\mathsf{D}_\lambda$ removed, and endow it with the direction going from lower-right to upper-left. In algebraic terms, if $\lambda=(\lambda_1\geq\cdots\geq\lambda_l)$ and we denote as convention $\lambda_{m} = 0$ for all $m > l$, then $\Gamma_\lambda$ is the broken line consisting of
\begin{itemize}
    \item segments $(i,\lambda_i)\rightarrow(i-1,\lambda_i)$ going left for $i\in\mathbb{Z}_{\geq 1}$, and 
    \item segments $(\lambda_i',i-1)\rightarrow(\lambda_i',i)$ going up for $i\in\mathbb{Z}_{\geq 1}$, where $\lambda'$ is the transpose of $\lambda$.
\end{itemize}

\subsection{Poles of $\Phi_\delta(\lambda,\mu)$ via intersection jets}
Given partitions $\lambda$ and $\mu$, $\Phi_\delta(\lambda,\mu)$ can be re-written in terms of certain intersection configurations of translated boundaries $\Gamma_\lambda$ and $\Gamma_\mu$, which we now explain.

Fix integers $u,v\in\mathbb{Z}$, and denote by $$\Gamma_{\mu}^{(u,v)}$$ the directed broken line translated along $(u,v)$. We define an intersection jet of $\Gamma_\lambda$ and $\Gamma_{\mu}^{(u,v)}$ to be the data of an intersection point $P$ along with two adjacent directed segments $S_1\subset \Gamma_\lambda$ and $S_2\subset \Gamma_{\mu}^{(u,v)}$ of length $1$. At a fixed $P\in \Gamma_\lambda \cap \Gamma_{\mu}^{(u,v)}$, there are $16$ possible intersection jet as each of the two segments can be incoming/outgoing and up/left. The pole of $\Phi_\delta(\lambda,\mu)$ at $\delta = u\epsilon_1+v\epsilon_2$ are related to the following $4$ specific types of intersection jets.

A {\em type-A intersection jet} is one where $S_1$ is outgoing left and $S_2$ is outgoing up. This implies that the intersection point $P$ has coordinates $(i,\lambda_i) = (\mu_j'+u, j-1+v)$ for some $i,j\geq 1$, or equivalently that there exists a box $s = (i,j)$ satisfying $u=-L_\mu(s), v = A_\lambda(s)+1$. Such $s$ contributes to the pole if and only if $s\in\lambda$, or equivalently $v = \lambda_i-j+1>0$, and in such case the contribution is $1/(-\delta+u\epsilon_1+v\epsilon_2)$. We denote by $$[\Gamma_\lambda \cap \Gamma_{\mu}^{(u,v)}]_A$$ the number of type-A intersection jets.

A {\em type-B intersection jet} is one where $S_1$ is incoming left and $S_2$ is incoming up. This implies that the intersection point $P$ has coordinates $(i-1,\lambda_i) = (\mu_j'+u, j+v)$ for some $i,j\geq 1$, or equivalently that there exists a box $s = (i,j)$ satisfying $u=-L_\mu(s)-1, v = A_\lambda(s)$. Such $s$ contributes to the pole if and only if $s\in\lambda$, or equivalently $v = \lambda_i-j\geq 0$, and in such case the contribution is $1/(\delta-u\epsilon_1-v\epsilon_2)$. We denote by $$[\Gamma_\lambda \cap \Gamma_{\mu}^{(u,v)}]_B$$ the number of type-B intersection jets.

A {\em type-C intersection jet} is one where $S_1$ is incoming up and $S_2$ is incoming left. This implies that the intersection point $P$ has coordinates $(\lambda_i',i) = (j-1+u, \mu_j+v)$ for some $i,j\geq 1$, or equivalently that there exists a box $s = (i,j)$ satisfying $u=L_\lambda(s)+1, v = -A_\mu(s)$. Such $s$ contributes to the pole if and only if $s\in\mu$, or equivalently $v = i-\mu_j \leq 0$, and in such case the contribution is $1/(-\delta+u\epsilon_1+v\epsilon_2)$. We denote by $$[\Gamma_\lambda \cap \Gamma_{\mu}^{(u,v)}]_C$$ the number of type-C intersection jets.

A {\em type-D intersection jet} is one where $S_1$ is outgoing up and $S_2$ is outgoing left. This implies that the intersection point $P$ has coordinates $(\lambda_i',i-1) = (j+u, \mu_j+v)$ for some $i,j\geq 1$, or equivalently that there exists a box $s = (i,j)$ satisfying $u=L_\lambda(s), v = -A_\mu(s)-1$. Such $s$ contributes to the pole if and only if $s\in\mu$, or equivalently $v = i-\mu_j-1 < 0$, and in such case the contribution is $1/(\delta-u\epsilon_1-v\epsilon_2)$. We denote by $$[\Gamma_\lambda \cap \Gamma_{\mu}^{(u,v)}]_D$$ the number of type-D intersection jets.

We collect the information on the four types in the chart below. Note that one intersection point may contribute several intersection jets.
\begin{center}
\small
\setlength{\tabcolsep}{4pt}
\renewcommand{\arraystretch}{1.4}
\begin{tabular}{c c c c c}
\hline
Type & $S_1\subset\Gamma_\lambda$ & $S_2\subset\Gamma_{\mu}^{(u,v)}$ & In $\Phi_\delta(\lambda,\mu)$ & Condition \\
\hline
A
& outgoing left
& outgoing up
& $1/(-\delta+u\epsilon_1+v\epsilon_2)$ & $v>0$ \\

B
& incoming left
& incoming up
& $1/(\delta-u\epsilon_1-v\epsilon_2)$ & $v\geq 0$ \\

C
& incoming up
& incoming left
& $1/(-\delta+u\epsilon_1+v\epsilon_2)$ & $v\leq 0$ \\

D
& outgoing up
& outgoing left
& $1/(\delta-u\epsilon_1-v\epsilon_2)$ & $v<0$ \\
\hline
\end{tabular}
\end{center}
In this way, we have
\begin{equation}
\begin{split}
\Phi_\delta(\lambda,\mu) = & \sum_{u\in\mathbb{Z}, v\in\mathbb{Z}_{>0}} \frac{-[\Gamma_\lambda \cap \Gamma_{\mu}^{(u,v)}]_A+[\Gamma_\lambda \cap \Gamma_{\mu}^{(u,v)}]_B}{\delta-u\epsilon_1-v\epsilon_2} \\
+ & \sum_{u\in\mathbb{Z}, v=0} \frac{[\Gamma_\lambda \cap \Gamma_{\mu}^{(u,v)}]_B-[\Gamma_\lambda \cap \Gamma_{\mu}^{(u,v)}]_C}{\delta-u\epsilon_1-v\epsilon_2} \\ 
+ & \sum_{u\in\mathbb{Z}, v\in\mathbb{Z}_{<0}} \frac{-[\Gamma_\lambda \cap \Gamma_{\mu}^{(u,v)}]_C+[\Gamma_\lambda \cap \Gamma_{\mu}^{(u,v)}]_D}{\delta-u\epsilon_1-v\epsilon_2},
\end{split}
\end{equation}
and thus
\begin{equation} \label{equation: proof lemma cov Phi K middle step 1}
\begin{split}
& P^2(z)\langle \Phi_{a-b}(\lambda,\mu) \cdot (K_a(\lambda)+K_b(\mu))\rangle_2 \\
= & \sum_{u\in\mathbb{Z}, v\in\mathbb{Z}_{>0}}\frac{1}{\delta-u\epsilon_1-v\epsilon_2} \sum_{\lambda,\mu} z^{|\lambda|+|\mu|} \left(-[\Gamma_\lambda \cap \Gamma_{\mu}^{(u,v)}]_A+[\Gamma_\lambda \cap \Gamma_{\mu}^{(u,v)}]_B\right) \left(K_{\frac{c+\delta}{2}}(\lambda)+K_{\frac{c-\delta}{2}}(\mu)\right) \\
+ & \sum_{u\in\mathbb{Z}, v=0} \frac{1}{\delta-u\epsilon_1-v\epsilon_2} \sum_{\lambda,\mu} z^{|\lambda|+|\mu|} \left([\Gamma_\lambda \cap \Gamma_{\mu}^{(u,v)}]_B-[\Gamma_\lambda \cap \Gamma_{\mu}^{(u,v)}]_C\right) \left(K_{\frac{c+\delta}{2}}(\lambda)+K_{\frac{c-\delta}{2}}(\mu)\right) \\
+ & \sum_{u\in\mathbb{Z}, v\in\mathbb{Z}_{<0}} \frac{1}{\delta-u\epsilon_1-v\epsilon_2} \sum_{\lambda,\mu} z^{|\lambda|+|\mu|} \left(-[\Gamma_\lambda \cap \Gamma_{\mu}^{(u,v)}]_C+[\Gamma_\lambda \cap \Gamma_{\mu}^{(u,v)}]_D\right) \left(K_{\frac{c+\delta}{2}}(\lambda)+K_{\frac{c-\delta}{2}}(\mu)\right),
\end{split}
\end{equation}
where we use the change of variables (\ref{equation: change of var ab to cdelta}) for the RHS.

\subsection{Cancellation of poles}
In the RHS of (\ref{equation: proof lemma cov Phi K middle step 1}), for each fixed $(u,v)$, cancellation happens in the associated sum over $\lambda,\mu$, which we describe separately for the three cases $v>0, v<0, v=0$ below. 

\begin{figure}[htbp]
  \centering
  \includegraphics[width=0.7\linewidth]{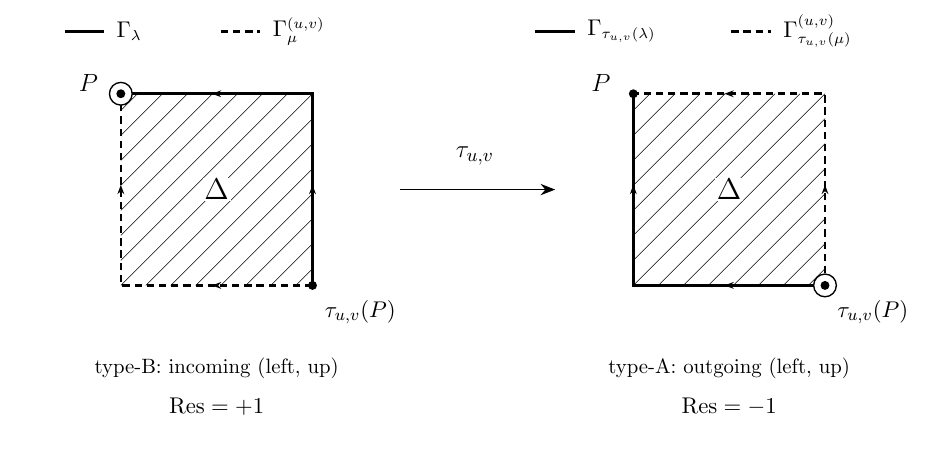}
  \caption{$\tau_{u,v}$ for $v>0$}
  \label{fig:gw1-excursion-switch}
\end{figure}

For $v>0$, we define an injective map
\begin{equation}
\tau_{u,v}:\ \left\{(\lambda,\mu,j) \ \middle|\ \begin{aligned}
j \text{ is a type-B intersection} \\
\text{jet of } \Gamma_\lambda \text{ and } \Gamma_\mu^{(u,v)}
\end{aligned}\right\} \ 
\longrightarrow \ 
\left\{(\lambda,\mu,j)\ \middle|\ \begin{aligned}
j \text{ is a type-A intersection} \\
\text{jet of } \Gamma_\lambda \text{ and } \Gamma_\mu^{(u,v)}
\end{aligned}\right\}.
\end{equation}
That is, given partitions $\lambda,\mu$ and type-B intersection jet $j = (P,S_1,S_2)$ of $\Gamma_\lambda$ and $\Gamma_\mu^{(u,v)}$, we define partitions $\tau_{u,v}(\lambda),\tau_{u,v}(\mu)$, and type-A intersection jet $$\tau_{u,v}(j) := (\tau_{u,v}(P),\tau_{u,v}(S_1),\tau_{u,v}(S_2)).$$ First, we define $\tau_{u,v}(P)$ as the intersection point of $\Gamma_\lambda$ and $\Gamma_\mu^{(u,v)}$ immediately preceding $P$. More precisely, if we regard both $\Gamma_\lambda$ and $\Gamma_\mu^{(u,v)}$ as parametrized by $t := y-x\in \mathbb{R}$, then their directions are both given by the direction where $t$ grows. Let $t_0:=t(P)$ and $t_{-1}:= \max\{t<t_0 \ |\ \Gamma_\lambda(t) = \Gamma_\mu^{(u,v)}(t)\}$. Then, we define $\tau_{u,v}(P):= \Gamma_\lambda(t_{-1}) = \Gamma_\mu^{(u,v)}(t_{-1})$. Such $t_{-1}$ always exists and is different from $t_0$, since $\mathrm{proj}_{y}\Gamma_\lambda(t_0-\epsilon) > \mathrm{proj}_{y}\Gamma_\mu^{(u,v)}(t_0-\epsilon)$ for $\epsilon$ small as $P$ is of type-B, and $\mathrm{proj}_{y}\Gamma_\lambda(-\infty) = 0 < v = \mathrm{proj}_{y}\Gamma_\mu^{(u,v)}(-\infty)$, by continuity of the function $\mathrm{proj}_{y}$ which takes the $y$-coordinate and thus determines a point along with $t = y-x$. Second, we define $\tau_{u,v}(\lambda)$ and $\tau_{u,v}(\mu)$ such that $\Gamma_{\tau_{u,v}(\lambda)}$ and $\Gamma_{\tau_{u,v}(\mu)}^{(u,v)}$ are given respectively by $\Gamma_{\lambda}$ and $\Gamma_\mu^{(u,v)}$ but with the piece between $P$ and $\tau_{u,v}(P)$ interchanged. Finally, we define $\tau_{u,v}(S_1)$ as the outgoing left segment at $\tau_{u,v}(P)$ in $\Gamma_\mu^{(u,v)}$ and thus in $\Gamma_{\tau_{u,v}(\lambda)}$, and $\tau_{u,v}(S_2)$ as the outgoing up segment at $\tau_{u,v}(P)$ in $\Gamma_\lambda$ and thus in $\Gamma_{\tau_{u,v}(\mu)}^{(u,v)}$. It is not hard to see that we obtain a type-A jet in this way, and that $\tau_{u,v}$ so defined is injective. See Figure \ref{fig:gw1-excursion-switch} for an illustration.

The pole contribution to the RHS of (\ref{equation: proof lemma cov Phi K middle step 1}) from $(\lambda,\mu,j)$ and $\tau_{u,v}(\lambda,\mu,j)$ cancel, since $|\lambda|+|\mu| = |\tau_{u,v}(\lambda)|+|\tau_{u,v}(\mu)|$ and $K_{\frac{c+\delta}{2}}(\lambda)+K_{\frac{c-\delta}{2}}(\mu) \equiv K_{\frac{c+\delta}{2}}(\tau_{u,v}(\lambda))+K_{\frac{c-\delta}{2}}(\tau_{u,v}(\mu)) \ \mathrm{mod} \ \delta-u\epsilon_1-v\epsilon_2$. The latter is because
\begin{equation}
\begin{split}
& K_{\frac{c+\delta}{2}}(\tau_{u,v}(\lambda))+K_{\frac{c-\delta}{2}}(\tau_{u,v}(\mu)) - K_{\frac{c+\delta}{2}}(\lambda)-K_{\frac{c-\delta}{2}}(\mu) \\
= & - \sum_{s=(i,j)\in\Delta} \left(\frac{c+\delta}{2} - (i-1)\epsilon_1-(j-1)\epsilon_2\right) + \sum_{s=(i,j)\in\Delta} \left(\frac{c-\delta}{2} - (i-u-1)\epsilon_1-(j-v-1)\epsilon_2\right) \\
= & - |\Delta| \cdot (\delta-u\epsilon_1-v\epsilon_2),
\end{split}
\end{equation}
where $\Delta$ is the region enclosed by $\Gamma_\lambda$ and $\Gamma_\mu^{(u,v)}$ between $P$ and $\tau_{u,v}(P)$.

\begin{figure}[htbp]
  \centering
\includegraphics[width=0.9\linewidth]{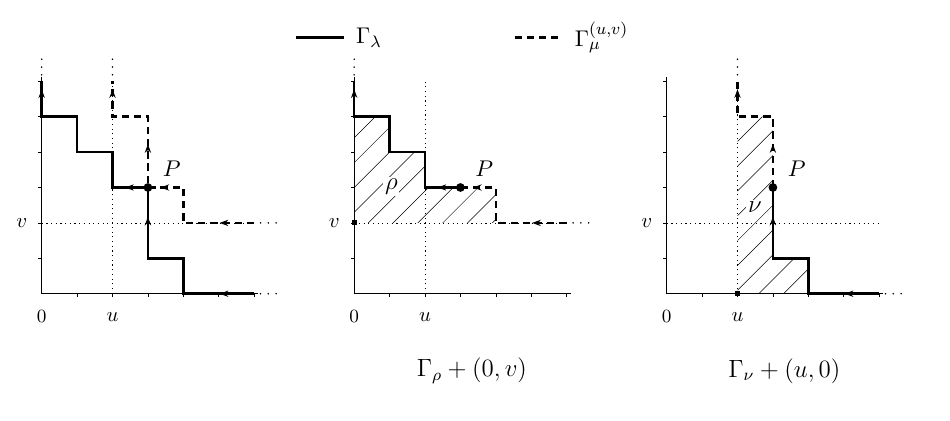}
\caption{$\tau_{u,v}$ for $v=0$}
  \label{fig:gw1-terminal-tail-exchange}
\end{figure}

$\tau_{u,v}$ is surjective for $u\leq 0$, in which case for each type-A intersection jet, there always exists an immediately following intersection point due to the continuity of $\mathrm{proj}_{x}$, which allows for building a map inverse to $\tau_{u,v}$. For $u>0$, however, $\tau_{u,v}$ is not surjective. A triple $(\lambda,\mu,j)$ is not in the image of $\tau_{u,v}$ if and only if $j$ is a type-A intersection jet of $\Gamma_\lambda$ and $\Gamma_\mu^{(u,v)}$ with no following intersection points. Denote by $P$ the intersection point where the jet $j$ is located. Let $\rho$ be the partition such that $\Gamma_\rho + (0,v)$ is given gluing the piece of $\Gamma_\mu^{(u,v)}$ to the lower-right of $P$ with the piece of $\Gamma_\lambda$ to the upper-left of $P$, and $\nu$ be the partition such that $\Gamma_\nu + (u,0)$ is given by the other two pieces. See Figure \ref{fig:gw1-terminal-tail-exchange} for an illustration.

This gives a bijection
\begin{equation}
\iota_{u,v}: \ 
\left\{(\lambda,\mu,j) \ \middle|\ \begin{aligned}
j \text{ is a type-A intersection} \\
\text{jet of } \Gamma_\lambda \text{ and } \Gamma_\mu^{(u,v)}
\end{aligned}\right\}\Big\backslash\mathrm{Im}(\tau_{u,v})
\ \longrightarrow \ 
\mathcal{P}_2 = \{(\rho,\nu)\}.
\end{equation}
The injectivity is obvious; for the surjectivity, we note that for any pair of partitions $(\rho,\nu)$, $\Gamma_\rho + (0,v)$ and $\Gamma_\nu + (u,0)$ always intersect as $u,v>0$, so there exists an inverse of $\iota_{u,v}$ given by cutting $\Gamma_\rho + (0,v)$ and $\Gamma_\nu + (u,0)$ at their last intersection point and cross-gluing the four pieces. For $(\rho,\nu) = \iota_{u,v}(\lambda,\mu,j)$, it is not hard to see that $|\rho|+|\nu| + uv \ = \ |\mu|+|\lambda|$ and
\begin{equation}
K_{\frac{c+\delta}{2}}(\lambda)+K_{\frac{c-\delta}{2}}(\mu) \equiv
K_{\frac{c+\delta}{2}}(\rho)+K_{\frac{c-\delta}{2}}(\nu)+\frac{uv}{2}(\epsilon_1+\epsilon_2+c) + v\epsilon_2(|\nu|-|\rho|) \quad \mathrm{mod} \ \delta-u\epsilon_1-v\epsilon_2.
\end{equation}
Summarizing the above argument, we see that the $v>0$ part in the RHS of (\ref{equation: proof lemma cov Phi K middle step 1}), which we denote by $W_{>0}$ for simplicity, eventually gives
\begin{equation}
\begin{split}
\Res_{\delta = u\epsilon_1+v\epsilon_2} W_{>0} = & - \sum_{\rho,\nu} z^{|\rho|+|\nu|+uv} \left( K_{\frac{c+\delta}{2}}(\rho)+K_{\frac{c-\delta}{2}}(\nu)+\frac{uv}{2}(\epsilon_1+\epsilon_2+c) + v\epsilon_2(|\nu|-|\rho|) \right) \\
= & -\frac{uv}{2}(c+\epsilon_1+\epsilon_2)\cdot z^{uv} P(z)^2 - \left((c+\epsilon_1+\epsilon_2)E_2(z)-(\epsilon_1+\epsilon_2)E_3(z)\right)z^{uv}P(z)^2
\end{split}
\end{equation}
for $u>0$, where the last step uses the norm of $K$ computed in Lemma \ref{lemma:mean-K}, and
\begin{equation}
\Res_{\delta = u\epsilon_1+v\epsilon_2} W_{>0} = 0
\end{equation}
for $u\leq 0$.

Similarly, for $v<0$, we define an injective map 
\begin{equation}
\tau_{u,v}:\ \left\{(\lambda,\mu,j)\ \middle|\ \begin{aligned}
j \text{ is a type-C intersection} \\
\text{jet of } \Gamma_\lambda \text{ and } \Gamma_\mu^{(u,v)}
\end{aligned}\right\} \ 
\longrightarrow \ 
\left\{(\lambda,\mu,j)\ \middle|\ \begin{aligned}
j \text{ is a type-D intersection} \\
\text{jet of } \Gamma_\lambda \text{ and } \Gamma_\mu^{(u,v)}
\end{aligned}\right\}
\end{equation}
by finding the intersection point immediately preceding $j$ and interchanging the piece of $\Gamma_\lambda$ and that of $\Gamma_\mu^{(u,v)}$ in between the two intersection points. The pole contribution to the RHS of (\ref{equation: proof lemma cov Phi K middle step 1}) from $(\lambda,\mu,j)$ cancels with that from $\tau_{u,v}(\lambda,\mu,j)$. $\tau_{u,v}$ is surjective unless $u<0$, in which case we define a bijection
\begin{equation}
\iota_{u,v}: \ 
\left\{(\lambda,\mu,j)\ \middle|\ \begin{aligned}
j \text{ is a type-D intersection} \\
\text{jet of } \Gamma_\lambda \text{ and } \Gamma_\mu^{(u,v)}
\end{aligned}\right\}\Big\backslash\mathrm{Im}(\tau_{u,v})
\ \longrightarrow \ 
\mathcal{P}_2 = \{(\rho,\nu)\}.
\end{equation}
again by interchanging the piece of $\Gamma_\lambda$ and that of $\Gamma_\mu^{(u,v)}$ after the last intersection point. In this way, we see that the $v<0$ part in the RHS of (\ref{equation: proof lemma cov Phi K middle step 1}), which we denote by $W_{<0}$ for simplicity, eventually gives
\begin{equation}
\Res_{\delta = u\epsilon_1+v\epsilon_2} W_{<0} = \frac{uv}{2}(c+\epsilon_1+\epsilon_2)\cdot z^{uv} P(z)^2 + \left((c+\epsilon_1+\epsilon_2)E_2(z)-(\epsilon_1+\epsilon_2)E_3(z)\right)z^{uv}P(z)^2
\end{equation}
for $u<0$, and 
\begin{equation}
\Res_{\delta = u\epsilon_1+v\epsilon_2} W_{<0} = 0
\end{equation}
for $u\geq 0$.

In the end, for $v=0$, we define a bijection
\begin{equation}
\tau_{u,v}:\ \left\{(\lambda,\mu,j)\ \middle|\ \begin{aligned}
j \text{ is a type-B intersection} \\
\text{jet of } \Gamma_\lambda \text{ and } \Gamma_\mu^{(u,v)}
\end{aligned}\right\} \ 
\longrightarrow \ 
\left\{(\lambda,\mu,j)\ \middle|\ \begin{aligned}
j \text{ is a type-C intersection} \\
\text{jet of } \Gamma_\lambda \text{ and } \Gamma_\mu^{(u,v)}
\end{aligned}\right\}
\end{equation}
by interchanging the pieces of $\Gamma_\lambda$ and $\Gamma_\mu^{(u,v)}$ preceding the intersection point. More precisely, put $j=(P,S_1,S_2)$. Then, we define $\tau_{u,v}(\lambda)$ as the partition such that $\Gamma_{\tau_{u,v}(\lambda)}$ is obtained from gluing the piece of $\Gamma_\lambda$ upper-left to $P$ with that of $\Gamma_\mu^{(u,v)}$ lower-right to $P$, $\tau_{u,v}(\mu)$ such that $\Gamma_{\tau_{u,v}(\mu)}^{(u,v)}$ is obtained from the other two pieces, and $$\tau_{u,v}(P,S_1,S_2) = (P,S_2,S_1).$$ $\tau_{u,v}$ is only well-defined when $v=0$ and thus both $\Gamma_\lambda$ and $\Gamma_\mu^{(u,v)}$ coincide with the $x$-axis at the lower-right infinity. Similar to the previous cases, the pole contribution to the RHS of (\ref{equation: proof lemma cov Phi K middle step 1}) from $(\lambda,\mu,j)$ and $\tau_{u,v}(\lambda,\mu,j)$ cancel. Hence, we see that the $v=0$ part in the RHS of (\ref{equation: proof lemma cov Phi K middle step 1}), which we denote by $W_{=0}$ for simplicity, gives for any $u$
\begin{equation}
\Res_{\delta = u\epsilon_1+v\epsilon_2} W_{=0} = 0.
\end{equation}

\subsection{The limit $\delta\rightarrow \infty$}
Substituting all results so far into (\ref{equation: proof lemma cov Phi K middle step 1}), we have that for any $u,v\in\mathbb{Z}$,
\begin{equation} \label{equation: proof lemma cov Phi K middle step 2}
\begin{split}
& \Res_{\delta = u\epsilon_1+v\epsilon_2}\langle \Phi_{\delta}(\lambda,\mu) \cdot (K_{\frac{c+\delta}{2}}(\lambda)+K_{\frac{c-\delta}{2}}(\mu))\rangle_2 \\
= & \Res_{\delta = u\epsilon_1+v\epsilon_2} \left(\frac{c+\epsilon_1+\epsilon_2}{2} \cdot z\frac{d\mathcal{C}_{\delta}(z)}{dz} + \left( (c+\epsilon_1+\epsilon_2)E_2(z) - (\epsilon_1+\epsilon_2)E_3(z)\right) \cdot \mathcal{C}_{\delta}(z)\right).
\end{split}
\end{equation}
Recall that by (\ref{eq:Phi-delta}), the only possible denominators of $\Phi_{\delta}(\lambda,\mu) \cdot (K_{\frac{c+\delta}{2}}(\lambda)+K_{\frac{c+\delta}{2}}(\mu))\rangle_2$ are of the form $\delta-u\epsilon_1-v\epsilon_2$, so (\ref{equation: proof lemma cov Phi K middle step 2}) implies
\begin{equation}
\begin{split}
& \frac{c+\epsilon_1+\epsilon_2}{2} \cdot z\frac{d\mathcal{C}_{\delta}(z)}{dz} + \left( (c+\epsilon_1+\epsilon_2)E_2(z) - (\epsilon_1+\epsilon_2)E_3(z)\right) \cdot \mathcal{C}_{\delta}(z) \\
= & \langle \Phi_{\delta}(\lambda,\mu) \cdot (K_{\frac{c+\delta}{2}}(\lambda)+K_{\frac{c+\delta}{2}}(\mu))\rangle_2 + \sum_{k\geq 0} z^k F_k(\delta;\epsilon_1,\epsilon_2,c),
\end{split}
\end{equation}
where $F_k$ are regular functions. For each fixed $k$, the $z^k$-coefficient of the LHS vanishes when $\delta\rightarrow\infty$ by definition of $\mathcal{C}_\delta$, and that of $\langle \Phi_{\delta}(\lambda,\mu) \cdot (K_{\frac{c+\delta}{2}}(\lambda)+K_{\frac{c+\delta}{2}}(\mu))\rangle_2$ vanishes when $\delta\rightarrow\infty$ as well since $\Phi_{\delta}(\lambda,\mu)$ gives $O(\delta^{-2})$ and $K_{\frac{c+\delta}{2}}(\lambda)+K_{\frac{c+\delta}{2}}(\mu)$ gives $O(\delta)$. Hence, $F_k$ vanishes when $\delta\rightarrow\infty$ as well, which forces it to be the constant zero. This completes the proof of the lemma.

\end{document}